\documentclass[abstracton, paper=a4, fontsize=11pt,DIV=12,bibliography=totoc,final]{scrartcl}

\pdfoutput=1
\usepackage{enumitem}
\usepackage{standalone}
\usepackage[utf8]{inputenc}
\usepackage[T1]{fontenc}
\usepackage{lmodern}
\usepackage[french,english]{babel}
\usepackage{csquotes}
\usepackage{amsthm, amssymb, amsmath, amsfonts, mathrsfs, dsfont, esint,textcomp}
\usepackage[colorlinks=true, pdfstartview=FitV, linkcolor=blue, citecolor=blue, urlcolor=blue,pagebackref=false]{hyperref}
\usepackage{mathtools}
\usepackage{tikz}

\usepackage{caption}
\usepackage{subcaption}

\usepackage{microtype}

\usepackage{url}
\usepackage[notcite,notref,color]{showkeys} 
\definecolor{labelkey}{gray}{.8}
\definecolor{refkey}{gray}{.8}

\usepackage[backend=biber,giveninits,style=alphabetic,maxnames=4,isbn=false, doi=false, eprint= true, url= false, sorting=ynt, sortcites = true]{biblatex}

\bibliography{bib.bib}

\definecolor{darkblue}{rgb}{0,0,0.7} 
\definecolor{darkred}{rgb}{0.9,0.1,0.1}
\definecolor{darkgreen}{rgb}{0,0.5,0}

\newcommand{\am}[1]{{{#1}}}
\newcommand{\amnew}[1]{{{#1}}}

\newcommand{\rh}[1]{{{#1}}}
\newcommand{\rhnew}[1]{{{#1}}}

\newcommand{\ml}[1]{{{#1}}}
\newcommand{\mlnew}[1]{{{#1}}}

\newtheorem{thm}{Theorem}[section]
\newtheorem{prop}[thm]{Proposition}

\newtheorem{lem}[thm]{Lemma}
\newtheorem{cor}[thm]{Corollary}
\theoremstyle{remark}
\newtheorem{rem}[thm]{Remark}
\theoremstyle{definition}
\newtheorem{defi}[thm]{Definition}

\renewcommand{\leq}{\leqslant}
\renewcommand{\geq}{\geqslant}

\renewcommand{\subset}{\subseteq}

\newcommand{\B}{\mathcal{B}}

\newcommand{\E}{\mathbb{E}}

\renewcommand{\L}{\mathcal{L}}

\newcommand{\N}{\mathbb{N}}
\newcommand\norm[1]{\left|\left| #1 \right|\right|}

\newcommand{\Ll}{\left}
\newcommand{\Rr}{\right}

\newcommand{\R}{\mathbb{R}}
\renewcommand{\S}{\mathbb{S}}

\newcommand{\Z}{\mathbb{Z}}
\newcommand{\mR}{\mathcal{R}}

\renewcommand{\H}{\mathcal{H}}

\renewcommand{\P}{\mathbb{P}}

\newcommand{\eps}{\varepsilon}
\renewcommand{\d}{{\mathrm{d}}}

\newcommand{\mP}{\mathcal{P}}
\newcommand{\mW}{\mathcal{W}}
\newcommand{\mS}{\mathcal{S}}

\newcommand{\dd}{\, d}

\newcommand{\Id}{\mathrm{Id}}

\renewcommand{\Re}{\mathrm{Re}}
\newcommand{\app}{\mathrm{app}}

\newcommand{\dmin}{d_{\min}}
\newcommand{\loc}{\mathrm{loc}}
\newcommand{\De}{\mathrm{De}}

\DeclareMathOperator{\dist}{dist}

\DeclareMathOperator{\dv}{div}
\DeclareMathOperator{\supp}{supp}
\DeclareMathOperator{\Sym}{Sym}
\DeclareMathOperator{\curl}{curl}

\DeclareMathOperator{\Tr}{Tr}

\numberwithin{equation}{section}

\mathtoolsset{showonlyrefs}

 \renewcommand{\arraystretch}{1.4}
\usepackage{microtype}

\DeclareMathOperator{\sym}{sym}

\def\XXint#1#2#3{{\setbox0=\hbox{$#1{#2#3}{\int}$ }
\vcenter{\hbox{$#2#3$ }}\kern-.6\wd0}}

 \usepackage{authblk}

\setkomafont{date}{\large}

\title{\LARGE Derivation of the viscoelastic stress in  Stokes flows induced by  non-spherical Brownian rigid particles through homogenization}{\normalfont\bfseries}

\author{Richard M. H\"ofer\thanks{hoefer@imj-prg.fr, Institut de Mathématiques de Jussieu -- Paris Rive Gauche, Université Paris Cité, 8 Place Aurélie Nemours, 75013 Paris, France}, Marta Leocata\thanks{mleocata@luiss.it, Luiss University, Viale Romania 32, 00197, Roma RM, Italy}, Amina Mecherbet\thanks{mecherbet@imj-prg.fr, Institut de Mathématiques de Jussieu -- Paris Rive Gauche, Université Paris Cité, 8 Place Aurélie Nemours, 75013 Paris, France}}

\begin{document}

\maketitle

\begin{abstract}
    We consider a microscopic model of $n$ identical axis-symmetric rigid Brownian particles suspended in a Stokes flow. We rigorously derive in the homogenization limit of many small particles a classical formula for the viscoelastic stress that appears in so-called Doi models which couple a Fokker-Planck equation to the Stokes equations. We consider both Deborah numbers of order $1$ and very small Deborah numbers.
    Our microscopic model contains several simplifications, most importantly, we neglect the time evolution of the particle centers as well as hydrodynamic interaction for the evolution of the particle orientations.
    
    \rh{The microscopic fluid velocity is modeled by the Stokes equations with given torques at the particles in terms of Stratonovitch noise. We give a meaning to this PDE in terms of an infinite dimensional Stratonovitch integral. This requires the analysis of the shape derivatives of the Stokes equations in perforated domains, which we accomplish by the method of reflections.
    }
\end{abstract}

{\small 
\noindent Keywords: Homogenization, Stokes flows, viscoelasticity, Doi model, Brownian particles

\smallskip

\noindent MSC: 76M50, 76D07, 35R60, 35Q70
}
\tableofcontents

 \section{Introduction} 
\label{sec:aim}
Suspensions of non-spherical rigid Brownian particles in viscous fluids are prototypes of viscoelastic fluids.  
It is classical (see e.g.  \cite{DoiEdwards88, KimKarilla13, Graham}) to model such complex fluids by an effective system consisting of a Fokker-Planck equation coupled to the Navier-Stokes equations which feature a viscoelastic stress $\sigma$ that depends on the particle density. 
In the absence of external forces, such a model for rod-like particles reads in dimensionless form
\begin{equation}
\label{eq:full.model} \def\arraystretch{1.8}
\left\{
\begin{array}{l}
\displaystyle \partial_t f + \dv (u f) + \dv_{\xi} (P_{\xi^\perp} \nabla_x u \xi f) 
	= \frac 1 {\De} \Delta_\xi f + \rh{\frac {\lambda_1} {\De}} \dv_x((\Id + \xi \otimes \xi) \nabla_x f), \\
\displaystyle	\Re (\partial_t u + (u \cdot \nabla) u ) -  \Delta u + \nabla p - \dv \sigma = 0, \qquad \dv u = 0\\
\displaystyle	\sigma = \sigma_v + \sigma_e = \rh{\lambda_2} \int_{\S^2} (D u : \xi \otimes \xi ) \xi \otimes \xi \dd \xi +  \rh{\frac {\lambda_3}{\De}} \int_{\S^2} (3 \xi \otimes \xi - \Id) f \dd \xi, \\
\displaystyle	u(0,\cdot) = u_0, \quad f(0,\cdot) = f_0.
\end{array}
\right.
\end{equation}
Here $u(t,x)$ and $p(t,x)$ are the fluid velocity and pressure, and $f(t,x,\xi)$ is the density of particles at time $ \rhnew{t\geq 0} $, position $ \rhnew{x \in \R^3} $  and orientation $\xi \in \S^2$. \rhnew{Moreover $P_{\xi^\perp}$ denotes the orthogonal projection in $\R^3$ to the subspace $\xi^\perp$.} 
\rh{\rhnew{Furthermore}, $\Re$ $\De, \lambda_1, \lambda_2, \lambda_3$ are dimensionless parameters, where $\Re$ is the Reynolds number and $\De$ is the Deborah number which is the ratio between the observation time scale and the diffusion time scale for the particle orientation. The parameter $\lambda_1 \ll 1$ depends only on the aspect ratio of the rod-like particles.
The parameters $\lambda_2, \lambda_3$ also depend  on the diluteness of the suspension.  It is usually argued that it is necessary for the validity of the Doi model that the particles can freely rotate which means $n \ell^3 \ll 1$ (so-called dilute regime in the terminology of \cite{GuazzelliHinch11}), where $n$ is the number density of the particles and $\ell$ the length of the rod-like particles. 
In this dilute regime $\lambda_2 \ll 1$, $\lambda_3 \ll 1$.}

\rh{The two parts of the viscoelastic stress $\sigma_v$ and $\sigma_e$ are sometimes referred to as the viscous part and the elastic part, respectively. 
The viscous part already occurs for suspensions of non-Brownian spherical particles for which it was first studied theoretically by Einstein \cite{Ein06} and later for ellipsoidal particles by Jeffery \cite{Jeffery22}. 
The viscous stress $\sigma_v$ as in \eqref{eq:full.model} can be obtained from Jeffery's computation as the leading order term for very elongated ellipsoids. 
Jeffery also computed the periodic orbits (Jeffery orbits) of such particles in constant shear flow by taking into account non-spherical inertialess particles are partly affected by pure straining motion of the fluid. In the case of rod-like particles the symmetric and skew-symmetric part of the gradient of the fluid velocity contribute  equally to the particle rotation (asymptotically for infinite aspect ratio) as reflected by the full gradient in the third term in the Fokker-Planck equation in \eqref{eq:full.model}.

The elastic part of the stress only arises for Brownian non-spherical particles. Theoretical studies on such elastic stresses go back to the 1940s, see e.g. \cite{Shima40,KuhnKuhn45,KirkwoodRiseman48} and also the later works \cite{HinchLeal71, HinchLeal72, Brenner74} and the references therein. 
To our knowledge the elastic stress $\sigma_e$ has so far not been obtained rigorously from a corresponding microscopic model. On the contrary the  scaling up has been described as \emph{mysterious} by Constantin and Masmoudi in \cite{constantin2008global}.
In the present paper, we provide such a mathematically rigorous derivation of  the elastic stress starting from a simplified microscopic system. We will treat both the cases of Deborah numbers $\De$ of order one, and very small Deborah numbers.}

\subsection{Previous mathematical results}
In the mathematical literature, the model \eqref{eq:full.model} is often called Doi model and has attracted a lot of attention over the last years.
Global well-posedness of such models have been studied in different contexts regarding the solution concept, the space dimension and the model for the fluid where in addition to  the imcompressible Navier-Stokes equations also the Stokes equations and the compressible Navier-Stokes equations have been considered. We refer to   \cite{constantin2005nonlinear, constantin2007regularity, constantin2008global,constantin2009holder, constantin2010global, OttoTzavaras08, BaeTrivisa12, BaeTrivisa13}. 
In \cite{HelzelTzavaras17, helzel2006multiscale} some further insights on rod-like suspensions are presented when effects of gravity are included. In \cite{saintillan2008instabilities}, \cite{saintillan2008instabilities2} a generalization of the Doi model for active particle is introduced. Existence of global weak entropy solution for this generalization of the Doi model is studied in \cite{chen2013global}
\rh{All the previously mentioned papers neglect the presence of the additional viscous stress $\sigma_v$, i.e., they pretend $\sigma = \sigma_e$.
Global well-posedness for the Doi model including the full viscoelastic stress $\sigma$ is treated in \cite{LionsMasmoudi07, ZhangZhang08, La19}.}

\medskip

 \rh{Regarding the derivation of the Doi model as a rigorous mean-field limit, there are only partial results so far. On the one hand, the viscous stress has been obtained in the quasistatic case. On the other hand, fully coupled systems have been obtained for non-Brownian spherical particles. In all these results the fluid is modeled by the stationary Stokes equations instead of the Navier-Stokes equations.

When the time evolution of the particles is neglected, the Stokes equations with a viscous stress corresponding to $\sigma_v$ above have recently been  obtained rigorously as homogenization limits for spherical particles in \cite{NiethammerSchubert19} and  subsequently for arbitrary shapes of the particles in \cite{HillairetWu19}. These results have been refined in \cite{ DuerinckxGloria20, DuerinckxGloria21, Duerinckx20,Gerard-VaretHillairet19, Gerard-VaretHoefer21, Gerard-VaretMecherbet20}.}

Dynamical models regarding the sedimentation of  inertialess non-Brownian spherical particles have been rigorously derived in \cite{JabinOtto04, Hofer18MeanField,Mecherbet18, Hofer&Schubert}

\medskip

\am{
\rh{Another prominent example of viscoelastic fluids are suspensions of flexible polymers, which are typically modelled by chain of monomers connected by springs, see e.g. \cite{BAH87, BCAH87, DoiEdwards88}.}
We refer to \cite{LL07,LL11} for a \rh{mathematical} introduction to the modelling of such fluids.
The corresponding \rh{viscoelastic} stress tensor $\sigma$ is given by the so-called Kramers expression which takes into account the polymer chain configurations and the force needed to extend or to compress the springs.
In particular, two models can be found in the literature regarding the force modelling the polymers extension: the potential associated to  Hooke's law which states that the force is linear to the length of the chain and the FENE (Finitely Extensible Nonlinear Elastic) pontential which takes into account the finite extensibility of the chain. A simplified model, the so-called dumbbell model, consists in assuming that the polymer is constituted of only two monomers, see \cite[Subsections 2.3 and 2.4]{LL11} for more details. We refer to \cite{JLL04,JLL02,JLLO06} and the references therein for well posedness, long time behaviour and numerical investigations of such models.
}


\medskip
Another popular model for the evolution of particles depending on their orientation are so-called Vicsek models.
In these models, the particle evolution is also described by a (kinetic) Fokker-Planck equation. The particles (which might be microswimmers, but also birds, fish or pedestrians) move in the direction of their orientation. Moreover, the particles change their orientation due to Brownian motion and interaction with each other. This interaction is given in terms of an interaction kernel instead of the interaction through the fluid in the Doi model. 
Viscek models have been studied mathematically including well-posedness, rigorous mean-field results and long-time behavior. One could refer for instance to  \cite{FigalliKangMorales15, GambaKang16, BriantDiezMerino21}.

\subsection{Heuristics}
\label{sec:heuristics}

Let us give a heuristic argument for the production of the elastic stress due to Brownian rod-like particles (general non-spherical particles are analogous). 
A similar argument can be found in \cite[p. 309]{DoiEdwards88}.

\begin{figure}
\centering
\begingroup
\tikzset{every picture/.style={scale=0.7}}%
\begin{subfigure}[b]{.3\textwidth}
\centering
\begin{tikzpicture}
\draw (-2,-0.5) node{};
\draw (3,7.5) node{};
\begin{scope} [shift={(0.5,0.5)}]

\draw (0,0) ellipse (0.2 and 0.1);
\draw (0,2) ellipse (0.2 and 0.1);
\draw (0.2,0)--(0.2,2);
\draw (-0.2,0)--(-0.2,2);
\draw[thick, red, bend left = 30, ->] (0,2) to (0.8,1.7);
\draw[thick, red, bend left = 30, ->] (0,0) to (-0.8,0.3);
\draw[thick, blue,  ->] (-1.3,2.5)--(1.3,2.5);
\draw[thick, blue,  ->] (1.3,-.5)--(-1.3,-.5);
\draw[thick, blue,  ->] (-1.5,2.3)--(-1.5,-.3);
\draw[thick, blue,  ->] (1.5,-0.3)--(1.5,2.3);

\draw[shift={(0,4)}] (0,0) ellipse (0.2 and 0.1);
\draw[shift={(0,4)}] (0,2) ellipse (0.2 and 0.1);
\draw[shift={(0,4)}] (0.2,0)--(0.2,2);
\draw[shift={(0,4)}] (-0.2,0)--(-0.2,2);
\draw[shift={(0,4)}, thick, red, bend left = 30, ->] (0,2) to (0.8,1.7);
\draw[shift={(0,4)}, thick, red, bend left = 30, ->] (0,0) to (-0.8,0.3);
\draw[shift={(0,4)}, thick, blue,  ->] (-1.3,2.5)--(1.3,2.5);
\draw[shift={(0,4)}, thick, blue,  ->] (1.3,-.5)--(-1.3,-.5);
\draw[shift={(0,4)}, thick, blue,  ->] (1.5,2.3)--(1.5,-.3);
\draw[shift={(0,4)}, thick, blue,  ->] (-1.5,-0.3)--(-1.5,2.3);
\end{scope}

;\end{tikzpicture}
\caption{Rotations due to fluid\\velocity gradients}
\label{fig:rotating.rod.gradient}
\end{subfigure}
\begin{subfigure}[b]{.3\textwidth}
\centering
\begin{tikzpicture}
\draw (-2,-0.5) node{};
\draw (3,7.5) node{};
\begin{scope} [shift={(0.0,0.5)},rotate=-20]
\draw (0,0) ellipse (0.2 and 0.1);
\draw (0,2) ellipse (0.2 and 0.1);
\draw (0.2,0)--(0.2,2);
\draw (-0.2,0)--(-0.2,2);
\draw[thick, red, bend left = 30, ->] (0,2) to (0.8,1.7);
\draw[thick, red, bend left = 30, ->] (0,0) to (-0.8,0.3);
\draw[thick, blue,  ->] (-1.3,2.5)--(1.3,2.5);
\draw[thick, blue,  ->] (1.3,-.5)--(-1.3,-.5);
\draw[thick, blue,  ->] (-1.5,2.3)--(-1.5,-.3);
\draw[thick, blue,  ->] (1.5,-0.3)--(1.5,2.3);
\end{scope}

\begin{scope} [shift={(0.75,5)},rotate=20] 
\draw (0,0) ellipse (0.2 and 0.1);
\draw (0,2) ellipse (0.2 and 0.1);
\draw (0.2,0)--(0.2,2);
\draw (-0.2,0)--(-0.2,2);
\draw[ thick, red, bend right = 30, ->] (0,2) to (-0.8,1.7);
\draw[ thick, red, bend right = 30, ->] (0,0) to (0.8,0.3);
\draw[ thick, blue,  ->] (1.3,2.5)--(-1.3,2.5);
\draw[ thick, blue,  ->] (-1.3,-.5)--(1.3,-.5);
\draw[ thick, blue,  ->] (1.5,2.3)--(1.5,-.3);
\draw[ thick, blue,  ->] (-1.5,-.3)--(-1.5,2.3);
\end{scope}

;\end{tikzpicture}
\caption{Fluid velocity caused by\\rotations}
\label{fig:rotating.rod.rotated}
\end{subfigure}
\begin{subfigure}[b]{.3\textwidth}
\centering
\begin{tikzpicture}

\draw (-2,-0.5) node{};
\draw (3,7.5) node{};
\begin{scope} [shift={(0.5,2.5)}]

\draw (0,0) ellipse (0.2 and 0.1);
\draw (0,2) ellipse (0.2 and 0.1);
\draw (0.2,0)--(0.2,2);
\draw (-0.2,0)--(-0.2,2);
\draw[thick, red, bend left = 30, <->] (-0.8,1.8) to (0.8,1.8);
\draw[thick, red, bend left = 30, <->] (0.8,0.2) to (-0.8,0.2);
\draw[ thick, blue,  ->] (0,-1)--(0,-0.2);
\draw[ thick, blue,  ->] (0,3)--(0,2.2);
\draw[ thick, blue,  ->] (1,1)--(1.8,1);
\draw[ thick, blue,  ->] (-1,1)--(-1.8,1);

\end{scope}

;\end{tikzpicture}
\caption{Average fluid velocity due to\\Brownian motion}
\label{fig:rotating.rod.stress}
\end{subfigure}
\caption{Heuristic explanation of the viscoelastic stress arising from rotational Brownian motion.  The blue arrows represent the motion of the fluid whereas the red arrows represent the motion of the rod.}
\endgroup
\label{fig:rotating.rod}
\end{figure} 
Consider a rod orientated in some direction as in Figure \ref{fig:rotating.rod.gradient}. Since the rod is inertialess, it follows the fluid flow and consequently rotates if the surrounding fluid is rotating as well. However, as we see in the lower figure, it also rotates under purely straining motion of the fluid. This is the reason why the full gradient of $u$ appears in the third term of the Fokker-Planck equation in \eqref{eq:full.model} and not only its skewsymmetric part.

This is directly related (though the symmetry of the resistance tensor, see Subsection \ref{sec:resistance}) to the fact that a rigid rod-like particle subject to a torque needs to exert a stresslet on the fluid in order to maintain its shape. Therefore,  Brownian torques  on the particles arising from thermal noise lead to corresponding Brownian stresslets. It can then be argued that there is an average net stresslet at each of the particles which gives rise to the elastic stress $\sigma_e$,  see Figure \ref{fig:rotating.rod.stress} where the generated stresslet is illustrated through the resulting fluid velocity.
There is one additional subtlety in the above argument: By linearity of the Stokes equations, the instantaneous stresslet produced by a torque which corresponds to a rotation to the left (of a vertically oriented rod)
is exactly the negative of the effect of  a torque on the same rod corresponding to a rotation to the right. Due to the random nature, such torques occur with equal probability and therefore seem to cancel out.

However, since the random torques produce a Brownian motion of the particle orientation one has to take into account quadratic effects: One has to consider the stresslets produced by such torques after the rod has already started to rotate as visualized in Figure \ref{fig:rotating.rod.rotated}. Summing these contributions leads  indeed to the effective stresslet as in Figure \ref{fig:rotating.rod.stress}.

We summarize that there are three key ingredients that cause the elastic part of the stress: 1) The Brownian torque on the fluid; 2) the corresponding Brownian motion of the particle orientation;  3) the particle anisotropy.

We refer to Subsection \ref{subsec:proof_strategy} for a more formalized version of this heuristic argument.


\subsection{Formal statement of the main results}

The rigorous derivation of the complete system \eqref{eq:full.model} from a microscopic description seems to be completely out of reach for the moment due to the highly singular interaction of the particles. 

Instead, we consider a simplified microscopic model keeping those aspects that produce the viscoelastic stress:
First, we model the fluid by the stationary Stokes equations with no-slip conditions and balance laws on each particle.
Second, we freeze the time evolution of the particle centers and assume that there is no exchange of net force between the particles and the fluid. Third, we model the evolution of the particle orientation as if they were all alone in infinite fluid. Correspondingly, the Brownian torques at each particle are independent of each other.
We also include an external torque acting on the particles which could be related to an external fluid flow, a magnetic force or chemotaxis. 
We mainly consider such a torque in order to have a nontrivial particle distribution for very small Deborah numbers $\De \to 0$. For simplicity we do not include the effect of these torques on the fluid.

After non-dimensionalization, this leads to the system
\begin{subequations} \label{eq:micro.T_D.intro}
\begin{equation} \label{eq:Stokes.micro.T_D.intro}
\left \{
\begin{array}{rcl}
\displaystyle - \Delta u_n+ \nabla p_n&=& 0 \quad  \text{ in }\mathbb{R}^3 \setminus \bigcup_{i=1}^n \B_i(t),\\
\displaystyle \dv u_n&=&0 \quad  \text{ in }\mathbb{R}^3 \setminus \bigcup_{i=1}^n \B_i(t),\\
\displaystyle D u_n&=&0 \quad  \text{ in } \bigcup_{i=1}^n  \B_i(t), \\
\displaystyle \int_{\partial \B_i(t)} \Sigma(u_n,p_n) \nu & = &0 \quad  \text{ for all } 1 \leq i \leq n,\\
\displaystyle \int_{\partial \B_i(t)} [\Sigma(u_n,p_n) \nu ] \times (x-x_i) &=& \frac {\phi_n} n \sqrt{2 \gamma_{rot}\mR_2(\xi_i(t))} \circ \dot B_i(t) \quad  \text{ for all } 1 \leq i \leq n,
\end{array}
\right.
\end{equation}
\begin{align} \label{eq:Particles.T_D.intro}
\left \{
\begin{array}{rcl}
\displaystyle	 \dd \xi_i(t) &=& \sqrt{2} \xi_i(t) \times    \circ \dd B_i(t) +  P_{\xi_i^\perp} h(t,\xi_i(t),x_i) \dd t ,  \\
\xi_i(0) &=& \xi_{i,0}.
\end{array}
\right.
\end{align}
\end{subequations}
Here, the particles $\B_i(t)$, $1 \leq i \leq n$, are obtained from an arbitrary axisymmetric reference particle (not necessarily rod-like) by translation, rotation and rescaling with a factor $r$ which only depends on the number of particles $n$.
Moreover,  $\phi_n = n r^3$ is the volume fraction of the particles, $x_i \in \R^3$, $\xi_i \in \S^2$ the particle centers and orientations, $\mR_2(\xi_i)$ and $\gamma_{rot}$ are related to the particle resistance to rotation,  $D u_n$ is the symmetric gradient of $u_n$, 
$\Sigma(u_n,p_n)$ is the fluid stress, and $\nu$ the outer unit normal. Moreover, $B_i$ are Brownian motions, $\dot B_i$ corresponding white noise and $\circ$ indicates a product in the Stratonovitch sense. Finally $h$ is a given function related to an external torque.

\rh{We emphasize that the equation for the fluid \eqref{eq:Stokes.micro.T_D.intro}  does not involve the time derivative $\partial_t u$ as we neglect the fluid inertia. Nevertheless, the fluid velocity $u = u(t,x)$}  implicitly depends on time trough the particle orientations as well as the random torques. The latter also prevents an interpretation of this equation pointwise in time. For the precise notion of solutions $u$ to \eqref{eq:Stokes.micro.T_D.intro} we refer to Subsection \ref{sec:well-posedness}.

For details about the modeling we refer to Section \ref{section2} where we obtain this system (after non-dimensionalization) as a simplified version of a more physically accurate microscopic system. We refer to Subsection \ref{subsec:discussion} for a discussion on the difficulties to treat such a more accurate system. 
A benefit of this simplified system is that it reveals very clearly that the elastic part of the stress arises because of the interplay of the rotational Brownian motion, the Brownian torque on the fluid and the particle anisotropy, as explained in the previous subsection.
 More precisely, the elastic stress $\sigma_e$ arises as a mean-field/Law-of-Large-Numbers term from the expectation of the stresslets produced at each particle through the random torque, whose time integrals read
\begin{align} \label{eq:exp.stresslet}
    \E \left[ \int_0^t \mS(\xi_i(s)) \sqrt{ 2\gamma_{rot}\mR_2(\xi_i(s))} \circ \dd B_i(s) \right] .
\end{align}
Here, $\mS(\xi_i(s)) \colon \R^3 \to \R^{3 \times 3}$ is the linear map that maps a given torque to the stresslet that a particle with orientation $\xi_i(s)$ exerts on the fluid when the particle is alone in an unperturbed fluid (see Subsection \ref{sec:resistance}). Since $\xi_i(s)$ itself solves the SDE \eqref{eq:Particles.T_D.intro} involving $B_i$, the It\^o-Stratonovitch conversion formula (see Definition \ref{def:Strato}) yields that the expectation in \eqref{eq:exp.stresslet} is non-vanishing but instead (see Lemma \ref{lem:Expectations})
\begin{align}
    \E \left[ \int_0^t \mS(\xi_i(s)) \sqrt{ 2\gamma_{rot}\mR_2(\xi_i(s))} \circ \dd B_i(s) \right] &= \gamma_E \E\int_0^t[3 \xi_i(s) \otimes \xi_i(s) - \Id] \dd s.
\end{align}
In this view, the deterministic elastic stress $\sigma_e$ and the deterministic orientational diffusion $\Delta_\xi f$ arise in a somewhat similar way due to the quadratic variation of Brownian motion.

We emphasize that the ratio between the viscoelastic time scale and the rotational diffusion time scale is of order of the volume fraction of the particles $\phi_n$ which corresponds to the parameter $\lambda_3$ in \eqref{eq:full.model}. Our methods restrict us to assume $\phi_n \to 0$ (see Subsection \ref{sec:assumptions} for the precise assumptions). 
However, we emphasize that we consider in this paper two fundamentally different scalings
\begin{enumerate}[label=(\roman*)., ref=(\roman*)]
    \item The case of Deborah numbers of order $1$: Here we consider the observation timescale to be comparable to the diffusion time scale; hence the rescaled viscoelastic stress term is small whereas the rescaled diffusion coefficient is of order one ($\gamma_{rot}$ is of order one). This case corresponds to the system \eqref{eq:micro.T_D.intro} above.
        \item Very small Deborah numbers, where we consider the observation timescale to be comparable to the viscoelastic time scale; hence the rescaled viscoelastic stress term is of order one whereas the rescaled diffusion coefficient is large, see \eqref{eq:Stokes.micro.T_u}, \eqref{eq:Particles.T_u} for the corresponding microscopic system.
\end{enumerate}

In the first case, we obtain that the empirical measure of the particles converges to the solution of the instationary Fokker-Planck equation
\begin{equation}\label{eq:Fokker-Planck.instationary}
\left\{ 
\begin{array}{rcl}
	\partial_t f + \dv_\xi (P_{\xi^\perp} h f - \nabla_\xi f ) &=& 0, \\
	f(0,\cdot) &=& f_0.
	\end{array}
	\right.
\end{equation}
Since in this case, the total viscoelastic stress, which is the only source term for the fluid equation, is of order $\phi_n$, we consider rescaled fluid velocities $\phi_n^{-1} u_n$. We show that this rescaled sequence $\phi_n^{-1} u_n$ converges to
the Stokes equations with an additional viscoelastic stress, namely
\begin{align} {\def\arraystretch{1.6}} \label{eq:viscoelastic}
\left\{  
\begin{array}{l}
\displaystyle	-\Delta u + \nabla p = \dv \sigma, \qquad \dv u = 0,\\
\displaystyle	\sigma(t,x) = \int \gamma_E (3 \xi \otimes \xi - \Id) f(t,x,\xi) \dd \xi.
	\end{array}
	\right.
\end{align}
The parameter $\gamma_E \in \R$ depends only on the shape of the reference particle.
We emphasize that the viscoelastic stress $\sigma$ here corresponds to the elastic part $\sigma_e$ in \eqref{eq:full.model}. 
The viscous part does not appear in this case because its effect is of order $\phi_n^2$. It is classical that the viscous stress produced by the particles is proportional to their volume fraction to leading order. In our case however, since the elastic stress as the only source term is also of order $\phi_n$, the total effect of the viscous stress is indeed quadratic in the volume fraction.

In the second case, which corresponds to very small Deborah numbers $\De \to 0$, we obtain in the limit the (quasi)-stationary Fokker-Planck equation for the particles
\begin{align} \label{eq:Fokker-Planck.stationary}
\left\{ \def\arraystretch{1.6}
\begin{array}{rcl}
		\dv_\xi (P_{\xi^\perp} h f - \nabla_\xi f ) &=& 0, \\
	\displaystyle	\int f(t,\cdot) \dd \xi &=& \int f_0 \dd \xi.
	\end{array}
	\right.
\end{align}
In this case, the fluid velocity $u_n$ itself converges to the solution to \eqref{eq:viscoelastic}. 

The Fokker-Planck equation \eqref{eq:Fokker-Planck.stationary}  depends on time only through  $h$. \rhnew{Since the equation is elliptic, the initial configuration $f_0$ only enters as a constraint on the spacial density $\int f(t,\cdot) \dd \xi$ which is constant in time because the particles do not move in space.
However, the solution $f$ to \eqref{eq:Fokker-Planck.stationary} does in general not satisfy $f(0,\cdot) = f_0$. This discrepancy in the initial data arises from the fast diffusion for very small Deborah numbers which creates an initial layer for the solution of the microscopic particle density. This initial layer can be related to the instationary Fokker Planck equation \eqref{eq:Fokker-Planck.instationary} in the sense that $f(0,\cdot)$ is given as long-time limit $\lim_{t\to \infty} \tilde f$ for the solution $\tilde f$ to \eqref{eq:Fokker-Planck.instationary} with $h$ replaced by $\tilde h(t,\cdot) = h(0,\cdot)$.}

\medskip
Elliptic equations with temporal noise like \eqref{eq:Stokes.micro.T_D.intro} seem to have been rarely studied in the mathematical literature so far. A similar model is analyzed numerically in \cite{Bao18}, though. Moreover, such equations seem to be (implicitly) underlying the derivation of the Langevin equations for Brownian particles suspended in a fluid (see e.g. \cite{Roux92},  \cite[Section 5.4.2]{GuazzelliMorris12} and \cite[Section 6.4]{Graham}).
We emphasize that already making sense of the equations for the fluid velocity $u_n$ in \eqref{eq:micro.T_D.intro} is non-trivial
due to the Stratonovitch white noise in the torque that acts as a boundary condition for the Stokes equations. 
We overcome this issue by using the linearity of the problem in order to define $u_n$  as the distributional derivative of an Hilbert space  valued Stratonovitch integral, see Subsection \ref{th:well-posed}.

For the regularity in time of $u_n$ we obtain the optimal regularity $H^{-s}$, $s<\frac{1}{2}$, which corresponds to the regularity of white noise.
For homogenization problems of the Stokes equations in perforated domains where the stresslets at the particles produce a nontrivial term in the limit, it is classical to obtains $L^p$-convergence in space of the fluid velocity for $p < 3/2$. Due to issues with Banach valued stochastic integrals, we will work in negative spaces $H^{-s}$, $s<\frac{1}{2}$, with respect to time and space.

Moreover, the Stratonovitch nature of the white noise makes it necessary to consider shape derivatives of the solution to the Stokes equations in perforated domains with prescribed torques.
Such estimates and approximations for this solution will be obtained through the method of reflections that has already been used in several related works (e.g. \cite{Hofer18MeanField, Mecherbet18, NiethammerSchubert19}) but will here be used for the first time to estimate shape derivatives.

\subsection{Organization of the rest of the paper}
The rest of the paper is organized as follows 
\begin{itemize}
    \item 
    \rh{Section 2 is devoted to the modeling of the microscopic system that eventually leads to \eqref{eq:micro.T_D.intro}. After specifying the assumptions on the particle shape and recalling properties of the grand-resistance tensor for a single (axisymmetric) particle in Subsections \ref{sec:particles} and \ref{sec:resistance}, we introduce a full microscopic model of inertialess Brownian particles in a Stokes flow in Subsection \ref{sec:dynamics}. Subsequent simplifications and nondimensionalization in Subsections \ref{sec:simplification} and \ref{sec:nondim} then leads to \eqref{eq:micro.T_D.intro}.}
    \item In Section \ref{section3} we first present the main assumptions on the initial particle configuration and specify the notion of solutions to \eqref{eq:micro.T_D.intro}. The  main convergence results  are then stated in Subsection \ref{sec:convergence}. In Subsection \ref{subsec:notatation} we introduce the main notations used in this paper. We summarize the key steps of the proof in Subsection \ref{subsec:proof_strategy}  and we discuss about the limitations and possible generalizations of our approach in Subsection \ref{subsec:discussion}.
    \item Sections \ref{sec:L_n}--\ref{section6} are devoted to the proof of the main results. For a more detailed outline of these sections and the appendices, we refer to Subsection \ref{subsec:proof_strategy}.
\end{itemize}

\section{The microscopic model}\label{section2}

\subsection{The particles} \label{sec:particles}
We consider $n$ identical particles $\B_i$ given by scaling, rotation and translation of a reference particle $\B$.
We assume that $0 \in \B \subset B(0,1)$ is a smooth compact set with rotational symmetry, i.e.
\begin{align}\label{rod_reference}
	R \B = \B \quad \text{for  all } R \in SO(3) \text{ with } R e_3 =e_3,
\end{align}
where $e_3 = (0,0,1) \in \R^3$.

We then consider $n$ identical particles 
 $\B_i = x_i + r R(\xi_i) \B$, where $r>0$ a scaling factor, $x_i \in \R^3$ is the position and $\xi_i \in \S^2$ 
the orientation of the $i$-th particle. The rotation matrix $R(\xi_i) \in SO(3)$ is chosen such that $R(\xi_i) e_3 = \xi_i$. Note that this constraint does not characterize $R(\xi_i)$ uniquely, but due to the rotational symmetry of $\B$, the choice of $R(\xi_i)$ does not affect $\B_i$.

\subsection{The resistance tensor} \label{sec:resistance}

For the setup of the model, we need to recall
the notion of the (Stokes) resistance tensor.
A detailed discussion on this topic can be found for example in \cite[Chapter 5]{KimKarilla13}.

In the following $\Sym_0(3)$ denotes the symmetric traceless matrices in $\R^{3 \times 3}$.
For an arbitrary (smooth) bounded domain $A \subset \R^3$,
the (grand) resistance tensor 
$\mR_A \in \R^{6 \times 6}$ relates the translational and angular velocities $V, \omega \in \R^3$ and the 
rate of strain  $E \in \Sym_0(3)$ of a particle in a quiescent Stokes flow to the force, torque and stresslet exerted on the fluid. More precisely, let $x_0 \in \R^3$ be a fixed reference point, and 
consider the solution $(v,p) \in (\dot H^1(\R^3),L^2(\R^3))$ (where $\dot{H}^1(\R^3)$ stands for the standard homogeneous Sobolev space) to the Stokes equations
with some viscosity $\mu > 0$
\begin{equation} \label{eq:resistance.problem}
\left\{
\begin{array}{ll}	- \mu \Delta v + \nabla p = 0, \quad \dv v = 0 &\quad  \text{in } \R^3 \setminus A, \\
	v = V + \omega \times (x - x_0) + E(x-x_0) & \quad \text{in } A.
	\end{array} \right.
\end{equation}
Then, the force and torque and stresslet exerted on the fluid by the particle are given as
\begin{align}
	F &= -\int_{\partial A} \Sigma[v,p] \nu, \\
	T &= -\int_{\partial A} \Sigma[v,p] \nu \times (x-x_0), \\
	S &= - \sym_0\left(\int_{\partial A} \Sigma[v,p] \nu \otimes (x-x_0) \right), 
\end{align}
with $\nu$ the outer normal at $\partial A$ and $\Sigma[v,p] = 2 \mu Dv - p \Id$ the fluid stress where $D v = \sym(\nabla v)$ and $\sym B$ and $\sym_0 B$ denote the symmetric and symmetric traceless part of a matrix $B \in \R^{d \times d}$, respectively. 


Then, the resistance tensor is defined through\footnote{By linearity, one can replace the velocities $V,\omega,E$ by the corresponding quantities relative to prescribed velocities of the fluid at infinity, which we set equal to zero by imposing $v \in \dot H^1(\R^3)$. Considering nonzero velocities at infinity, $E$ might be nonzero even for a rigid particle.}
\begin{equation} \label{def:R}
\begin{pmatrix}
  F \\
  T \\
  S
\end{pmatrix}= 
\mu \mR_A \begin{pmatrix}
  V \\
  \omega \\
  E 
\end{pmatrix}.
\end{equation}

By linearity of the Stokes equations, $\mR_A$ is well defined. Moreover, by some integration by parts, $\mR_A$ is symmetric and positive definite.
We denote
\begin{align*}
	\bar \mR := \mathcal R_{\mathcal B} :=  \begin{pmatrix}
  \bar \mR_1 & \bar \mR_{12} & \bar \mR_{13} \\  \bar \mR_{12}^T &  \bar \mR_2 & \bar \mR_{23} \\
  \bar \mR_{13}^T & \bar \mR^T_{23} & \bar \mR_3
\end{pmatrix},
\end{align*}
and find due to rotational symmetry
\begin{align}
	\bar \mR_1 &= \gamma_\perp (\Id - e_3 \otimes e_3) + \gamma_\parallel e_3 \otimes e_3,  \label{eq:translational.resistence.reference}\\
	 \bar \mR_2 &= \gamma_{rot} (\Id - e_3 \otimes e_3) + \gamma_{rot,\parallel} e_3 \otimes e_3, \label{eq:angular.resistence.reference} \\
	\bar \mR_{23}^T \omega &= \gamma_E \Sym\left( (\omega \times e_3) \otimes e_3\right ),
\end{align}
for some $\gamma_\perp,\gamma_\parallel, \gamma_{rot}, \gamma_{rot,\parallel} > 0$, $\gamma_E \in \R$. Formula for these quantities in the case of spheroids can be found in \cite[Subsection 3.3]{KimKarilla13}). In particular $\gamma_E \neq 0$ for all spheroids except for spheres.
{(We omit corresponding formula for $\bar \mR_3$, $\bar \mR_{12}$, $\bar \mR_{13}$ 
(see \cite[Subsection 5.5]{KimKarilla13}). In fact, $\bar \mR_{12}=0 = \bar
\mR_{13}$ according to \cite{Graham} by choosing $x_0$ as the so-called center of hydrodynamic reaction.)} 




By straightforward transformation arguments the resistance tensor of $\B_i$ with respect to $x_i$ is given by
\begin{equation} \label{def:R_B_i}
\mR_{\B_i} =
\begin{pmatrix}
r \mR_1 (\xi_i) & r^2 \mR_{12} (\xi_i) & r^2 \mR_{13} (\xi_i) \\
 r^2 \mR_{12}^T (\xi_i) &  r^3 \mR_2(\xi_i) & r^3 \mR_{23}(\xi_i) \\
 r^2 \mR_{13}^T (\xi_i) &  r^3 \mR_{23}^T(\xi_i) &  r^3 \mathcal{R}_{3}(\xi_i) 
\end{pmatrix},
\end{equation}
where
\begin{align} 
 	\mR_1(\xi_i) &= \gamma_\perp (\Id - \xi_i \otimes \xi_i) + \gamma_\parallel \xi_i \otimes \xi_i, \label{eq:R_1}\\
 	 \mR_2(\xi_i) &=  \gamma_{rot} (\Id - \xi_i \otimes \xi_i) + \gamma_{rot,\parallel} \xi_i \otimes \xi_i, \label{eq:R_2} \\
 	 \mR_{23}^T(\xi_i) \omega &= \gamma_E \Sym\left( (\omega \times \xi_i) \otimes \xi_i\right ).  \label{eq:R_23}
\end{align}

For an inertialess rigid particle, one is interested in the problem to determine the particle velocities $V,\omega$ as well as the stresslet $S$ when  force and torque $F,T$ are given, as well as the rate of strain $E$ (which corresponds to the rate of strain of the fluid far from  the particle). This is known as the mobility problem. 
We consider here only the case $E=0, \rhnew{F=0}$. In this case, we have
\begin{align*}
	\omega &=  \mu^{-1} r^{-3} \mR_2^{-1} T, \\
	S &= \mu r^3 (\mR_{23})^T \omega = \mR_{23}^T \mR_2^{-1} T 
\end{align*}

The mapping $\mR_{23}^T \mR_2^{-1}$ which relates the torque to the stresslet will play an important role in the analysis of the viscoelastic stress, hence we introduce the following operator 
\begin{equation} \label{def:mS_i}
\begin{array}{lccl} 
\mS :& \mathbb{S}^2 &\to& \L(\mathbb{R}^3, \Sym_0 (\mathbb{R}^3)),\\
& \xi & \mapsto & \mR_{23}^T(\xi) \mR_2^{-1}(\xi).
  \end{array}
\end{equation}
Note that $\mS$ is smooth. In particular, there exists a constant $C>0$ such that 
\begin{equation}\label{reg_mS_i}
    \|\mS \|_{{C}^1(\mathbb{S}^2;\mathcal{L}(\mathbb{R}^3, \Sym_0 (\mathbb{R}^3))} \leq C .
\end{equation}

\subsection{The dynamics} \label{sec:dynamics}

%
%

 We assume that the fluid satisfies the Stokes equations with no-slip condition at the particles:
\begin{equation} \label{eq:u_n}
\left \{
\begin{array}{rcll}
-\mu \Delta u_n+ \nabla p_n&=& 0 &\quad \text{ in }\mathbb{R}^3 \setminus \bigcup_{i=1}^n \B_i,\\
\dv u_n&=&0 &\quad\text{ in } \mathbb{R}^3 \setminus \bigcup_{i=1}^n \B_i,\\
u_n&=& v_i+\omega_i \times (x-x_i) &\quad \text{ in } \B_i.
\end{array}
\right.
\end{equation}
Here, $v_i, \omega_i \in \R^3$ are the translational and angular velocities on $\B_i$.
 Neglecting the particle inertia, the velocities $v_i, \omega_i$ are not given but they are determined through the following conditions, prescribing the total force and torque acting on each particle:
\begin{align} \label{eq:F_i.T_i}
\int_{\partial B_i} \Sigma(u_n,p_n) \nu &=F_i ,\\
  \int_{\partial B_i} [\Sigma(u_n,p_n) \nu ] \times (x-x_i)&=T_i. 
\end{align}
Since the particles are inertialess, these forces and torques balance the forces and torques acting on the particles
which are the sum of external forces and torques and of the random forces and torques acting on each particle due to collisions with fluid particles, $F_i = F_i^E + F_i^B$ and $T_i = T_i^E + T_i^B$.

\rh{
According to the fluctuation-dissipation theorem, the random forces and torques are given by 
\begin{align} \label{eq:full.Stokes.Einstein}
    	(F^B,T^B) = \sqrt{2 k_B \Theta \mu\mathscr R_n} \mlnew{\circ \dot B},
\end{align} 
see e.g. \cite{Roux92}. Here,
$k_B$ is the Boltzmann constant, $\Theta$ the absolute temperature, $B$ is a  $6n$-dimensional Brownian motion and $F^B, T^B \in \R^{3n}$ are the vectors containing all the forces  and torques $F_i^B, T_i^B$. Moreover, $\mathscr R_n \in \R^{6n\times6n}$ is the resistance matrix  of all the particles (excluding, stresslet/strain). More precisely, similar as in \eqref{def:R}, $\mathscr R_n$ relates given velocities $V_i,\omega_i$ at all particles to  forces and torques $F_i,T_i$ by solving the corresponding $n$ particle problem instead of \eqref{eq:resistance.problem}. In particular $\mathscr R_n$ depends on the positions and orientations of all the $n$ particles.
}



The fluid equations are complemented by the equations of motion for the particles
\begin{align} 
	\dot{x}_i &= v_i \label{eq:velocity},\\
	\dot{\xi}_i &= \omega_i \times \xi_i. \label{eq:angularVelocity}
\end{align}

 We emphasize once more that $u=u(t,x)$ implicitly depends on time through the time dependent fluid domain $ \mathbb{R}^3 \setminus \bigcup_{i=1}^n \B_i$ as well as the forces and torques $F_i,T_i$. Moreover, the notion of solutions to such fluid equations will be specified in Subsection \ref{sec:well-posedness}.

\subsection{Simplification of the model} \label{sec:simplification}


As outlined in the introduction, deriving the Doi model from \eqref{eq:u_n}--\eqref{eq:angularVelocity} seems presently out of reach. We now detail the simplifications that lead from \eqref{eq:u_n}--\eqref{eq:angularVelocity} to the model \eqref{eq:Stokes.micro.T_D.intro}--\eqref{eq:Particles.T_D.intro}.


Instead of \eqref{eq:velocity} and \eqref{eq:angularVelocity}, we fix the particle centers, and set the forces $F_i$ equal to $0$,
\begin{align}
		\dot{x}_i &= 0 \label{eq:velocity0}, \\
		F_i &= 0.
\end{align}
Moreover, instead of the equation of motion for the particle orientation
\eqref{eq:angularVelocity}, we assume
\begin{align}
		\dd \xi_i(t) &=  \xi_i(t) \times \sqrt{2 k_B \Theta \mu^{-1} r^{-3} \mR_{2}^{-1}} \circ \dd B_i(t) + P_{\xi_i^\perp} h(t,\xi_i(t),x_i) \dd t \label{eq:acceleration2},
\end{align}
where $\mR_2$ is as in \eqref{def:R_B_i}.
The first term on the right-hand side above corresponds to an angular velocity $\omega_i$ caused by random collisions with fluid particles as if the particle $\B_i$ was alone in the fluid, i.e. $\omega_i =  \mu^{-1} r^{-3} \mR_{2}^{-1} \rhnew{T_i^B}$ \rhnew{with $T_i^B = \sqrt{2 k_B \Theta \mu r^{3} \mR_2(\xi_i)}  \mlnew{\circ \dot B_i}$}. This is known as hydrodynamic decoupling
and is at least formally justified for small particle volume fraction.
Here, and in the following, \mlnew{$B_i$ is brownian motion in $\R^3$}.
The additional term in \eqref{eq:acceleration2} containing $h$ could be understood as arising from an external torque for example associated to an external fluid flow, a magnetic field or chemotaxis. 
By an action reaction principle, a corresponding torque should act on the fluid as well, but we will omit this for the sake of simplicity.
More precisely, regarding the fluid equations, we consider the random torques $T_i^B$ as the only
torques acting on the fluid.

Then, the simplified model at one glance is given by
\begin{subequations}  \label{eq:reduced}
\begin{align} \label{eq:Stokes.reduced}
\left \{
\begin{array}{rcl}
\displaystyle -\mu \Delta u_n+ \nabla p_n&=& 0 \quad \text{ in }\mathbb{R}^3 \setminus \bigcup_{i=1}^n \B_i(t),\\
\dv u_n&=&0\quad \text{ in }\mathbb{R}^3 \setminus \bigcup_{i=1}^n \B_i(t),\\
\displaystyle D u_n&=&0\quad  \text{ in } \bigcup_{i=1}^n  \B_i(t), \\
\displaystyle \int_{\partial \mathcal B_i(t)} \Sigma(u_n,p_n) \nu & = &0 , \\
\displaystyle \int_{\partial \mathcal B_i(t)} [\Sigma(u_n,p_n) \nu ] \times (x-x_i) & = & \sqrt{2 k_B \Theta \mu r^{3} \mR_2(\xi_i)}  \mlnew{\circ \dot B_i}(t), 
\end{array}
\right.
\end{align}
\begin{align} \label{eq:Particles.reduced}
\left \{
\begin{array}{rcl}
\displaystyle	\dd \xi_i(t) &=& \xi_i(t) \times \sqrt{2 k_B \Theta \mu^{-1} r^{-3} \gamma_{rot}^{-1}} \circ  \dd B_i(t) + P_{\xi_i^\perp} h(t,\xi_i(t),x_i) \dd t, \\
\xi_i(0) &=& \xi_{i,0}, \\
x_i(t) & = & x_{i,0}.
\end{array}
\right.
\end{align}
\end{subequations}
Here we used \eqref{eq:R_2} and properties of the cross product to simplify the equation for $\xi_i$.
Moreover, for the ease of notation, we replaced the condition $u_n = v_i + \omega_i \times (x - x_i)$ in $\B_i$ by the equivalent condition $D u_n = 0$ in $\B_i$.



\subsection{Nondimensionalization} \label{sec:nondim}
We examine the expected order of magnitude of various terms and perform a nondimensionalization of the equation. A similar reasoning (for a model including translations and also gravity) can be found in \cite{OttoTzavaras08}. 

We recall that we make the assumption that we are in the so-called dilute regime,
meaning
\begin{align} \label{eq:diluteness}
	\phi_n := \frac{n r^3}{L^3} \ll 1,
\end{align}
where $L$ is the characteristic length of the cloud of particles such that $n/L^3$ is the number density. This  assumption entails that the particles can freely rotate.
Note that $\phi_n$ is proportional to the particle volume fraction. For very elongated rods the reference particle $\B$ has a very small volume such that the  particle volume fraction might be much less than $\phi_n$.
However, we fix the reference particle $\B$ to be independent of $n$. Throughout the paper we will often simply refer to $\phi_n$ as the particle volume fraction.

%

From \eqref{eq:Particles.reduced} we obtain a rotational diffusion constant
$$D_r \sim \frac{k_B \Theta}{\mu \gamma_{rot} r^3},
$$
this gives rise to a typical timescale for diffusion
$$T_D =  \frac 1 {D_r} \sim \frac{\mu r^3 \gamma_{rot}}{k_B \Theta}.
$$

On the other hand, we look at the viscoelastic stress $\sigma$. 
Following the heuristics given in Subsection \ref{sec:heuristics}, we remind that a non-isotropic particle $\B_i$ induces a stresslet on the fluid proportional to the torque. 
The average stresslet produced by one particle arises then from the variation of the stresslet with respect to changes of the orientation.
Thus combining the the formula for the torque from \eqref{eq:Stokes.reduced} with the random part of the change of orientation in \eqref{eq:Particles.reduced}, formally leads to a stresslet of order
$$|S_{i}|  \sim \sqrt{k_B \Theta \mu^{-1} r^{-3}} \sqrt{k_B \Theta \mu r^3} = k_B \Theta,
$$
for each individual particle. For a rigorous argument on how the individual stress arises, we refer to Lemma \ref{lem:Expectations}.
To obtain the total \rhnew{viscoelastic} stress, we multiply with the number density $|\sigma| \sim \frac n {L^3} |S_i|$.
Since the induced fluid gradient is of order $
|\nabla u| \sim \frac{|\sigma|}{\mu}$ we arrive at the viscoelastic timescale
$$T_u = \frac 1 {|\nabla u|} \sim \frac{\mu \gamma_{rot} L^3}{n k_B \Theta}.
$$

In particular, we have 
$$\frac{T_D}{T_u} \sim \frac{n r^3}{L^3} = \phi_n \ll 1,
$$
which means that the diffusion happens on a much faster timescale than the viscoelastic stress.
Consequently, when nondimensionalizing, we can choose to rescale to the diffusive timescale $T=T_D$ or to rescale to the viscoelastic timescale $T= T_u$.



We first rescale \eqref{eq:reduced} 
with the characteristic time $T_D$ and the length $L$.
Keeping the same symbols for the rescaled quantities, lengthy but straightforward calculations\footnote{Here one needs to use that the Brownian motion scales as $B(Tt) \sim \sqrt T B(t)$.} yield  
\begin{subequations} \label{eq:micro.T_D}
\begin{equation}  \label{eq:Stokes.micro.T_D}
\left \{
\begin{array}{rcl}
- \Delta u_n+ \nabla p_n&=& 0 \quad \text{ in }\mathbb{R}^3 \setminus \bigcup_{i=1}^n \B_i(t),\\
\dv u_n&=&0 \quad  \text{ in }\mathbb{R}^3 \setminus \bigcup_{i=1}^n \B_i(t),\\
D u_n&=&0 \quad  \text{ in } \bigcup_{i=1}^n  \B_i(t), \\
\displaystyle \int_{\partial B_i(t)} \Sigma(u_n,p_n) \nu & = &0  \quad  \text{ for all } 1 \leq i \leq n, \\ 
\displaystyle  \int_{\partial B_i(t)} [\Sigma(u_n,p_n) \nu ] \times (x-x_i) &=&  r^3 \sqrt{ 2 \gamma_{rot} \mR_2(\xi_i(t))}  \mlnew{\circ \dot B_i(t)}  \quad  \text{ for all } 1 \leq i \leq n, 
\end{array}
\right.
\end{equation}
\begin{align} \label{eq:Particles.T_D}
\left \{
\begin{array}{rcl}
\displaystyle	  \dd \xi_i(t) &=& \sqrt{2} \xi_i(t) \times \circ \dd B_i(s) + P_{\xi_i^\perp} h(t,\xi_i(t),x_i) \dd t ,  \\
\xi_i(0) &=& \xi_{i,0}.
\end{array}
\right.
\end{align}
\end{subequations}
Note that from now on the nondimensional fluid stress is given as $\Sigma(u_n,p_n) = 2 D u_n - p_n \Id$.
Note that due to the rescaling of the lengthscale with $L$, the  ``volume fraction'' $\phi_n$ is now given by 
\begin{align}
	\phi_n = n r^3.
\end{align}
Finally, we remark that we dropped the \rhnew{trivial equation $x_i(s) = x_{i,0}$, and we consider instead the positions $x_i$ as given time-independent quantities.}

Similarly rescaling instead  
with the characteristic time $T_u$ yields

\begin{subequations} \label{eq:micro.T_u}
\begin{equation}  \label{eq:Stokes.micro.T_u}
\left \{
\begin{array}{rcl}
- \Delta u_n+ \nabla p_n&=& 0\quad  \text{ in }\mathbb{R}^3 \setminus \bigcup_{i=1}^n \B_i(t),\\
\dv u_n&=&0\quad  \text{ in }\mathbb{R}^3 \setminus \bigcup_{i=1}^n \B_i(t),\\
D u_n&=&0\quad  \text{ in } \bigcup_{i=1}^n  \B_i(t), \\
\displaystyle \int_{\partial B_i(t)} \Sigma(u_n,p_n) \nu & = &0 \quad  \text{ for all } 1 \leq i \leq n , \\
\displaystyle \int_{\partial B_i(t)} [\Sigma(u_n,p_n) \nu ] \times (x-x_i) &=& \frac 1 n \sqrt{2 \gamma_{rot} \phi_n  \mR_2(\xi_i(t))} \mlnew{ \circ \dot B_i(t)}  \quad  \text{ for all } 1 \leq i \leq n, 
\end{array}
\right.
\end{equation}
\begin{align} \label{eq:Particles.T_u}
\left \{
\begin{array}{rcl}
\displaystyle 	 \dd \xi_i(t) &=& \xi_i(t) \times \sqrt{\frac {2} {\phi_n}}   \circ \dd B_i(t) + \frac 1 \phi_n P_{\xi_i^\perp} h(t,\xi_i(t),x_i) \dd t ,  \\
\xi_i(0) &=& \xi_{i,0}.
\end{array}
\right.
\end{align}
\end{subequations}

\rh{Clearly \eqref{eq:micro.T_D} can be recoverd from \eqref{eq:micro.T_u} upon rescaling time with $\phi_n$.
Thus, in the limit $n\to \infty$ with $\phi_n \to 0$ one can 
interpret \eqref{eq:micro.T_D} as an initial layer of \eqref{eq:micro.T_u}.
We emphasize that in this sense, the function $h$ in \eqref{eq:micro.T_D} corresponds to $h(\phi_n t, \cdot)$ for $h$ as in \eqref{eq:micro.T_u}. 

}


\section{Main Results}\label{section3}

\subsection{Assumptions} \label{sec:assumptions}

\rh{We will impose the following assumptions for the rest of the paper.}
\ml{We work in the framework of a filtered probability space, denoted by $(\Omega, \mathcal{F}, \left\{\mathcal{F}_t\right\}, \mathbb{P})$.}
We recall that we consider $n$ particles occupying $\B_i(\xi_i) = x_i + r R_i(\xi_i) \B$, where $\B \subset \R^3$ is the reference particle specified in \eqref{rod_reference},  the positions $x_i \in \R^3$ are static random variables, the orientations $\xi_i \in \S^2$ are random and time-dependent, and the rotation matrix $R_i(\xi_i) \in SO(3)$ is any matrix which satisfies $R_i(\xi_i) e_3 = \xi_i$.
We emphasize that $\B_i, x_i, r, \xi_i$ all implicitly depend on $n$.

We assume that the particle volume fraction $\phi_n = nr^3$ tends to zero sufficiently fast, namely
\begin{align} \label{ass:phi.log.n}\tag{H1}
	\lim_{n \to \infty} \phi_n \log n = 0.
\end{align}
Moreover, we assume that the particle centers $x_i$ are well-separated in the sense that there exists $c > 0 $, independent of $n$, such that
\begin{align} \label{ass:well.separated} \tag{H2}
	\dmin := \min_{i \neq j} |x_i - x_j| \geq c n^{-\frac 1 3}.
\end{align}
We remark that the above assumptions together imply in particular that there exists $n_0 \in \N$
such that for all $n \geq n_0$
\begin{align} \label{eq:separation}
	 \min_{i \neq j} \inf_{(\xi_i, \xi_j) \in \S^2 \times \S^2} \dist(\B_i, \B_j) \geq cr.
\end{align} 
 Note that in the assumption above, we use the (abusive) notation $\xi_i$ to refer to any variable in $\S^2$.
We assume in addition that the particles are contained in a bounded domain uniformly in $n$
\begin{equation}\label{ass:uniform_bound}
    \underset{n}{\sup} \underset{1\leq i \leq n}{\max}|x_i| < +\infty. \tag{H3}
\end{equation}

We assume that $B_i$ are independent $\R^3$ valued Brownian motions ,all independent of $\xi_{i,0}$,  defined on the filtered probability space $(\Omega, \mathcal{F}, \left\{\mathcal{F}_t\right\}, \mathbb{P})$. 
Moreover, we assume that the initial particle orientations
$\xi_{i,0}$ are independent random variables 
  and we assume the convergence of the initial empirical measure to some $f_0 \in \mathcal P_1(\R^3 \times \S^2) \cap L^2(\R^3\times \S^2)$,\footnote{To see that for a given compactly supported function $f_0$, there exists $x_i, \xi_{i,0}$ which satisfy both assumptions \eqref{ass:well.separated} and \eqref{ass:initial.convergence}, one might first generate i.i.d variables $\tilde x_i, \xi_{i,0}$ with law  $f_0 \dd x$ and then  define the positions $x_i \in n^{-1/3}\Z^3$ to be the closest points to $\tilde x_i$.

} 
i.e.
\begin{align} \label{ass:initial.convergence} \tag{H4}
     \lim_{n \to \infty}\E\left[  \mathcal W_1\left(f_0, \frac 1 n  \sum_i \delta_{x_i,\xi_{i,0}}\right)  \right].
\end{align}


Here, for two probability measures $g_1,g_2 \in \mP_1(\R^3 \times \S^2)$ \rh{with bounded first moment}, we denote by $\mW_1(g_1,g_2)$
their $1$-Wasserstein distance. \rh{It is well known that the $\mW_1$-distance metrizes weak convergence of measures.}


Finally, we assume for simplicity that the external force $h$ is smooth, i.e.
\begin{align} \label{ass:h} \tag{H5}	h \in C^\infty([0,\infty) \times \R^3 \times \S^2).
\end{align}


\subsection{Notion of solution and well-posedness of the microscopic system}\label{sec:well-posedness}
We will now specify what are the notion of solutions to the systems \eqref{eq:micro.T_u} and \eqref{eq:micro.T_D}.
The uncoupled SDEs for the particle orientations are well known to be well-posed for $h$ as in \eqref{ass:h} (see Theorem \ref{thm:SDE.M}).
 Note that by introducing the  diffusion matrix $\sigma_D$,
\begin{equation}\label{eq:sigma}
 \sigma_D(\xi)  := \sqrt{2} \begin{pmatrix}
 0 & -\xi_3 & \xi_2\\
\xi_3 & 0 & - \xi_1 \\
-\xi_2 & \xi_1 & 0  
\end{pmatrix},
\end{equation}
and by the Ito-Stratonovich conversion, Remark \ref{rem:SDE.Strato.Ito}, we can rewrite \eqref{eq:Particles.T_D} and \eqref{eq:Particles.T_u} respectively as
\begin{equation}\label{eq:strato.xi_i.T_D}
\dd \xi_i =   \sigma_D(\xi_i) \dd B_i -2\xi_i \dd s+ P_{\xi_i^\perp} h(s,\xi_i,x_i) \dd s , 
\end{equation}
\begin{equation}\label{eq:strato.xi_i.T_u}
\dd \xi_i (s)= \sqrt{\frac {1} {\phi_n}}\sigma_D(\xi_i(s)) \dd B_i(s)- {\frac {2} {\phi_n}}\xi_i(s) \dd s + \frac 1 \phi_n P_{\xi_i^\perp} h(s,\xi_i(s),x_i) \dd s. 
\end{equation}


Since the fluid is modeled by the stationary Stokes equations, we can only expect the fluid velocity to be pathwise as (ir-)regular in time as the white noise that drives the fluid velocity through the prescribed torques. It is well known that white noise is in $H^{-s}$ for any $s > 1/2$.

To give a meaning to $u_n$ and to obtain this regularity, we write $u_n$ as the distributional derivative of a suitable stochastic integral.
To this end, we first introduce the operator $L_n \colon \R^{3n} \to \L(\R^{3n}, \dot H^1(\R^3))$. For any fixed set of particle positions $(x_1, \dots, x_n)\in \R^{3n}$ it associates to
every set of orientations $(\xi_1,\dots \xi_n)\in \S^{2n}$
the solution to the Stokes system  with given torques $(T_1, \dots, T_n)$. More precisely,
\begin{align}\label{def:L_n.1}
	v = L_n((\xi_1,\dots, \xi_n))(T_1, \dots, T_n) , 
\end{align}
is defined to be the solution to 
\begin{equation}  \label{def:L_n.2}
\left \{
\begin{array}{rcl}
- \Delta v+ \nabla q&=& 0 \quad \text{ in }\mathbb{R}^3 \setminus \bigcup_{i=1}^n \B_i,\\
\dv v&=&0 \quad \text{ in }\mathbb{R}^3 \setminus \bigcup_{i=1}^n \B_i,\\
D v&=&0 \quad \text{ in } \bigcup_{i=1}^n  \B_i, \\
\displaystyle \int_{\partial \B_i} \Sigma(v,q) \nu & = &0, \\
\displaystyle \int_{\partial \B_i} [\Sigma(v,q) \nu ] \times (x-x_i) &=& T_i,
\end{array}
\right.
\end{equation}
where $\B_i =x_i + r R(\xi_i) \B$,  we emphasize here the abusive notation $(\xi_1,\cdots,\xi_n)$ (resp. $(x_1,\cdots, x_n)$) to represent any configuration of $n$ orientations in $\S^{2}$ (resp. positions in $\R^3$).
 It is straightforward to see that  for non-overlapping particles $\B_i$, the linear problem \eqref{def:L_n.1} admits a unique weak solution $v \in \dot H^1(\R^3)$ (see e.g. \cite{NiethammerSchubert19}). 
In particular, $L_n$ is well-defined for all particle configurations $(x_1, \dots, x_n)$ which satisfy
\eqref{eq:separation}.  We can now define the velocity field $u_n$ as follows
\begin{defi}[Definition of the velocity field $u_n$]
Let $n \in \N^*$, $t\geq 0$ and $\Xi = (\xi_1 ,\dots \xi_n)$ with $\xi_i$ the solutions to \eqref{eq:Particles.T_u}. We define $U_n$ as the following stochastic integral
\begin{align} \label{def:U_n}
	U_n(t) := \frac {\sqrt{\am{2\gamma_{rot}}\phi_n}} n \int_0^t L_n(\Xi(s)) \sqrt{\mathfrak R_2(\Xi(s))} \circ \dd (B_1(s), \dots, B_n(s)),
\end{align}
where  $\mathfrak R_2(\Xi) \in \R^{3n \times 3n}$ is the block-diagonal matrix with blocks $\mR_2(\xi_i)$ (see \eqref{eq:R_2}).
The (distributional) solution to \eqref{eq:Stokes.micro.T_u} $u_n$ is defined as the distributional derivative of $U_n$,
\begin{align} \label{eq:u_n.U_n}
	u_n := U_n'.
\end{align}
Similarly, concerning the system \eqref{eq:micro.T_D}, we set
\begin{align} \label{def:U_n.T_D}
	 U_n(t) &:= r^3 \am{\sqrt{2\gamma_{rot}}} \int_0^t L_n(\Xi(s)) \sqrt{\mathfrak R_2(\Xi(s))} \circ \dd (B_1(s), \dots, B_n(s)), \\ 
	u_n &:= U_n',
\end{align}
where $\xi_i$ are the solutions to \eqref{eq:Particles.T_D}.
\end{defi}
We make precise the definition of the integrals in \eqref{def:U_n} and \eqref{def:U_n.T_D} in Appendix \ref{appendix_stochastics}, where we collect some statements about the definition and properties of such stochastic integrals with values in separable Hilbert spaces. Essentially, all standard results immediately carry over due to It\^o isometry. Therefore, $U_n$, and thereby $u_n$, is well-defined provided we show that $L_n$ is continuously differentiable with respect to the orientations $\xi_i$.

We will show the following global well-posedness result of the microscopic dynamics.

\begin{thm}	 \label{th:well-posed}
	Let \eqref{ass:phi.log.n}--\eqref{ass:well.separated} be satisfied. Then, for all $n \in \N$, there exists a unique solution $(\xi_1,\dots,\xi_n)$ to \eqref{eq:Particles.T_D} and \eqref{eq:Particles.T_u}, respectively. 
	
	Moreover, there exists $N_0 \in \N$ such that for all $n \geq N_0$ and all $s > 1/2$,
	the operator $L_n$ defined through \eqref{def:L_n.1}--\eqref{def:L_n.2} satisfies $L_n \in C^1((\S^2)^n; \L(\R^{3n},H^{-s}_{\mathfrak s}(\R^3)))$, where $H^{-s}_{\mathfrak s}(\R^3)$ denotes the subspace of divergence free functions in $H^{-s}(\R^3)$.
	
	In particular, the integral \eqref{def:U_n} (respectively \eqref{def:U_n.T_D}) is well defined and for all $T>0$
	\begin{align}
		 u_n \in L^2(\Omega;H^{-s}(0,T;H^{-s}_{\mathfrak s}(\R^3))).
	\end{align}	 
\end{thm}


\subsection{Convergence results} \label{sec:convergence}

We are finally prepared to state the main results of our paper, the convergence to the mean-field limits for \eqref{eq:micro.T_u} and \eqref{eq:micro.T_D}.

We denote by $S_n$ the empirical measure 
\begin{equation}\label{eq:empirical_measure}
	S_n(t) := \frac 1 n \sum_i \delta_{x_i,\xi_i(t)}.
\end{equation}


We first state the main result concerning system \eqref{eq:micro.T_D}.
\begin{thm} \label{th:main.T_D}
	Let assumptions \eqref{ass:phi.log.n}--\eqref{ass:h} be satisfied. For $n \geq N_0$ as in Theorem \ref{th:well-posed}, let $\xi_i$, $1 \leq i \leq n$ and $u_n$ be the unique solutions to \eqref{eq:micro.T_D}. 
Then, for all $t> 0$ and all $s > 1/2$ 
the following convergence in probability holds:
\begin{align}
	\forall \eps > 0 ~ \lim_{n \to \infty} \P\left( \|\phi_n^{-1} u_n - u\|_{H^{-s}((0,t), H^{-s}(\R^3))} + \sup_{\tau \in [0,t]} \mW_1(S_n(\tau),f(\tau)) > \eps\right)  = 0,
\end{align}
\am{where $f\in C([0,t], \mathcal P_1(\R^3 \times \S^2)\cap L^2(\R^3 \times \S^2))$ is the unique weak solution to \eqref{eq:Fokker-Planck.instationary} such that for almost all $x\in \R^3$, $f(\cdot,\cdot,x)\in L^2(0,t;H^1(\S^2))$, $f'(\cdot,\cdot,x)\in L^2(0,t;H^{-1}(\S^2))$ and $u\in L^2(0,t;\dot{H}^1_{\mathfrak{s}}(\R^3))$ the unique weak solution to \eqref{eq:viscoelastic}.}
 \end{thm}
 
Similarly, we show the following convergence result for
system \eqref{eq:micro.T_u}.
\begin{thm} \label{th:main.T_u}
	Let assumptions \eqref{ass:phi.log.n}--\eqref{ass:h} be satisfied. For each $n \geq N_0$  as in Theorem \ref{th:well-posed}, let $\xi_i$, $1 \leq i \leq n$ and $u_n$ be the unique solutions to \eqref{eq:micro.T_u}. 
Then, for all $s_1 > 1/2, s_2 > 0, s_3 > 3/2$
the following convergence in probability holds:
\begin{align}
	\forall \eps > 0 ~ \lim_{n \to \infty} \P\left(\|u_n - u\|_{H^{-s_1}((0,t) ;H^{-s_1}(\R^3))} + \|S_n - f\|_{H^{-s_2}(0,t;H^{-s_3}(\R^3 \times \S^2))}  > \eps\right)  = 0,
\end{align}
\am{where $f\in L^2((0,t)\times \S^2\times \R^3)$ is the unique weak solution to \eqref{eq:Fokker-Planck.stationary} such that for almost all $x\in \R^3$, $f(\cdot,\cdot,x)\in C^\infty((0,t)\times\S^2)$, and $u\in L^2(0,t;\dot{H}^1_{\mathfrak{s}}(\R^3))$ the unique weak solution to \eqref{eq:viscoelastic}.}
 \end{thm}

\subsection{Notation}\label{subsec:notatation}

\begin{itemize}
    \item 
Throughout the paper we will often deal with vectors $V \in (\R^3)^n$, where $n \in \N$ is the number of particles. We will use the convention to denote the components $V_i \in \R^3$, $1 \leq i \leq n$ by Latin indices and use Greek indices to denote the components of these vectors, e.g. $V_{i,\alpha} \in \R$, $1 \leq \alpha \leq 3$. With this convention, we allow for the slight abuse to write $V \in \R^{3n}$ instead of $V \in (\R^3)^n$.\am{ In particular we denote by $e_{i,\alpha}$, $1 \leq i \leq n$, $1 \leq \alpha \leq 3$ the canonical basis of $\R^{3n}$.
}

Moreover, in the special case of the particle positions and orientations we will use capital letters to denote the collections $X = (x_1, \dots, x_n) \in \R^{3n}$, and $\Xi = (\xi_1, \dots, \xi_n) \in (\S^2)^n$ in order to avoid confusion with variables $x,\xi$ that appear in the limit systems.

\item \am{For any vector $T \in \R^3$, we set \begin{align}\label{eq:[T]_M}
    [T]_M := \begin{pmatrix}
 0 & -T_3 & T_2\\
T_3 & 0 & - T_1 \\
-T_2 & T_1 & 0  
\end{pmatrix},
\end{align}
the unique skew-symmetric matrix associated to $T\in \R^3$ satisfying for any $x\in\R^3$ and any smooth test function $\phi\in {C}_c^\infty(\R^3;\R^3)$
\begin{align}\label{eq:skew_curl}
    [T]_M x = T \times x,&& [T]_M:\nabla \phi(x)= 2 T \cdot \curl \phi(x).
\end{align}
}
\item We denote by $\Sym_0(3)$ the space of symmetric traceless matrices $A \in \R^{3\times3}$ and by  $\sym B$ (resp. $\sym_0 B$) the symmetric (resp. symmetric traceless) part of $B \in \R^{3\times3}$. 
\item For any vector valued function $u$ we set $Du:= \sym \nabla u$ the symmetric gradient of $u$.
\item For a reflexive Banach space $X$, an open set $\mathcal O \subset \R^d$, $p \in (1,\infty)$, $s \in \R$, we denote by $W^{s,p}(\mathcal O;X)$ the usual fractional Sobolev space of $X$-valued functions. For $s = k+ \gamma$, $k \in \N$, $\gamma \in (0,1)$ the norm in $W^{s,p}(\mathcal O;X)$ is given by
\begin{align}
    \|u\|^p_{W^{s,p}(\mathcal O;X)} := \|u\|^p_{W^{k,p}(\mathcal O;X)} + [\nabla^k u]^p_{s,p}, \\ 
    [f]^p_{s,p} := \int_{\mathcal O \times \mathcal O} \frac{\|f(x) - f(y)\|_X^p}{|x-y|^{d+\gamma p}}  \dd y \dd x.
\end{align}
For $s > 0$, we define $W^{-s,p'}(\mathcal O;X) := (W^{s,p}_0(\mathcal O;X'))^\ast$.
We recall (see e.g. \cite{Amann00}) that if $X$ is reflexiv, then an equivalent norm in $W^{-k + \theta,p}(\mathcal O;X)$, $\theta \in (0,1)$ is given by
\begin{align} \label{eq:char.negative.sobolev}
    \|u\|_{W^{-k + \theta,p}(\mathcal O;X)} = \inf \left\{ \sum_{|\beta| \leq k} \|u_\beta\|_{{W^{\theta,p}(\mathcal O;X)}} : u= \sum \partial^\beta u_\beta \right\}.
\end{align}

We also introduce weighted fractional Sobolev spaces $W^{s,p}_w(\mathcal O)$ with a weight $w \geq 0$ through
\begin{align}
    \|v\|^p_{W^{s,p}_w(\mathcal O)} &:= \|v\|^p_{W^{k,p}_w(\mathcal O)} + [\nabla^k v]^p_{w,s,p} ,\\ 
    [f]^p_{w,s,p} &:= \int_{\mathcal O \times \mathcal O} \frac{|f(x) - f(y)|_X^p}{|x-y|^{d+\gamma p}} w(x) w(y) \dd y \dd x, \\
    \|v\|^p_{W^{k,p}_w(\mathcal O)} &:= \sum_{l = 0}^k \int_{\mathcal O} |\nabla^l v|^p w \dd x,
\end{align}
for $s = k + \gamma$ as above, and 
$W^{-s,p'}_{\frac 1 w}(\mathcal O) := (W^{s,p}_{0,w}(\mathcal O))^\ast$.

Throughout the paper, we will fix the appearing weight function  $w$ as
 \begin{align} \label{def:weight.function}
w(x)=(1+|x|)^a,  && \text{for some } 0<a<1.
\end{align}
\rh{As usual we denote the Hilbert spaces $H^s(\mathcal O) := W^{s,2}(\mathcal O)$ and $H^s_w(\mathcal O) = W^{s,2}_w(\mathcal O)$.}

\item We allow for the abuse of notation 
\begin{align}
    \|g\|_{W^{s,p}_{\loc}(\R^3)} \leq C, 
\end{align} 
to indicate that for all compact $K \in \R^3$ there exists $C$ depending on $K$ such that 
\begin{align}
    \|g\|_{W^{s,p}(K)} \leq C.
\end{align} 

Similarly, we write  $g \in {W^{s,p_-}(U)}$ and 
\begin{align}
    \|g\|_{W^{s,p_-}(U)} \leq C ,
\end{align} 
to indicate that for all $1 \leq q < p$ there exists $C$ depending on $q$ such that
\begin{align}
    \|g\|_{W^{s,q}(U)} \leq C.
\end{align} 
The notation $\|g\|_{W^{s_-,p}(U)} \leq C.$ should be understood analogously.

\item \rh{For an open set $\mathcal O \subset \R^3$, $p \geq 1$ and $w$ as in \eqref{def:weight.function}, we introduce
\begin{align}
   \|g\|_{L^{p,2}_{w,\mathcal O}} &:= \|g\|_{L^p(\mathcal O)} + \|g\|_{L^2_w(\R^3 \setminus \mathcal O)}, \\
    L^{p,2}_{w,\mathcal O} &:= \{ g \in L^1_\loc(\R^3) : \|g\|_{L^{p,2}_{w,\mathcal O}} < \infty \}. \label{def:L^p,2_w,K}
\end{align}
}
\item \rh{
    We use an index $\mathfrak s$ to denote subspaces of (distributionally) divergence free functions, i.e. we use the notation $L^p_{\mathfrak s}(\mathcal O)$, $H^s_{\mathfrak s}(\mathcal O)$, $H^s_{\mathfrak s,w}(\mathcal O)$, $W^{s,p}_{\mathfrak s}(\mathcal O)$,  $W^{s,p}_{\mathfrak s,w}(\mathcal O)$, $L^{p,2}_{\mathfrak s,w,\mathcal O}$ and all these subspaces are closed.
}

\item \am{Finally, in Sections \ref{sec:L_n} -- \ref{section6} we will adopt the usual convention to denote by $C>0$ any  constant that might change from line to line and that might depend on certain quantities which are fixed like the reference particle or the constant appearing in \eqref{ass:well.separated}. We emphasize that the constant $C$ will never depend on $n$ though.}
\end{itemize}

 \subsection{Outline of the proof of the main results}\label{subsec:proof_strategy}

Starting from \eqref{eq:micro.T_D}, the first step towards the derivation of the mesoscopic system \eqref{eq:Fokker-Planck.instationary}, \eqref{eq:viscoelastic} is the following approximation for the fluid equations \eqref{eq:Stokes.micro.T_D}:
\begin{align} \label{eq:u_n.tilde}
	- \Delta u_{n,app} + \nabla \tilde p_{n,app} = \frac{\phi_n} n  \sum_i ([T_i]_M + \mS(\xi_i)T_i)   \nabla \delta_{x_i},
\end{align}
where $T_i = \sqrt{\am{2\gamma_{rot}}\mR_2(\xi_i)} \circ \dot B_i$ and we recall that $\mS(\xi_i)$ from \eqref{def:mS_i} is the tensor that relates the torque to the stresslet and $[T_i]_M$ is defined in \eqref{eq:[T]_M}.
Thus the approximation \eqref{eq:u_n.tilde} consists in a formal superposition principle of point torques and stresslets; each particle contributes to $u_{n,app}$ as if it was alone in the fluid and only its total torque and stresslet acts on the fluid at the particle center.

The rigorous justification of this approximation is the content of Section \ref{sec:L_n}. More precisely,
Section \ref{sec:L_n} deals with the approximation of the solution $L_n T$ to \eqref{def:L_n.2} by the solution $ L_{n,app} T$ to \eqref{eq:u_n.tilde} for any given $T \in \R^{3n}$ and any given collection of particle positions $x_i$ and orientations $\xi_i$, $1 \leq i\leq n$ which satisfy assumotion \eqref{ass:phi.log.n}--\eqref{ass:uniform_bound}. 
To this end, we will first introduce an intermediate approximation $L^{im}_{n,app}$ defined as the superposition of single particle problems, but with the full  boundary conditions (no-slip and balance equations) for the single particle. 

Estimates for approximations similar to \eqref{eq:u_n.tilde}  have been given for example in \cite{HoferVelazquez18, Hofer18MeanField, Gerard-Varet19,  NiethammerSchubert19, Gerard-VaretHillairet19, HillairetWu19}. 
The main novelty here are estimates for $\nabla \xi_i (L_n - L_{n,app}^{im}) T$,
which are needed because of the Stratonovitch integral in \eqref{def:U_n}.
Since the fluid domain $\R^3 \setminus \cup_i \B_i$ depends on the particle orientations, we enter here the topic of shape derivatives. Such shape derivatives for the Stokes equations with different boundary conditions have been considered for example in \cite{Simon91,BadraCaubetDambrine11}.  
Instead of identifying the boundary value problem solved by $\nabla_{\xi_i} L_n$ as in \cite{Simon91,BadraCaubetDambrine11}, we rely on the method of reflections to analyze these shape derivatives in Subsection \ref{sec:reflections}. The method of reflections has been used for related problems, see for example \cite{HoferVelazquez18, Hofer18MeanField,  NiethammerSchubert19, Gerard-VaretHillairet19, HillairetWu19, Mecherbet18, Mecherbet19Cluster}.
It allows to express $L_n$ in terms of single particle problems only. Since these single particle operators have an explicit dependence on the particle orientation,
this yields a useful expression for the shape derivative.

As reflected by the approximation \eqref{eq:u_n.tilde} the given torques at each particle position $x_i$ perturb the fluid velocity in a singular way as $|x-x_i|^{-2}$. This is the reason why we only obtain sufficient estimates under the assumptions \eqref{ass:phi.log.n}--\eqref{ass:well.separated}. Moreover, due to the singular nature of the perturbation, these estimates are locally only in $L^p$ for $p < 3/2$.
Since the stochastic integral \eqref{def:U_n} does not seem compatible with such $L^p$ spaces (see e.g. \cite{van2015stochastic}),  
we work instead in weighted negative Sobolev spaces $H^{-s}_w(\R^3)$, $s > 1/2$, which are Hilbert spaces. In Appendix \ref{sec:weighted sobolev space}, we show that ${L^{p,2}_{w,K}}$ (see \eqref{def:L^p,2_w,K}) embeds into $H^{-s}_w(\R^3)$ for a compact set $K$.
Therefore, in Section \ref{sec:L_n}, we will always work in the space ${L^{p,2}_{w,K}}$ where we choose $K$ to contain all the particles thanks to assumption \ref{ass:uniform_bound}.

The results in Section \ref{sec:L_n} immediately imply that the stochastic integral \eqref{def:U_n} that defines $u_n$ is well-defined which yields Theorem \ref{th:well-posed}.

The approximation \eqref{eq:u_n.tilde} suggests 
through an appropriate version of the Law of Large Numbers that the limit of $u_n$ is given in terms of
the expectation of the  torque and stresslet exerted on the fluid by each particle. 
The following lemma states that these expectation indeed correspond to the formula for the viscoelastic stress in \eqref{eq:viscoelastic}. The proof of the lemma is a straightforward calculation which we postpone to Appendix \ref{sec:Expectations}.

\begin{lem} \label{lem:Expectations}
Let $\xi_i$ be the solution to \eqref{eq:Particles.T_D}.
\am{For any $A\in {C}^\infty_c((0,t), \Sym_0(3))$ and $b\in {C}^\infty_c((0,t), \R^{3})$ we have
\begin{align}
    \E \left[\int_0^t b(s) \cdot   \sqrt{ \mR_2(\xi_i(s))} \circ \dd B_i(s) \right] &= 0 , \label{eq:torque.average} \\
     \E \left[ \int_0^t A(s): \mS(\xi_i(s)) \sqrt{ \mR_2(\xi_i(s))} \circ \dd B_i(s) \right] &= \frac{\gamma_E}{\sqrt{2\gamma_{rot}}} \E\int_0^t[3 \xi_i(s) \otimes \xi_i(s) - \Id]:A(s) \dd s.\\ \label{eq:stress.average}
\end{align}
}

\end{lem}


The proof of the main convergence results, Theorems \ref{th:main.T_D} and \ref{th:main.T_u}, is completed in Sections \ref{section5} and \ref{section6}, respectively. 
Since the particle dynamics is uncoupled, we do not face here the problem of propagation of chaos as in classical mean field systems. However, a slight complication arises because the particle position are not assumed to be independent (which is impossible due to assumption \eqref{ass:well.separated}). If they were independent, the law of each particle would be given by the Fokker-Planck equation \eqref{eq:Fokker-Planck.instationary} in the case of system \eqref{eq:Particles.T_D}, and we would deduce convergence immediately by the Law of Large Numbers. 
Instead, we follow here a compactness approach used for example in \cite{meleard1996asymptotic,oelschlager1985law,kipnis1998scaling,flandoli2021navier}. 
This approach is roughly described as follows. First, one shows tightness in a suitable functional spaces of the laws of the empirical measure $S_n$ and the fluid velocity $u_n$ in suitable spaces. Then one introduces functionals such that, on the one hand,  $u_n$ and $S_n$ concentrate on the zeroes of these functionals, and on the other hand, these zeroes are precisely the unique solutions of the desired limit systems \eqref{eq:viscoelastic} together with \eqref{eq:Fokker-Planck.instationary} and \eqref{eq:Fokker-Planck.stationary}, respectively.

The choice of the functionals correspond to distributional solutions of the limit system. We therefore need to show that such solutions are unique in the negative Sobolev spaces we use.  The proof of these uniqueness results, which we give in  Appendix \ref{appendixA} for completeness, is based on well-posedness and regularity of the corresponding dual problems. 
This is completely standard for the Stokes equations and the instationary Fokker-Planck equations. The slightly more involved proof for the stationary Fokker-Planck system we carry out in more detail.

\subsection{Limitations and possible generalizations}\label{subsec:discussion}


In this subsection we comment on open questions related to limitations and possible generalizations of the analysis in this paper.

\medskip

We dropped the evolution of the particle translations. It seems very challenging to include translations with the current techniques because they require well-separated particles in the sense of assumption \eqref{ass:well.separated}. Although this condition has been shown to propagate in time in \cite{Hofer18MeanField,Mecherbet18, Hofer&Schubert} under suitable assumptions for sedimenting inertialess non-Brownian rigid spherical particles, the presence of Brownian forces and torques in the current model could break such propagation.

\medskip

As discussed in Subsection \ref{sec:dynamics} it seems physically more accurate to prescribe random forces and torques given by \eqref{eq:full.Stokes.Einstein} and the  evolution of the particle evolution through \eqref{eq:angularVelocity}. For vanishing particle volume fraction as in assumption \eqref{ass:phi.log.n}, it seems not unrealistic to treat such a microscopic model (still ignoring the evolution of the particle centers). One might worry about rotations caused by hydrodynamic interactions between the particles which have a singularity like $|x|^{-3}$: a  stresslet $S$ at a particle at $x_i$ creates a fluid velocity roughly like $\nabla \Phi(x-x_i) : S$ where $\Phi$ is the fundamental solution of the Stokes equations. The induced rotation at another particle then scales like the gradient of this function. 
However, as the viscoelastic stress is of smaller order then the diffusion rate by a factor $\phi_n$, this singular interaction could  still be neglegligible if $\phi_n \to 0$. In the Doi model \eqref{eq:full.model} this corresponds to $\lambda_2 = \lambda_3 = 0$ for $\phi_n \to 0$. 
The situation changes completely when one considers non vanishing volume fractions $\phi_n$. In that case, the singular interaction cannot be neglected. In fact, the critical singularity $|x|^{-3}$ of the interaction leaves little hope that one can pass to the formal mean-field limit that would give rise to the term $\dv_{\xi} (P_\xi^\perp \nabla_x u \xi f)$ in (a version of) the Doi model \eqref{eq:full.model}. Indeed, discrepancies between discrete and continuous convolutions which such singular convolution kernel have been studied in a related setting in \cite{Gerard-VaretHillairet19}.

\medskip

On the technical side, we are restricted to vanishing volume fractions $\phi_n$, more precisely to assumption \eqref{ass:phi.log.n}, because of the factor $\phi_n \log n$ appearing in the estimates for the shape derivatives (see Proposition \ref{pro:L_n}).
It would be desirable to analyze whether this bound is optimal or could be improved under suitable assumptions on the particle configuration.


\medskip

We dropped the fluid inertia and model the fluid by the stationary Stokes equations.
An interesting question would be to investigate the instationary Navier-Stokes equations. Heuristically, the result should be unchanged since on the microscopic scale of the particles the fluid inertia does not matter. Mathematically though, already the problem of giving  a meaning to the microscopic fluid velocity changes completely.

To our knowledge, even in the case of quasistatic non-Brownian particles, no rigorous homogenization results regarding the effective viscosity of suspensions are available for the Navier-Stokes equations. This has been obtained only for prescribed particle velocities (i.e. inertial particles) in \cite{Feireisl2016} where the effective equation contains the typical Brinkman term.

\medskip

In this paper the particles are all obtained from (isotropic) rescaling of a fixed reference particle. 
In order to model rod-like particles, it would be interesting to consider particles which become more and more slender as $n \to \infty$. 
We refer to \cite{HoeferPrangeSueur22} for the analysis of the evolution of such filaments in a fluid.

 \medskip
 
Finally, an important open problem is the rigorous derivation of  models (expressions for the viscoelastic stress to begin with) of flexible polymers mentioned in the introduction.
In the simplest setting one could start with a dumbbell model, where the flexible polymer is just modeled by two rigid balls connected by an (infinitesimally small) spring.

 \section{Approximation of the operator \texorpdfstring{$L_n$}{Ln}} \label{sec:L_n}

In this section, we study the solution operator $L_n$ which has been defined in Subsection \ref{sec:well-posedness} as the solution operator for the Stokes problem \eqref{def:L_n.2}. Throughout this section, we denote by $x_i,\xi_i$, $1 \leq i \leq n$, a generic (time-independent) collection of positions and orientations satisfying \eqref{ass:phi.log.n}--\eqref{ass:uniform_bound}.

We introduce the following explicit approximate solution operator 
\begin{equation} 
\begin{array}{lccc}
	{L}_{n,app}  \colon& \R^{3n}  &\to & \L(\R^{3n},  L^{\frac 3 2_-}_{\mathfrak s, \loc }(\R^3)), \\
 & \Xi:=(\xi_1,\cdots,\xi_n) & \mapsto & L_{n,app}(\Xi),
 \end{array}
\end{equation}
defined for all $T:=(T_1,\cdots,T_n)\in \R^{3n}$ by
\begin{align} \label{def:L_n,app}
\left({L}_{n,app}((\xi_1,\dots, \xi_n))T\right)(x) = -\sum_i ([T_i]_M + \mS(\xi_i) T_i) : \nabla \Phi(x-x_i),
\end{align}
where $\mS(\xi_i)$ is defined in \eqref{def:mS_i}, $[T_i]_M$ is defined in \eqref{eq:[T]_M} and
\begin{align*}
	\Phi(x) = \frac{1}{8 \pi} \left(\frac 1 {|x|} + \frac{x \otimes x}{|x|^3} \right),
\end{align*}
is the fundamental solution of the Stokes equations.
Moreover, in \eqref{def:L_n,app} and later on, we use the convention that for a   matrix $M \in \R^{3\times3}$ we write
\begin{align}
  (M : \nabla \Phi)_\alpha = \sum_{\beta, \gamma}  \partial_{x_\beta}M_{\gamma \beta} \Phi_{\alpha \gamma}.
\end{align}

From the definition, we seen that $v= {L}_{n,app}((\xi_1,\dots, \xi_n))T$ is a distributional solution to
\begin{equation} \label{eq:distributional.v}
-\Delta v +\nabla p =  -\dv\left( \underset{i}{\sum} ([T_i]_M+\mS(\xi_i)T_i)\delta_{x_i} \right), \quad \dv v=0, \: \text{ in } \mathbb{R}^3.
\end{equation}

We emphasize that in contrast to $L_n$ the operator $ L_{n,app}$ does only depend on the particle positions but not on the scaling factor $r$. The only dependence on the particle shape and orientations is through the function $\mS(\xi_i)$.
Due to the asymptotic behavior of $\Phi$, $L_{n,app} T$ fails to be in $L^p$ for $p< 3/2$. To capture the decay at infinity we therefore consider $L^{p,2}_{\mathfrak s,w,B(0,M)}$ (see \eqref{def:L^p,2_w,K}) with weight $w$ as in
\eqref{def:weight.function} and for 
\begin{align} \label{def:M}
    M := \sup_n \max\{1,2\, \, \underset{i}{\max}|x_i|,8r\},
\end{align}
which is finite thanks to assumption  \eqref{ass:uniform_bound}.
We make the following observation regarding $L_{n,app}$.
\begin{lem} \label{lem:L.tilde} 
	Let $1<p<3/2$, $M$ be as in \eqref{def:M} and $w$ as in \eqref{def:weight.function}. Then,  the operator $ L_{n,\app}$ satisfies $ L_{n,\app} \in C^1((\S^2)^n;\L(\R^{3n}, L^{p,2}_{\mathfrak s,w,B(0,M)}))$  with
	\begin{align}
		\| L_{n,\app} \|_{C^1((\S^2)^n;\rh{(L^{p,2}_{w,B(0,M)})^{3n})}} \leq C \sqrt{n}.
	\end{align}
\end{lem}
\begin{rem}
    Note that we identified $\L(\R^{3n}, L^{p,2}_{\mathfrak s,w,B(0,M)})$ with $ (L^{p,2}_{\mathfrak s,w,B(0,M)})^{3n}$ in order to clarify that we consider the norm
\begin{align}
    \|L_{n,app}\|^2_{C^1((\S^2)^n;(L^{p,2}_{w,B(0,M)})^{3n})} = \sum_{i,\alpha} \| L_{n,app} e_{i,\alpha}\|^2_{L^{p,2}_{w,B(0,M)}} + \sum_{i,j,\alpha} \|\nabla_{\xi_j} L_{n,app} e_{i,\alpha}\|^2_{L^{p,2}_{w,B(0,M)}}.
\end{align}    
\end{rem}

\begin{proof}
Since $\mathcal S$ is a smooth function of $\xi$, the assertion follows immediately from the decay of $\nabla \Phi$ which implies $\nabla \Phi(\cdot - x_i) \in L_{\loc}^{\frac 3 2_-}(\R^3)$ and $\nabla \Phi(\cdot - x_i) \in L^2_w(\R^3\setminus B(0,M))$.
\end{proof}

The main result of this section is the following proposition.
\begin{prop}\label{pro:L_n}
	Let $M$ be as in \eqref{def:M} and $w$ as in \eqref{def:weight.function}. Then, there exists $N \in \N$ such that for all $n \geq N$ 
 and all $1< p<3/2$  the operator $ L_{n}$ defined in Subsection \ref{sec:well-posedness} satisfies $ L_{n} \in C^1((\S^2)^n;\L(\R^{3n}, L^{p,2}_{\mathfrak s,w,B(0,M)}))$
		\begin{equation*}
		\|L_n -  L_{n,\app}\|_{C^1((\S^2)^n;\rh{(L^{p,2}_{w,B(0,M)})^{3n})}}  \leq C \sqrt{n} (\phi_n \log n+r^{-2+\frac{3}{p}}) .
	\end{equation*}
\end{prop}
The proof of this  proposition is given at the end of this section, see subsection \ref{section:proof_main_prop_L_n}. In subsection \ref{section:proof_main_prop_L_n} we also show how this implies Theorem \ref{th:well-posed}.

\subsection{An intermediate semi-explicit approximation}

For the proof of Proposition \ref{pro:L_n}, we introduce 
yet another approximation for $L_n$ denoted $ L_{n,\app}^{im}$.
For this approximation $L_{n,\app}^{im}$, we neglect interactions between the particles, but we treat each particle by solving a Stokes problem in the exterior domain of that particles. In this sense, $L_{n,\app}^{im}$ can be seen as intermediate between $L_n$ and $ L_{n,\app}$.

More precisely, we define $L_{n, \app}^{im} \colon \R^{3 n} \to  \L(\R^{3 n }, \dot H^1_{\mathfrak s}(\R^3))$ by
\begin{align} \label{eq:L_n,app.representation}
	L_{n,\app}^{im}(\xi_1,\dots,\xi_n) T = \sum_i  U_i[T_i],
\end{align}
where $U_i[T_i]$ is defined as the solution to  
\begin{equation}  \label{eq:U_i.single.rod}
\left \{
\begin{array}{rcl}
- \Delta w_i + \nabla p_i &=& 0 \quad \text{ in }\mathbb{R}^3 \setminus  \B_i,\\
\dv w_i &=&0 \quad \text{ in }\mathbb{R}^3 \setminus   \B_i,\\
D w_i &=& 0  \quad \text{ in } \B_i, \\
\displaystyle \int_{\partial  \B_i} \Sigma(w_i,p_i) \nu &=& 0, \\
\displaystyle \int_{\partial \B_i} [\Sigma(w_i,p_i) \nu ] \times (x - x_i) &=&  T_i.
\end{array}
\right.
\end{equation}

The dependence on the particle orientation of $L_{n,\app}^{im}$ can be made explicit.
Indeed, consider the problem
\begin{equation}  \label{eq:single.rod}
\left \{
\begin{array}{rcl}
- \Delta w + \nabla p &=& 0, \quad \text{ in }\mathbb{R}^3 \setminus r \B,\\
\dv w &=&0 \quad \text{ in }\mathbb{R}^3 \setminus  r \B,\\
D w &=& 0 \quad \text{ in } r \B, \\
\displaystyle \int_{\partial r \B} \Sigma(w,p) \nu &=& 0, \\
\displaystyle \int_{\partial r \B} [\Sigma(w,p) \nu ] \times x &=&  T,
\end{array}
\right.
\end{equation}
where $T \in \R^3$ is a given torque.
We denote by
$$
	(U[T],P[T])
$$
the unique solution $(w,p) \in \dot H^1_{\mathfrak s}(\R^3) \times L^2(\R^3)$ to \eqref{eq:single.rod}.
Then, for $R_i \in SO(3)$ such that $R_i e_3 = \xi_i$, the solution $(U_i[T],P_i[T])$ to \eqref{eq:U_i.single.rod} is given by
\begin{equation}\label{def:U_i}
	\left(U_i[T](x),P_i[T](x)\right)    
	= \left(R_i U[R_i^T T](R_i^T(x-x_i)),  P[R_i^T T](R_i^T(x-x_i)) \right).
\end{equation}
Indeed, 
\begin{align*}
	& \int_{\partial \B_i} R_i \left[\Sigma\left(U[R_i^T T], P[R_i^T T] \right)(R_i^T(x-x_i)) R_i^T \nu \right]\times (x - x_i) \\
	&= \int_{\partial r \B} R_i \left[\Sigma\left(U[R_i^T T], P[R_i^T T] \right)(y)  \nu \right] \times (R_i y) \\
 &= T,
\end{align*}
where we used that the normal $\nu$ also gets rotated and that $(R_i a) \times (R_i b) = R_i (a \times b)$.

\medskip

We show the following estimates between the explicit approximation $L_{n,app}$ and the semi-explicit approximation $L_{n,app}^{im}$.
\begin{prop} \label{pro:L_n,app}
Let $i\in \{1, \cdots,n \}$ and $T_i \in \R^3$.
Then, for all $1 < p < 3/2$, $U_i[T_i] \in C^1(\S^2;\L(\R^{3n}, L^{p,2}_{\mathfrak s,w,B(0,M)}))$
where $w$ and $M$ are as in \eqref{def:weight.function} and \eqref{def:M}, respectively,
 and for all $x \in \R^3 \setminus \B_i$ and all $l \in \N$
\begin{align}  \label{est:nabla.U_i.pointwise}
	|\nabla^l (U_i[T_i] )(x)|  + |\nabla^l  \nabla_{\xi_i} (U_i[T_i] )(x)| \leq \frac {C|T_i| }{|x - x_i|^{l+2}},
	\end{align}
and 
\begin{align} \label{est:nabla.U_i.L^p}
	\|U_i[T_i] \|_{C^1(\S^2;\L(\R^{3n}, L^{p,2}_{\mathfrak s,w,B(0,M)}))} \leq C |T_i|.
\end{align}

In particular,  $L^{im}_{n,\app} \in C^1((\S^2)^n;\L(\R^{3n},L^{p,2}_{\mathfrak s,w,B(0,M)}))$ and 
\rh{
\begin{align} \label{est:L.n.app.L.n.app.im}
		\| L_{n,app}^{im}-{L}_{n,app}\|_{C^1((\S^2)^n;\rh{(L^{p,2}_{w,B(0,M)})^{3n})}} \leq C \sqrt{n}r^{-2 + \frac 3 p} .
	\end{align}
}
\end{prop}

The proof relies on the following expansion of $U[T]$ (see  e.g. \cite[Proposition~2.2.]{HillairetWu19} for a similar result).
\begin{prop}\label{prop_HW} 
There exists $\mathcal{H}[T]$ such that for and  all $x \in \R^3 \setminus r\B$
$$
U[T](x)=-([T]_M + \mathcal S\rhnew{(e_3)} T) : \nabla \Phi(x)+\mathcal{H}[T],
$$
The error term satisfies for all $l\in\N $   and  all $x \in \R^3 \setminus r\B$
\begin{align} \label{est:H}
\left| \nabla^l \mathcal{H}[T](x) \right| \leq C \frac{r|T|}{|x|^{3+l}}, 
\end{align}
where the constant $C$ depends on $l$ and the reference particle $\mathcal{B} $. 
\end{prop}

\rh{
\begin{proof}
    By scaling, it suffices to show the assertion for $r = 1$. Let $|x| > 4$. Then, using the force free condition for  $U[T]$  as well as that $T$ and $\mathcal{S}(e_3) T$ are by definition the torque and stresslet associated to $U[T]$
    \begin{align}
        &U[T](x) + ([T]_M + \mathcal{S}(e_3) T) : \nabla \Phi(x) \\
        &= -\int_{\partial  \B} \Sigma[U[T],P[T]](y) \nu \cdot  \left( \Phi(x-y) - \Phi(x) + y \cdot \nabla \Phi(x) \right)  \dd y \\
        &\leq 2 \|D U[T]\|_{L^2(\R^3)} \|D \psi\|_{L^2(\R^3)} 
    \end{align}
   for any divergence free function $\psi \in \dot H^1(\R^3)$ with $\psi(y) = \Phi(x-y) - \Phi(x) + y \cdot \nabla \Phi(x) $ on $\partial \B$.
By Lemma \ref{lem:extension} below and the decay of $\nabla ^2\Phi$, we find such a $\psi$ with
\begin{align}
    \|D \psi\|_{L^2(\R^3)}  \leq C \frac{1}{|x|^3}.
\end{align}
   Moreover, by some integration by parts, we have
   \begin{align}
        2 \|D U[T]\|^2_{L^2(\R^3)} = T \cdot \mR_2^{-1} T \leq C .
   \end{align}
Collecting these estimates yields the assertion for $|x| > 4$.
On $\partial \B$, \eqref{est:H} holds as well since $U[T](x) = \omega \times x = (\mR_2^{-1} T) \times  x$ on $\partial \B$.
Thus, by standard regularity theory, \eqref{est:H} also holds in $B(0,4) \setminus \B$.
\end{proof}
}

The following Lemma that we used above is standard. For spherical particles it can be found for example in \cite[Lemma 4.4]{NiethammerSchubert19}, and the proof given there also applies for general smooth particles considered here.
\begin{lem} \label{lem:extension}
	Let $\varphi \in H^1_{\mathfrak s}(\B)$. Then, there exists  $\psi \in H^1_{\mathfrak s, 0}(B(0,2))$
	such that $D \psi = D \varphi$ in $\B$ and
	\begin{align}
		\|D \psi\|_{L^2(B(0,2))} \leq C \|D \varphi \|_{L^2(\B)},
	\end{align}
	where $C$ depends only on $\B$.
\end{lem}
\begin{rem}
	By translation  and scaling, the same statement holds for $\B$ replaced by $\B_i$, where the constant $C$ is independent of $i$. 
\end{rem}




\begin{proof}[Proof of Proposition \ref{pro:L_n,app}]
\rh{
We focus on the estimates of derivatives in $\xi_i$ in \eqref{est:nabla.U_i.pointwise}, \eqref{est:nabla.U_i.L^p} and  \eqref{est:L.n.app.L.n.app.im}. The estimates for the functions themselves can be obtained analogously.
The right-hand side of the representation \eqref{def:U_i} allows to view $U_i$ as dependent of $R_i$ instead of $\xi_i$. We recall that this is true since the right-hand side is the same for all $R_i \in SO(3)$ with $R_i e_3 = \xi_i$.  It is then sufficient to consider the derivative with respect to $R_i$.\footnote{To see this, observe that we might locally fix the choice of $R_i$ such that $R_i[\xi_i]$ is differentiable in $\xi$ with $|\nabla_{\xi_i} R| \leq C$.}
}
Fix $ 1 \leq i \leq n$  \am{and let $T_i \in \R^3$}. \rh{It suffices to consider the case $x_i = 0$. (Note that the weighted norm $L^2_w(\R^3 \setminus B(0,M))$ is not translation invariant. However, by definition of $M$ in \eqref{def:M} we have $|x - x_i| \geq \frac 1 2 |x|$ for all $x \in \R^3 \setminus B(0,M)$ and thus the position of $x_i$ does not matter.)}

Combining \eqref{def:U_i} and Proposition \ref{prop_HW} and using how $\mathcal S$ and $\Phi$ transform under rotations, we have
\begin{align}
    U_i[T_i](x) + ([T_i]_M + \mathcal S(\xi_i)T_i) : \nabla \Phi(x) = R_i \mathcal H[R_i^T T_i](R_i^T x),
\end{align}
and thus by the chain rule, \eqref{est:H}
\begin{align} \label{eq:nabla.U_i.1}
  \left|\nabla_{\xi_i} \left( ([T_i]_M + \mathcal S(\xi_i)T_i) : \nabla \Phi(x) + U_i[T_i]\right)(x)\right| \leq 
  C \frac{r|T_i|}{|x|^{3}} \quad \text{in } \R^3 \setminus \mathcal B_i.
\end{align}
Using the decay of $\nabla \Phi$, this implies \eqref{est:nabla.U_i.pointwise} for $l=0$. The estimate for $l \geq 1$ is analogous.

Estimate \eqref{eq:nabla.U_i.1} directly implies  $\|\nabla_{\xi_i} U_i[T_i]\|_{L^2_{w}(\R^3 \setminus B(0,M))} \leq C |T_i|$. Moreover,  we have $U_i[T_i](x) = r^{-3} (\mR_2^{-1} T_i) \times x$ in $\B_i$, and thus
\begin{align} \label{eq:nabla.U_i.2}
  \left|\nabla_{\xi_i} \left( ([T_i]_M+ \mathcal S(\xi_i)T_i) : \nabla \Phi(x) + U_i[T_i]\right)(x)\right| \leq 
  C \frac{|T_i|}{|x|^{2}} \quad \text{in }  \mathcal B_i.
\end{align}
This together with \eqref{est:nabla.U_i.pointwise} implies $\|\nabla_{\xi_i} U_i[T_i]\|_{L^p_{\loc}( B(0,M))} \leq C |T_i|$ for $p < 3/2$ and therefore \eqref{est:nabla.U_i.L^p} holds.

It remains to show  \eqref{est:L.n.app.L.n.app.im}. By linearity \rhnew{of $L_{n,app}$ and $L_{n,app}^{im}$ in the particles} it suffices to show
\begin{align} \label{est:L.n.app.im.one.particle}
  \|U_i[T_i] + ([T_i]_M + \mathcal S(\xi_i)T_i) : \nabla \Phi\|_{C^1(\S^2; L^{p,2}_{w,B(0,M)})} \leq  C r^{-2 + 3/p} |T_i|.
\end{align}
However, the pointwise estimates  \eqref{eq:nabla.U_i.1}--\eqref{eq:nabla.U_i.2}
 directly imply \eqref{est:L.n.app.im.one.particle}.
 \end{proof}

\subsection{Estimates for \texorpdfstring{$L_n - L_{n,\app}^{im}$}{Ln} through the method of reflections}

\label{sec:reflections}

\rh{In order to estimate $L_n - L_{n,app}^{im}$, we write $L_n$ in terms of single particle operators relying on the so-called method of reflections similarly as in \cite{Hofer18MeanField, Hoefer19}.
We emphasize that $L_n - L_{n,app}^{im}$ can be also estimated by energy and duality methods (see e.g. \cite{Gerard-VaretHoefer21}). More precisely, arguing as in \cite{Gerard-VaretHoefer21} one can show that
\begin{align} \label{est:L_n-L_n,app.energy}
    \|(L_n - L_{n,app}^{im})\|_{L^p_{\loc}} \leq C \sum_i \|D L_{n,app}^{im} T\|_{L^1(\cup_i \B_i)}.
\end{align}
Using the decay estimate from Proposition \ref{prop_HW} together with assumptions \eqref{ass:phi.log.n},\eqref{ass:well.separated}, the right-hand side is bounded by $\phi_n \log n \sum_i |T_i|$ which we recover through the method of reflections, see Proposition \ref{pro:L_n.L_n,app.l^1}. 

The reason why we rely on the method of reflections in this section are the derivatives with respect to the particle orientations $\xi_i$. Here the method of reflections is very useful because the single particle problems occurring in the method of reflections have an explicit dependence on $\xi_i$, giving direct access to these shape derivatives.
}
More precisely, following the notation from \cite{Hofer18MeanField, Hoefer19}, we introduce $Q_i$ as the solution operator that maps a  function $w \in H^1_{\mathfrak s}(\B_i)$ to the solution $v \in \dot H^1_{\mathfrak s}(\R^3)$ to 
\begin{equation} \label{eq:Q}
\left\{
\begin{array}{ll}
	- \Delta v + \nabla p = 0, \quad \dv v =  0& \quad \text{ in } \R^3 \setminus \B_i, \\
	D v = D w &\quad \text{ in } \ml{\B_i}, \\
\displaystyle	0 = \int_{\partial \B_i} \Sigma[v,p] n = \int_{\partial \B_i} \Sigma[v,p] n \times (x- x_i).
	\end{array}\right.
\end{equation} 

Then, we claim that
\begin{align} \label{eq:MOR}
	L_n = \lim_{k \to \infty} (1 - \sum_i Q_i)^k L_{n,\app}^{im}.
\end{align}
in the sense of convergence of operators $\R^{3n} \to \dot H^1_{\mathfrak s}(\R^3)$.
Indeed, this has been proven  for spherical particles  in \cite{Hofer18MeanField}  under
the assumptions \eqref{ass:phi.log.n},\eqref{ass:well.separated} and in  \cite{Hoefer19} assuming only \eqref{ass:well.separated}.
Since we impose \eqref{ass:phi.log.n} in order to control the right-hand side of \eqref{est:L_n-L_n,app.energy}, we follow here the (simpler) approach of \cite{Hofer18MeanField}.
The adaptation to non-spherical particles is rather straightforward,
and  based on the following decay estimates for $Q_i$.

\begin{lem} \label{lem:decay.Q}
Let $w \in H^1_{\mathfrak s}(\mathcal{B}_i) $ such that $ \nabla w \in L^\infty(\mathcal{B}_i)$. There exists a universal constant $C>0$ such that for all $x \in \R^3 \setminus B(x_i,2r)$
\begin{align} \label{est:Q_i.pointwise}
	|\nabla^l (Q_i w)(x)| \leq \frac {C r^3}{|x - x_i|^{l+2}} \|D w\|_{L^\infty(\B_i)} .
\end{align}
Moreover, 
	\begin{align} \label{est:Q_i.H^1}
		\|Q_i w\|_{\dot H^1(\R^3)} \leq C r^{3/2} \|D w\|_{L^\infty(\B_i)},
	\end{align}
and for all $p< 3/2$, with $w$ and $M$ as in \eqref{def:weight.function} and \eqref{def:M}, repsectively,
\begin{align} \label{est:Q_i.L^p}
 \|Q_i w\|_{L^{p,2}_{w,B(0,M)}} \leq C r^3 \|D w\|_{L^\infty(\B_i)}.
\end{align}
\end{lem}
\begin{proof}
	We observe that $Q_i w$ minimizes the $\dot H^1$ norm among all divergence free functions with $D v = D w$. By Lemma \ref{lem:extension} we thus immediately obtain \eqref{est:Q_i.H^1}.
	
	Now, let $p < 3/2$, $K \subset \R^3$ be compact and let $g \in L^{p'}$ with $\supp g \subset K$. Let $\varphi$ be the solution to
	\begin{align*}
		-\Delta \varphi + \nabla p = g , \quad \dv \varphi = 0 \qquad \text{in } \R^3.
	\end{align*}
	Then, by standard regularity theory and Sobolev embedding
	\begin{align*}
		\|\nabla \varphi\|_\infty \leq C_{K,p} \|g\|_{L^{p'}}.
	\end{align*}
	Consequently, using the equation that $Q_i w$ solves
	\begin{align*}
		\int g \cdot Q_i w &= \int D \varphi : D Q_i w \\
		& = \int_{\B_i} D \varphi : D Q_i w + \int_{\R^3 \setminus \B_i} D \varphi : D Q_i w  \\
		& \leq C  r^{3} \|D w\|_{L^\infty(\B_i)}  \|g\|_{L^{p'}} + 
	C	r^{3/2} \|D w\|_{L^\infty(\B_i)}  \|D \psi\|_{L^2(\B_i)}
	\end{align*}
	for all $\psi \in \dot H^1_{\mathfrak s}(\R^3)$ such that $D \psi = D \varphi$ in $\B_i$.
	Appealing again to Lemma \ref{lem:extension}, such a function $\psi$ exists with
	\begin{align*}
		\|D \psi\|_{\dot H^1(\R^3)} \leq C r^{3/2} \|D \varphi\|_{L^\infty(\B_i)} \leq  C r^{3/2} \|g\|_{L^{p'}}.
	\end{align*}
	Combination of the above estimates with $K=B(0,M)$ yields the $L^p(B(0,M))$ estimate in \eqref{est:Q_i.L^p}.

The proof of \eqref{est:Q_i.pointwise} is similar to the proof of Proposition \ref{prop_HW} in the sense that we have for $|x-x_i|\geq 2r$  
$$
|Q_i w(x)| \leq C \frac{r^{3/2}}{|x-x_i|^2} \| D w \|_{L^2(\mathcal{B}_i)}
$$
In particular this yields the $L^2_w(\R^3\setminus B(0,M))$ estimate in \eqref{est:Q_i.L^p}. Indeed since the above decay is valid for $x \not \in  B(x_i,2r)$, we have in particular $|x-x_i|\geq \frac{1}{2}|x| $ for $x \not \in  B(0,M)$, $1 \leq i \leq n$,  this yields
\begin{align*}
\|Q_i w\|_{L^2_w(\R^3\setminus B(0,M))}&\leq C r^3 \|D w\|_{L^\infty(\mathcal{B}_i)} \left (\int_{\R^3\setminus B(0,M)} \frac{(1+|x|)^{a}}{|x|^4}dx\right)^{1/2}
\end{align*}
this concludes the proof since the above integral is finite for \amnew{$a<1$}.
\end{proof}

\begin{prop} \label{pro:L_n.L_n,app.l^1}
\rhnew{There exists $N_0 \in \N$ such that for all $n \geq N_0$}
	\begin{align} \label{eq:convergence.MOR}
	L_n = \lim_{k \to \infty} (1 - \sum_i Q_i)^k L_{n,\app}^{im}
\end{align}
in the sense of convergence of linear operators operators from $\R^{3n}$ to  $\dot H^1_{\mathfrak s}$.
Moreover, 	with $M$ and $w$ as in \eqref{def:M}  and \eqref{def:weight.function}, respectively, for  all $1\leq p < 3/2$,
\rh{
\begin{equation} \label{est:L_n.L_n,app.l^1}
	\|L_n -  L_{n,\app}^{im}\|_{C((\S^2)^n;(L^{p,2}_{w,B(0,M)})^{3n})} \leq C \sqrt{n} \phi_n \log n .
\end{equation}
}
\end{prop}

\begin{proof}
	The convergence \eqref{eq:convergence.MOR} can be proved exactly as in \cite{Hofer18MeanField}. We therefore only give the proof of \eqref{est:L_n.L_n,app.l^1}.


	Let $p< 3/2$ and $T \in \R^{3n}$. We denote $v = L_n T$  and $v_k =  (1 - \sum_i Q_i)^k L_{n,\app}^{im} T$.
	To prove \eqref{est:L_n.L_n,app.l^1}, we apply Lemma \ref{lem:decay.Q} to see that
	\begin{align} \label{est:v_k+1.v_k_L^p_loc}
	\|v_{k+1} - v_k\|_{L^2_w(\R^3\setminus B(0,M))}+	\|v_{k+1} - v_k\|_{L^p(B(0,M))}& \leq  \sum_i \|Q_i v_k\|_{L^p(B(0,M))} + \|Q_i v_k\|_{L^2_w(\R^3\setminus B(0,M))} \notag\\
	&	\leq  C \sum_i r^3 \|D v_k\|_{L^\infty({\B_i})}. 
	\end{align}
			
	Now, using the fact that $D(Q_i v_k) = D v_k$ in $B_i$ we get $$Dv_{k+1}= D v_k - \underset{j}{\sum} D Q_j v_k=- \underset{j \neq i}{\sum} D Q_j v_k, \text{ in } B_i.  $$ Thus, using again Lemma \ref{lem:decay.Q} and assumption \eqref{ass:well.separated}, we get
	\begin{align*} 
		\sum_i \|D v_{k+1}\|_{L^\infty(\B_i)} &\leq \sum_i \sum_{j \neq i}  \|D(Q_j v_k)\|_{L^\infty(\B_i)} \\
		& \leq \sum_i \sum_{j \neq i} \frac{r^3}{|x_i - x_j|^3}  \|D v_{k}\|_{L^\infty(\B_j)}\\
		& \leq C \phi_n \log n \sum_j \|D v_{k}\|_{L^\infty(\B_j)}.
	\end{align*}
	Hence, by iteration
	\begin{align} \label{est:Dv_k}
		\sum_i \|D v_{k}\|_{L^\infty(\B_i)} \leq (C \phi_n \log n)^k \sum_i\|D v_0\|_{L^\infty( \B_i)}.
	\end{align}
	By \eqref{eq:L_n,app.representation}, we can write $v_0 = \sum_i U_i(T_i)$, and thus the decay estimates from Proposition \ref{prop_HW} yield
	\begin{align} \label{est:Dv_0}
	r^3	\sum_i \|D v_0\|_{L^\infty(\B_i)} \leq C  \sum_i \sum_{j \neq i} \frac{|T_j| r^3}{|x_i - x_j|^3} \leq \phi_n \log n |T|_{l^1}.
	\end{align}
	Combining \eqref{est:v_k+1.v_k_L^p_loc}, \eqref{est:Dv_k} and \eqref{est:Dv_0} yields
	\begin{align*}
		\|v_{k+1} - v_k\|_{L^p_{\loc}} + \|v_{k+1} - v_k\|_{L^2_w(\R^3\setminus B(0,M))}\leq (C \phi_n \log n)^{k+1} |T|_{l^1}.
	\end{align*}
	Summing these errors in $k$ yields \eqref{est:L_n.L_n,app.l^1}. 
\end{proof}

We now turn to the derivatives in $\xi_i$. We begin by estimating the derivative of $Q_i$ with respect to $\xi_i$.

\begin{lem} \label{lem:decay.nabla.Q}
Let $w \in H^2_{\mathfrak s}(B(x_i,2r))$ such that $ \nabla w \in W^{1,\infty}(B(x_i,2r)) $.
Then,  $\nabla_{\xi_i} (Q_i w)\in L^2_\loc(\R^3)$ and for all $x \in \R^3 \setminus B(x_i,2r)$ and all $l \in \N$
\begin{align}  \label{est:nabla.Q_i.pointwise}
	|\nabla^l  \nabla_{\xi_i} (Q_i w)(x)| \leq \frac {C r^3}{|x - x_i|^{l+2}} \left(\|D w\|_{L^\infty(\B_i)} + r \|\nabla D w\|_{L^\infty(\B_i)} \right).
	\end{align}
Moreover, for all $p< 3/2$ 	and with $M$ and $w$ as in \eqref{def:M}  and \eqref{def:weight.function}, respectively,
\begin{align} \label{est:nabla.Q_i.L^p}
\|\nabla_{\xi_i} Q_i w\|_{L^{p,2}_{w,B(0,M)}} \leq C r^3 \left(\|D w\|_{L^\infty(\B_i)} + r \|\nabla D w\|_{L^\infty(\B_i)} \right).
\end{align}
\end{lem}

\begin{proof}
	We begin by dropping all indices $i$ and assume $x_i = 0$.
	Moreover, we set $Q = Q[\xi]$ to denote the dependence on $\xi$.
	By considering the defining equation for $Q$, \eqref{eq:Q}, we observe that for any $R \in SO(3)$ with $R \xi = e_3$
	\begin{align} \label{eq:Q[xi]}
		(Q[\xi] w)(x) = R (Q[e_3] \bar w) (R^T x),
	\end{align}
	where $\bar w(x) = R^T w(R x)$. Note that this corresponds to the way how $U_i$ is obtained from $U$ in \eqref{def:U_i}. Analogously, as argued at the beginning of the proof of Proposition \ref{pro:L_n,app}, it suffices to view $Q$ as a function of $R$ to derive estimates for the derivative.
	
By the assumptions on $w$ and the chain rule
$\nabla_\xi \bar w \in H^1(B(0,2r))$ with
\begin{align} \label{est:bar.w}
	|\nabla_\xi D_x \bar w|(x) \leq C |\nabla w|(x) +  C |x||\nabla^2 w|(x).
\end{align}
Consequently, since $w \in H^2(B(0,2r))$ and $Q$ is a linear operator with values in $\dot{H}^1$,
the representation \eqref{eq:Q[xi]} implies that $\nabla_\xi (Q[\xi] w) \in L^2_{\loc}(\R^3)$.
Moreover, we can combine \eqref{eq:Q[xi]} and \eqref{est:bar.w} with Lemma \ref{lem:decay.Q} to obtain for $x \in \R^3 \setminus B(0,2r)$ 
\begin{align} \label{eq:expansion.nabla.Q}
	\begin{aligned}
		|\nabla_{\xi} (Q w)(x)| &\leq C |\nabla^l(Q[e_3] \bar w) (R^T x)|
		+ C |x| |\nabla(Q[e_3] \bar w) (R^T x)|
		+  C |Q[e_3]( \nabla_\xi \bar w) (R^T x)| \\
		& \leq C \frac{r^3}{|x|^{2}} \|D w\|_{L^\infty(\B_i)} + C \frac{r^4}{|x|^{2}} \|\nabla D w\|_{L^\infty(\B_i)}.
	\end{aligned}
\end{align}
This establishes \eqref{est:nabla.Q_i.pointwise} for $l=0$ and the $L^2_w(\R^3\setminus B(0,M))$ estimate in \eqref{est:nabla.Q_i.L^p}. The estimate for $l \geq 1$ is analogous.

Moreover, 
\eqref{eq:expansion.nabla.Q} implies
\begin{align}
	\|\nabla_{\xi} (Q w)(x)\|_{L^p(B(0,M)\setminus B(0,2r)} \leq C r^3 \|D w\|_{L^\infty(\B_i)} + C r^4 \|\nabla D w\|_{L^\infty(\B_i)}.
\end{align}

It remains to estimate the $L^p$ norm inside of $B(0,2r)$.
We note that the first inequality in \eqref{eq:expansion.nabla.Q} also holds for all $x \in B(0,2r)$. 
Together with Lemma \ref{lem:decay.Q} and \eqref{est:bar.w}, it implies for all $p < 3/2$
\begin{align*}
	\|\nabla_{\xi} (Q w)\|_{L^p( B(0,2r))} &\leq C r^3 \|D \bar w\|_{L^\infty (B(0,2r))} +  \|x\|_{L^6(B(0,2r))} \|\nabla Q w\|_{L^2(\R^3)} + r^3 \|\nabla_\xi D_x \bar w\|_{L^\infty(\B_i)} \\
	& \leq C r^3 \|D w\|_{L^\infty(\B_i)} + C r^4 \|\nabla D w\|_{L^\infty(\B_i)}.
\end{align*}
This finishes the proof.
\end{proof}

To estimate $\nabla_{\xi_i} (L_n - L_{n,\app}^{im})$ the idea is to proceed similarly as in Proposition \ref{pro:L_n.L_n,app.l^1}. However, this is more delicate due to the loss of regularity when taking the derivative in $\xi_i$. More precisely, by Lemma \ref{lem:decay.nabla.Q} we only have that $\nabla_{\xi_i} (Q_i w)$ lies in $L^2_\loc$. In fact it is easy to see that $\nabla_{\xi_i} (Q_i w)$ will in general not possess a weak derivative regardless how much regularity we impose on $w$.

A consequence of this loss of regularity is that $Q_i \nabla_{\xi} Q_j$ is not defined a priori since $Q_i$ needs data in $H^1(\B_i)$.
However Lemma \ref{lem:decay.nabla.Q} ensures that $Q_i \nabla_{\xi} Q_j$ is well defined for $i \neq j$.
On the other hand, one can exploit that $Q_i Q_i = Q_i$ to avoid any appearance
of $Q_i \nabla_{\xi} Q_i$. To this end, we rewrite the series expansion \eqref{eq:MOR}.
Namely, adopting the notation $v_k =  (1 - \sum_i Q_i)^k L_{n,\app}^{im} T$, we have
\begin{multline} \label{eq:MOR.expansion}
	v_k = v_0 - \sum_{i_1} Q_{i_1} v_0 + \sum_{i_1} \sum_{i_2 \neq i_1} Q_{i_2} Q_{i_1} v_0 - \sum_{i_1} \sum_{i_2 \neq i_1} \sum_{i_3 \neq i_2}  Q_{i_3} Q_{i_2} Q_{i_1} v_0 + \dots \\
	+ 
	(-1)^k \sum_{i_1} \sum_{i_2 \neq i_1} \dots \sum_{i_k \neq i_{k-1}}  Q_{i_k} Q_{i_{k-1}} \dots Q_{i_1} v_0.
\end{multline}
This representation can be directly deduced from \eqref{eq:MOR} by just using
$Q_i Q_i = Q_i$ (see \cite[Section 2]{HoferVelazquez18} for details).

\begin{prop}\label{pro:nabla.L_n.L_n,app}
\rhnew{There exists $N_0 \in \N$ such that for all $n \geq N_0$ the following holds.}
\rh{Let $M$ be as in \eqref{def:M} and $w$ as in \eqref{def:weight.function},  $1 \leq i \leq n$ and $1 \leq \alpha \leq 3$ Then, for all $1<p < 3/2$, $L_n e_{i,\alpha}$ is differentiable in $\xi_j$ for all $1 \leq j \leq n$ as a function in $L^{p,2}_{\mathfrak s,w,B(0,M)}$  and  we have 
\begin{align} \label{est:nabla.L_n.L_n,app}
	\|\nabla_{\xi_i} (L_n - L_{n,\app}^{im}) e_{i,\alpha}\|_{L^{p,2}_{w,B(0,M)}} \leq C \phi_n \log n .
\end{align}
Moreover, for $j\neq i$  
\begin{align} \label{est:nabla.L_n.L_n,app.different}
	\|\nabla_{\xi_j} (L_n - L_{n,\app}^{im}) e_{i,\alpha}\|_{L^{p,2}_{w,B(0,M)}} \leq C \frac{r^3 }{|x_i-x_j|^3}.
\end{align}
}
\end{prop}
\begin{proof}
In the following we will assume $i=1$. Since the $L^2_w(\R^3 \setminus B(0,M)) $ and $L^p(B(0,M))$ estimates are analogous we only treat below the $L^p(B(0,M))$ estimates.
	
let $v_k =  (1 - \sum_i Q_i)^k L_{n,\app}^{im} e_{i,\alpha}$, then by virtue of \eqref{eq:MOR.expansion} we have
	\begin{align*}
		v_{k+1}  - v_{k} = (-1)^{k+1} \sum_{i_1 \neq 1} \sum_{i_2 \neq i_1} \dots \sum_{i_{k+1} \neq i_{k}}  Q_{i_{k+1}} Q_{i_{k}} \dots Q_{i_1} v_0.
	\end{align*}
	Here we used that $D v_0 = D   L_{n,\app}^{im} e_{i,\alpha} = D U_1[e_{\alpha}] = 0$ in $\B_1$
	to deduce that the first sum only runs over $i_1 \neq 1$.  
	
	Thanks to $\nabla_{\xi_1} Q_i = 0$ for $i \neq 1$, taking the derivative in $\xi_1$ yields
	\begin{align} \label{eq:MOR.expansion.nabla}
	\begin{aligned}
		&\nabla_{\xi_1} (v_{k+1}  - v_{k}) \\
		&= \sum_{i_1 \neq 1} \sum_{i_2 \neq i_1} \dots \sum_{i_{k+1} \neq i_{k}}  Q_{i_{k+1}} Q_{i_{k}} \dots Q_{i_1} \nabla_{\xi_1} v_0 \\
&+		\sum_{l =2}^{k+1} (-1)^{k+1} \sum_{i_1 \neq 1} \sum_{i_2 \neq i_1} \dots
		 \sum_{\substack{i_{l -1} \neq 1 \\ i_{l-1} \neq i_{l-2}}}  \sum_{i_{l+1} \neq 1} \dots \sum_{i_{k+1} \neq i_{k}}  Q_{i_{k+1}} \dots Q_{i_{l+1}} \nabla_{\xi_1} Q_1  Q_{i_{l-1}} \dots Q_{i_1} v_0 \\
		 &=:  \Psi_{k,1} + \Psi_{k,2}.
	\end{aligned}
	\end{align}
	Although the right-hand side looks very complicated, it can be estimated analogously  as in the proof of Proposition \ref{pro:L_n.L_n,app.l^1}.
	Indeed, inductive application of Lemma \ref{lem:decay.Q} and eventually application of Proposition \ref{pro:L_n,app} yields
	\begin{align*}
		\|\Psi_{k,1}\|_{L^p_\loc} & \leq   C r^3 \sum_{i_1 \neq 1} \sum_{i_2 \neq i_1} \dots \sum_{i_{k+1} \neq i_{k}} \|D  Q_{i_{k}} \dots Q_{i_1} \nabla_{\xi_1} v_0\|_{L^\infty(\B_{k+1})} \\
				& \leq  C r^3 \sum_{i_1 \neq 1} \sum_{i_2 \neq i_1} \dots \sum_{i_{k+1} \neq i_{k}} \frac{r^3}{|x_{i_k} - x_{i_{k+1}}|^3} \| D  Q_{i_{k-1}} \dots Q_{i_1} \nabla_{\xi_1} v_0\|_{L^\infty(\B_{k})} \\
		& \leq C r^3 (C \phi_n \log n)^k \sum_{i_1 \neq 1} \|D \nabla_{\xi_1} v_0\|_{L^\infty(\B_{i_1})} \\
		& \leq C r^3 (C \phi_n \log n)^k \sum_{i_1 \neq 1}\frac{1}{|x_1 - x_{i_1}|^3} \\ 
		& \leq (C \phi_n \log n)^{k+1}.
	\end{align*}
	For the second term on the right-hand side of \eqref{eq:MOR.expansion.nabla},
	we proceed similarly. We observe that the combination of Lemmas \ref{lem:decay.Q} and \ref{lem:decay.nabla.Q} implies for $i,j \neq 1$ and any function $\psi \in H^1_{\mathfrak s}(\B_i)$ such that $\nabla \psi \in L^\infty(\mathcal{B}_i)$.

	\begin{align*}
		\|D \nabla_{\xi_1} Q_1 Q_i \psi \|_{L^\infty(\B_j)} &\leq C \frac{r^3}{|x_j - x_1|^3}  \left( \|D Q_i \psi\|_{L^\infty(\B_1)} + r \|D\nabla  Q_i \psi\|_{L^\infty(\B_1)} \right) \\
		& \leq  C \frac{r^3}{|x_j - x_1|^3} \frac{r^3}{|x_i - x_1|^3} \left(1 + \frac{r}{|x_i - x_1|} \right) \|D \psi\|_{L^\infty(\B_i)} \\
		& \leq C  \frac{r^3}{|x_j - x_1|^3} \| D\psi\|_{L^\infty(\B_i)}.
	\end{align*}
	Similarly in the case $l=k+1$ we have
	\begin{align*}
		\| \nabla_{\xi_1} Q_1 Q_i \psi \|_{L^p_{loc}} \leq  C  \frac{r^3}{|x_i - x_1|^3} \| D\psi\|_{L^\infty(\B_i)}.
	\end{align*}
	In this way, we can estimate the second term on the right-hand side of \eqref{eq:MOR.expansion.nabla} by 
	\begin{align*}
	\|\Psi_{k,2}\|_{L^p_\loc} & \leq  k (C \phi_n \log n)^{k+1}
	\end{align*}
	where the factor $k$ originates from the sum over $l$. Combining these estimates for $\Psi_{k,1}$ and $\Psi_{k,2}$ yields
	\begin{align*}
		\|\nabla_{\xi_1} (v_{k+1}  - v_{k})\|_{L^p_\loc} \leq  (k+1) (C \phi_n \log n)^{k+1}.
	\end{align*}
	Summation over $k$ yields the first assertion.
	
For the second estimate, we set $j=2$ and remark that since $\nabla_{\xi_2} v_0=\nabla_{\xi_2}U_1[e_\alpha]=0 $ we have
\begin{align*}
&\nabla_{\xi_2}(v_{k+1}-v_k) \\
&=\sum_{l =1}^{k+1} (-1)^{k+1} \sum_{i_1 \neq 1} \sum_{i_2 \neq i_1} \dots  \sum_{\substack{i_{l -1} \neq 2 \\ i_{l-1} \neq i_{l-2}}}  \sum_{i_{l+1} \neq 2} \dots \sum_{i_{k+1} \neq i_{k}}  Q_{i_{k+1}} \dots Q_{i_{l+1}} \nabla_{\xi_2} Q_2  Q_{i_{l-1}} \dots Q_{i_1} v_0 
\end{align*}	
with the convention that for $l=1$ the term corresponds to
$$
 \sum_{i_2 \neq 2} \dots   \sum_{i_{k+1} \neq i_{k}}  Q_{i_{k+1}}   \dots Q_{i_2} \nabla_{\xi_2} Q_2 v_0.
$$
We have then
\begin{align}
\begin{aligned}
&\| \nabla_{\xi_2}(v_{k+1}-v_k)\|_{L^p_\loc} \leq r^3 \sum_{l =1}^{k+1} \sum_{i_1 \neq 1}  \dots  \sum_{\substack{i_{l -1} \neq 2 \\ i_{l-1} \neq i_{l-2}}}  \sum_{i_{l+1} \neq 2}\sum_{i_{l+2} \neq i_{l+1}} \dots \sum_{i_{k+1} \neq i_{k}}\\
&  \frac{Cr^3}{|x_{i_{k+1}}-x_{i_k}|^3}\cdots \frac{Cr^3}{|x_{i_{l+2}}-x_{i_{l+1}}|^3}\frac{Cr^3}{|x_{i_{l+1}}-x_2|^3 } \frac{Cr^3}{|x_{2}-x_{i_{l-1}}|^3}\cdots \frac{C }{|x_{i_1}-x_1 |^3}\\
& \leq r^3\underset{l=1}{\overset{k+1}{\sum}} (C\phi_n \log n)^{k-l+1} \sum_{i_1 \neq 1}  \dots  \sum_{\substack{i_{l -1} \neq 2 \\ i_{l-1} \neq i_{l-2}}}  \frac{Cr^3}{|x_{2}-x_{i_{l-1}}|^3}\cdots \frac{C}{|x_{i_1}-x_1 |^3}\\
&:= \underset{l=1}{\overset{k+1}{\sum}} (C\phi_n \log n)^{k-l+1} I_{l-1} ,
\end{aligned} \label{est:nabla.xi_2}
\end{align}	
where for $l\geq 1$
$$
I_l:= \underset{i_1\neq 1}{\sum} \dots \sum_{\substack{i_{l} \neq 2 \\ i_{l} \neq i_{l-1}}} \frac{Cr^3}{|x_{2}-x_{i_l}|^3} \frac{Cr^3}{|x_{i_l}-x_{i_{l-1}}|^3}\dots \frac{Cr^3}{|x_{i_1}-x_1|^3},
$$
and $I_0:=\frac{C r^3}{|x_2-x_1 |^3}$. 
Now  we aim to show by induction that for some constant $\bar C>C$ we have for $l\geq 1$
\begin{align} \label{est:I_k}
I_l  \leq ( \bar C \phi_n \log n)^{l} \frac{Cr^3}{|x_1-x_2|^3},
\end{align}
indeed, for $l>1$, we have by separating the term with $i_{l-1}=2 $
\begin{align*}
I_l&= \underset{i_1\neq 1}{\sum} \dots  \sum_{i_{l-1} \neq i_l, 2} \,\sum_{i_{l} \neq 2 }\frac{Cr^3}{|x_{2}-x_{i_l}|^3}\frac{Cr^3}{|x_{i_l}-x_{i_{l-1}}|^3} \frac{Cr^3}{|x_{i_{l-1}}-x_{i_{l-2}}|^3}\dots \frac{Cr^3}{|x_{i_1}-x_1|^3}\\
&+  \underset{i_1\neq 1}{\sum} \dots \sum_{i_{l-2} \neq 2 }\, \sum_{i_{l} \neq 2} \frac{Cr^3}{|x_{2}-x_{i_l}|^3} \frac{Cr^3}{|x_{i_l}-x_{2}|^3} \frac{Cr^3}{|x_{2}-x_{i_{l-2}}|^3}\dots \frac{Cr^3}{|x_{i_1}-x_1|^3}\\
&\leq 8 (C\phi_n \log n ) I_{l-1}+  (C \phi_n \log n)^2I_{l-2},
\end{align*}
where we used $\log n \geq 1$ (for $n \geq 3$) for the second term and for the first term that for any $i_l\neq 2 \neq i_{l-1}$ we have
$$
\frac{1}{|x_2-x_{i_l}|} \frac{1}{|x_{i_l}-x_{i_{l-1}}|} \leq \frac{1}{|x_2-x_{i_{l-1}}|} \left(\frac{1}{|x_2-x_{i_l}|} +\frac{1}{|x_{i_l}-x_{i_{l-1}}|}  \right).
$$
This yields \eqref{est:I_k} using the induction hypothesis and taking $ 8C \bar C+C^2\leq \bar C^2$. Moreover, using the same arguments as above, one can show that the induction hypothesis is satisfied for $l=0,1$.
	
Inserting \eqref{est:I_k} into \eqref{est:nabla.xi_2} and summing over $k$ yields \eqref{est:nabla.L_n.L_n,app.different}.
\end{proof}

\subsection{Proofs of Proposition \ref{pro:L_n} and Theorem \ref{th:well-posed}}\label{section:proof_main_prop_L_n}
\begin{proof}[Proof of Proposition \ref{pro:L_n}]
The proof of the Proposition is a direct consequence of Propositions \ref{pro:L_n,app}, \ref{pro:L_n.L_n,app.l^1} and \ref{pro:nabla.L_n.L_n,app}. Indeed we have thanks to Propositions \ref{pro:L_n,app} and \ref{pro:L_n.L_n,app.l^1} for all $1<p<3/2$
\begin{align*}
&\left \| L_n-L_{n,app}\right\|_{(L^\infty((\mathbb{S}^2)^n; (L^{p,2}_{w,B(0,M)})^{3n})} \\
&\leq C \rh{\sqrt n} (\phi_n \log n + r^{\frac 3 p - 2}).
\end{align*}
For the derivative we get from Proposition \ref{pro:nabla.L_n.L_n,app} for all $1<p<3/2$
\begin{align*}
&\underset{i}{\sum}\underset{j,\alpha}{\sum} \left \|\nabla_{\xi_i} (L_n-L_{n,app}) e_{j,\alpha}\right\|^2_{L^\infty((\mathbb{S}^2)^n, L^{p,2}_{w,B(0,M)})} \\
&\leq C \underset{i}{\sum}  \underset{\alpha}{\sum}\left(\sum_{j\neq i} \frac{r^6}{|x_i-x_j|^6}  + (\phi_n \log n)^2 \right)  \\
&\leq C n \phi_n^2(1+\log^2 n).
\end{align*}
where we used that for all $i$,
$ \underset{j \neq i}{\sum} \frac{r^6}{|x_i-x_j|^6} \leq C \frac{r^6}{d_{\min}^6} \leq C \phi_n^2$ thanks to assumption \ref{ass:well.separated}.

Combining again with Proposition  \ref{pro:L_n,app} yields the assertion.
\end{proof}

For the proof of Theorem \ref{th:well-posed}, we first deduce the following corollary which is a direct consequence of Lemma \ref{lem:L.tilde} and Proposition \ref{pro:L_n} combined with Lemma \ref{lem:embedding.Lebesque.weighted}.
\begin{cor} \label{cor:L_n.H^-s}
For all $n$ sufficiently large and for all $1/2<s<1$ and $p=\frac{6}{3+2s}$
    	\begin{align}
    	\|L_n \|_{C^1((\S^2)^n;( H^{-s}_w(\R^3))^{3n})} \leq C \sqrt{n}, \\
		\|L_n -  L_{n,\app}\|_{C^1((\S^2)^n;( H^{-s}_w(\R^3))^{3n})} \leq C \sqrt{n} (\phi_n \log n+r^{-2+\frac{3}{p}}) \label{est:L_n.L_n,app.cor}.
	\end{align}
\end{cor}
    
\begin{proof}[Proof of Theorem \ref{th:well-posed}]
The assertion follows immediately from Corollary \ref{cor:L_n.H^-s} combined with Proposition \ref{pro:Strato.infinite} (see also Remark \ref{rem:manifold}).
\end{proof}

 \section{Passage to the limit for Deborah numbers of order  \texorpdfstring{$1$}{1}
}\label{section5}


In this section we prove Theorem \ref{th:main.T_D}. We recall the strategy from Subsection \ref{subsec:proof_strategy}. First, we show  that the laws of the empirical measure of the particles and of the fluid velocity field respectively are tight in suitable function spaces. For $S_n$ this is classical but we include the proof for completeness in Subsection \ref{sec:tightness.S_n}.
Tightness of $\phi_n^{-1} u_n$, which we show in Subsection \ref{sec:tightness.u_n}, follows from the estimates in Section \ref{sec:L_n} which also allow us to replace $u_n$ by  more explicit functions $u_{n,app}$.

The tightness of the laws implies weak convergence along subsequences by the Prokhorov Theorem.
To conclude the proof of Theorem \ref{th:main.T_D} we introduce certain functionals in Subsection \ref{sec:sub.proof.main.T_D}. As we show later in Appendix \ref{appendixA}, these functionals vanish precisely on the solution of the desired limit system \eqref{eq:Fokker-Planck.instationary}, \eqref{eq:viscoelastic}. For the proof of Theorem \ref{th:main.T_D} it therefore  remains to show that the laws of the microscopic system
concentrate on the zeroes of these functionals as $n \to \infty$.

\subsection{Tightness of  \texorpdfstring{$S_n$}{Sn}}\label{sec:tightness.S_n}

\begin{lem}\label{lem:tightness_Sn}
    Let $T> 0$ and let $S_n$ be the empirical measure defined in \eqref{eq:empirical_measure}. 
Then, the family of laws $\{Q^{S_n}\}$ of the empirical measures $S_n$ is tight in the space $C([0,T],\mathcal{P}_1(\R^3\times \S^2))$.
\end{lem}


\begin{proof}
By standard arguments, tightness follows from the uniform bounds
\begin{equation}\label{eq:lem_bound_Sn_1}
  \E\left[\sup_{t\in[0,T]}\int_{\R^3}\int_{\S^2}(|x|+|\xi|)S_n(t)(dx,d\xi)\right]\leq C,  
\end{equation}
\begin{equation}\label{eq:lem_bound_Sn_2}
    \E\left[\int_0^T\int_0^T \frac{\mathcal{W}_1(S_n(t_1),S_n(t_2))^p}{|t_1-t_2|^{1+s p}} \dd t_1 \dd t_2\right]\leq C,
\end{equation}
for some $s\in (0,\frac12)$, $p > s^{-1}$. Indeed, by the Ascoli-Arzelà  theorem in metric spaces, 
the set
\begin{multline*}
\mathcal{K}_{M,R}=\left\{f\in C([0,T],\mathcal{P}_1(\R^3\times \S^2))\,\,s.t.\sup_{t\in[0,T]} \int_{\R^3}\int_{\S^2}(|x|+|\xi|)f_t \dd \xi \dd x \leq M,\right.\\\left.
\int_0^T\int_0^T \frac{\mathcal{W}_1(f(t_1),f(t_2))^p}{|t_1-t_2|^{1+s p}}\dd t_1 \dd t_2\leq R\right\}.
\end{multline*}
is relatively compact  in $C([0,T],\mathcal{P}_1(\R^3\times \S^2))$. Indeed, for $f \in \mathcal K_{M,R}$ and $\varphi \in C^1(\R^3 \times \S^2)$, we have for $\theta  = s - 1/p$ due to Sobolev inequality for fractional spaces (see e.g. \cite[Section 8]{di2012hitchhiker})
 \[|\langle f(t),\varphi\rangle-\langle f(s),\varphi\rangle|\leq  C |t-s|^\theta[\langle f,\varphi\rangle]_{s,p}^p \leq C R |t-s|^\theta \| \nabla \varphi \|^p_\infty. \] 
 Taking the supremum in $\varphi$ yields equicontinuity and thus $\mathcal K_{M,R}$ is compact.
 
 

For any $\eps > 0$,  by Chebyshev's inequality and \eqref{eq:lem_bound_Sn_1}--\eqref{eq:lem_bound_Sn_2}, choosing $M$ and $R$ big enough,
\begin{align*}
Q^{S_n}(\mathcal{K}^c_{M,R})< \epsilon,
\end{align*}
which yields tightness of the laws $\{Q^{S_n}\}$. Thus, it remains to show \eqref{eq:lem_bound_Sn_1}--\eqref{eq:lem_bound_Sn_2}.

We recall that the particle positions $x_i$ do not evolve in time and stay in a compact set $K \subset \R^3$ due to assumption \eqref{ass:uniform_bound}.
Moreover the particle orientations lie on the sphere $\S^2$ which is also compact, thus the first estimate \eqref{eq:lem_bound_Sn_1} is trivial.

By definition of the Wasserstein distance, by Jensen inequality and by \eqref{eq:strato.xi_i.T_D}
 \begin{align}\label{eq:bound_Sn_1}
     \E\Ll[\mathcal{W}_1(S_n(t_1),S_n(t_2))^p\Rr]&\leq \E\Ll[\left(\frac1n \sum_{i=1}^n \Ll|\xi_i(t_1)-\xi_i(t_2)\Rr|\right)^p\Rr] \notag\\
    &\leq C \frac1n \sum_{i=1}^n  \E\Ll[\bigg|\int_{t_1}^{t_2} \sigma_D(\xi_i(\tau)) \dd B_i(\tau) \bigg|^p
    + \bigg| \int_{t_1}^{t_2} A(\tau,x_i,\xi_i(\tau)) \dd \tau  \bigg|^p\Rr], 
 \end{align}
where $A(t,x,\xi):=P_{\xi_i^\perp} h(t,\xi,x) -\am{2\xi}  $. Then, by the Burkholder-Davis-Gundy inequality, (\ml{see e.g. Theorem 3.28 in \cite{karatzas2012brownian}}) 
we deduce
\begin{align}\label{eq:bound_Sn_2}
    \E\Ll[\mathcal{W}_1(S_n(t_1),S_n(t_2))^p\Rr]&
    \leq C \frac1n \sum_{i=1}^n\E\Ll[\Ll(\int_{t_1}^{t_2} \norm{\sigma_D(\xi_i(\tau))}_{HS}^2\dd \tau \Rr)^{p/2}\Rr]+C |t_1-t_2|^p \notag \\ 
    &\leq C|t_1-t_2|^{p/2},
\end{align}
where HS stands for the Hilbert-Schmidt norm.
Therefore,
\[\E\left[\int_0^T\int_0^T \frac{\mathcal{W}_1(S_n(t_1),S_n(t_2))^p}{|t_1-t_2|^{1+s p}}\dd t_1 \dd t_2\right]\leq C \int_0^T\int_0^T |t_1-t_2|^{p(1/2 - s) -1}\dd t_1 \dd t_2 \leq C_2\]
since $s < \frac 1 2$.
\end{proof}

\subsection{Tightness of  \texorpdfstring{$u_n$}{un}}
\label{sec:tightness.u_n}

Recall from \eqref{eq:u_n.U_n}--\eqref{def:U_n.T_D} the definition of $u_n$.
Analogously, we define
\begin{align} \label{def:u_n,app}
	u_{n,app} &:= U_{n,app}', \\
	U_{n,app}(t) &:= \frac {{\phi_n}\am{\sqrt{2\gamma_{rot}}}} n \int_0^t L_{n,app}(\Xi(s)) \sqrt{\mathfrak R_2(\Xi(s))} \circ \dd (B_1(s), \dots, B_n(s)),
\end{align}
where $L_{n,app}$ is the operator defined at the beginning of Section \ref{sec:L_n}.
 \begin{lem}\label{lem:bound_un.T_D} 
 Let $s > \frac 1 2 $. Then, for all $n \in \N$, the stochastic integral in \eqref{def:u_n,app}  is well defined and $u_{n,app} \in L^2(\Omega; H^{-s}(0,T; H_{\mathfrak s}^{-s}(\R^3)))$. 
Moreover, there exists $N_0 \in \N$ such that for all for all  $n \geq N_0$ and all $s>\frac12$
\begin{align}\label{eq:lem_bound_un.T_D_1}
  \E\Ll[\norm{\phi_n^{-1}u_n}^2_{H^{-s}((0,T); H^{-s}_w(\R^3))}\Rr]&\leq C, \\
\label{eq:lem_bound_un.T_D_2}
    \lim_{n\to\infty}\E\Ll[\norm{\phi_n^{-1}(u_n-u_{n,app})}^2_{H^{-s}((0,T); H^{-s}_w(\R^3))}\Rr]&=0.
\end{align}
\end{lem}
\begin{proof}
To estimate the ${H^{-s}((0,T); H^{-s}_w(\R^3))}$ norm of $u_n$ and $u_n - u_{n,app}$, by \eqref{eq:char.negative.sobolev}, it suffices to estimate the corresponding $H^{-s+1}((0,T),H^{-s}_w(\R^3))$ norms of $U_n$ and $U_n - U_{n,app}$.
Thus, by assumption \eqref{ass:phi.log.n}, to show \eqref{eq:lem_bound_un.T_D_1}--\eqref{eq:lem_bound_un.T_D_2}, it suffices to prove
\begin{align}\label{eq:lem_bound_Un.T_D_1}
  \E\Ll[\norm{\phi_n^{-1}U_{n,app}}^2_{H^{1-s}((0,T);H^{-s}_w(\R^3))}\Rr]\leq C, \\
    \E\Ll[\norm{\phi_n^{-1}(U_n-U_{n,app})}^2_{H^{1-s}((0,T);H^{-s}_w(\R^3))}\Rr]\leq C (\phi_n \log n \rh{+ r^{3/p -2}})^2, \label{eq:lem_bound_Un.T_D_2}
\end{align}
\rh{for $p= \frac{6}{3+2s}$}. 

To estimate these norms, we appeal to \eqref{est:Z.H^s} and Remark \ref{rem:manifold}. To apply this estimate, we  
 first recall  from \eqref{eq:Particles.T_D} that for the vector $\Xi=(\xi_1,\dots, \xi_n)$, 
\[d\Xi(t)=\Sigma_D(\Xi(t))\circ \dd B(t)+H(t,\Xi(t),x) \dd t\]
where $H(t,\Xi,x)=\Ll(P_{\xi_1^\perp} h(t,\xi_1,x_1),\dots, P_{\xi_n^\perp} h(t,\xi_n,x_n)\Rr)$, $B=\Ll(B_1,\dots,B_n\Rr)$ and \am{$\Sigma_D(\Xi)$ is a block diagonal matrix in $\R^{3n\times 3n}$ whose blocks are $\sigma_D(\xi_i)$}, defined in \eqref{eq:sigma}. In particular
$\|\Sigma_D(\Xi)\|_{HS} \leq C \sqrt n$. 
Then, by \eqref{est:Z.H^s} and Remark \ref{rem:manifold}
\begin{align}
    &\left(\E\Ll[\norm{\phi_n^{-1}U_{n,app}}^2_{H^{1-s}((0,T);H^{-s}_w(\R^3))}\Rr]\right)^{\frac 1 2} \\
    &\leq \frac C n \left(\|L_{n,app}\|_{L^\infty((\S^2)^n;(H^{-s}_w(\R^3))^{3n})} + \sqrt n \|\nabla L_{n,app} \|_{L^\infty((\S^2)^n;(H^{-s}_w(\R^3))^{3n\times3n})}\right),
\end{align}
Applying Corollary \ref{cor:L_n.H^-s} yields \eqref{eq:lem_bound_Un.T_D_1}.
The proof of estimate \eqref{eq:lem_bound_Un.T_D_2} is completely analogous.
\end{proof}

\begin{lem}\label{lem:tightness_un}
For all $s > 1/2$, the family of laws $\{Q^{u_n}\}$ of $u^n$, defined in \eqref{def:U_n.T_D}, is tight in the space $H^{-s}((0,T),H^{-s}_{\mathfrak s}(\R^3))$.
\end{lem}
\begin{proof}
For $s>z>\frac12$, let
\[\mathcal{K}_M=\{v \in H^{-z}((0,T),H^{-z}_{\mathfrak s,w}(\R^3)) : \norm{v}_{H^{-z}((0,T), H^{-z}_w(\R^3))}\leq M\}.\]
By Lemma \ref{lem:compact.weighted} and \cite{Amann00}[Theorem 5.1], this set is relatively compact  in $H^{-s}((0,T),H^{-s}(\R^3))$. Thus, for all $\eps > 0$, by Chebyshev's inequality, by Lemma \ref{lem:bound_un.T_D} and choosing $M$ big enough,
\[
Q^{u_n}(\mathcal{K}^c_M)=\P\Ll(\norm{u_n}_{H^{-s}((0,T), H^{-s}_w(\R^3))}>M\Rr)\leq \frac{\E\Ll[\norm{u_n}_{H^{-s}((0,T), H^{-s}_w( \R^3))}\Rr]}{M}<\eps,
\]
which concludes the proof.
\end{proof}

\subsection{Proof of Theorem \texorpdfstring{\ref{th:main.T_D}}{main TD}}
\label{sec:sub.proof.main.T_D}
As outlined at the beginning of this Section, we now introduce functionals whose zeroes are the solutions of the system \eqref{eq:Fokker-Planck.instationary}, \eqref{eq:viscoelastic}. For $(v,g) \in H^{-s}((0,T); H_{\mathfrak s}^{-s}(\R^3))\times C([0,T];\mathcal{P}_1(\R^3\times \S^2))$, we define for each $\varphi\in C^\infty_c((0,T)\times  \R^3)$ with $\dv \varphi = 0$ and $\psi\in C^\infty_c([0,T]\times \R^3\times\S^2)$,
\begin{align}\label{eq:empiricalmeasure_weak}
\Psi_\psi(g)&= \langle f_0,\psi(0)\rangle - \langle g(T),\psi(T)\rangle + \int_0^T\langle g,\partial_t \psi  + \nabla_\xi\psi \cdot P_{\xi^\perp} h +  \Delta_\xi \psi\rangle \dd t,\\
\label{eq:formal.approximation_weak}
\Phi_\varphi(v,g)&=
\langle v,\Delta \varphi\rangle - {\gamma_E}\langle g , (3\xi \otimes \xi - \textrm{Id}):\nabla \varphi \rangle.  
\end{align}
Note that $\langle \cdot, \cdot \rangle$ in \eqref{eq:empiricalmeasure_weak} denotes the pairing in $\R^3 \times \S^2$ whereas \rhnew{the first pairing $\langle \cdot, \cdot \rangle$ in \eqref{eq:formal.approximation_weak} denotes the pairing in $(0,T) \times \R^3$ and the second pairing $\langle \cdot, \cdot \rangle$ in \eqref{eq:formal.approximation_weak} denotes the pairing in $(0,T) \times \R^3 \times \S^2$.} 

 \begin{lem} \label{lem:identity_S_n}
For every test function $\psi\in C^\infty_c([0,T]\times\R^3\times \S^2)$, the empirical measures $S_n$ satisfies the following identity,
\begin{align*}
    \Psi_\psi(S_n) = \rhnew{\langle f_0 - S_n(0),\psi(0)\rangle -} \int_0^T \frac{1}{n}\sum_{i=1}^n \nabla \psi(t,x_i,\xi_i(t))\sigma_D(\xi_i(t)) dB_i(t)
\end{align*}
with $\sigma_D$ as in \eqref{eq:sigma}.
\end{lem}
\begin{proof}
The proof follows by applying Ito's formula   (see \cite[Theorem 7.4.3]{Kuo} ) to $\psi(t,x_i,\xi_i(t))$ using \eqref{eq:strato.xi_i.T_D}. 
To avoid the issue of having to apply Ito's formula on the manifold $\S^2$, one might first extend $\psi$ to a smooth function $\tilde \psi \in C_c^\infty([0,T] \times \R^3 \times \R^3\setminus \{0 \})$.
Then, observing that $\frac{1}{2}\sigma_D(\xi)\sigma_D(\xi)^T = \Id - \xi \otimes \xi = \rh{\Id - \frac{\xi \otimes \xi}{|\xi|^2} =}  P_{\xi^\perp}$ \rh{on $\S^2$}, we find
\begin{align*}
    \tilde \psi(T,x_i,\xi_i(T)) &= \tilde \psi(0,x_i,\xi_{i,0}) +  \int_0^T (\partial_t  + P_{\xi^\perp} h \cdot \nabla_\xi - 2 \xi \cdot \nabla_\xi +   P_{\xi^\perp} : \nabla_\xi^2) \tilde \psi(t,x_i,\xi_i(t)) \dd t \\
    &+ \int_0^T \frac{1}{n}\sum_{i=1}^n \nabla_\xi \tilde \psi(t,x_i,\xi_i(t))\sigma_D(\xi_i(t)) dB_i(t),
\end{align*}
We then observe that
\begin{align}
P_{\xi^\perp} : \nabla_\xi^2 \tilde \psi -   \rh{ 2} \xi \cdot \nabla_\xi \tilde \psi = \dv ( P_{\xi^\perp} \nabla \tilde \psi) = \Delta_{\S^2} \tilde \psi,
\end{align}
where $\Delta_{\S^2}$
is the Laplace-Beltrami operator. This concludes the proof.
\end{proof}

\begin{prop}\label{prop:Q_sol_delta}
Denote by $Q^n$ the law of $(\phi_n^{-1}u_n,S_n)$ on the space $H^{-s}((0,T); H^{-s}_{\mathfrak s}(\R^3))\times C([0,T],\mathcal{P}_1(\R^3\times \S^2))$.
Then for all $\delta>0$, all $\varphi\in C^\infty_c((0,T),{\R^3})$ with $\dv \varphi = 0$ and for all $\psi\in C^\infty_c([0,T]\times \R^3\times\S^2)$
\begin{multline*}
\lim_{n\to\infty} Q^{n}\Big((v,g)\in H^{-s}((0,T); H^{-s}_{\mathfrak s}(\R^3))\times C([0,T],\mathcal{P}_1(\R^3\times \S^2)) :\\
\,\, |\Phi_\varphi(v,g)|+|\Psi_\psi(g)|>\delta\Big)=0.
\end{multline*} 
\end{prop}
\begin{proof}
By Chebyshev's inequality, it suffices to show that
\begin{align}\label{eq:Q_step0}
\lim_{n \to \infty} \E\Ll[ |\Phi_\varphi(\phi_n^{-1}u_{n},S_{n})| +  |\Psi_\psi(S_{n})|\Rr] = 0.
\end{align}
We start by studying the first term.
With $u_{n,app}$ as in \eqref{def:u_n,app} we have
\begin{align}
    \langle \phi_n^{-1}u_{n,app},\Delta\varphi \rangle &= \rhnew{\phi_n^{-1}}\int_0^T \langle U_{n,app}(t), \partial_t \Delta  \varphi(t,\cdot) \rangle \dd t \\
    &= \frac{\sqrt{2 \gamma_{rot}}}{n} \int_0^T \left\langle \int_0^t  L_{n,app}(\Xi(s)) \sqrt{\mathfrak R_2(\Xi(s))} \circ \dd B(s), \partial_t \Delta  \varphi(t,\cdot)  \right \rangle \dd t .
    \end{align}
    \rh{Now let us recall that since $L:=L_{n,app}(\Xi(s)) \sqrt{\mathfrak R_2(\Xi(s))} \in \mathcal{L}(\R^{3n}, H)$ with $H:=H^{-s}_{\mathfrak s}(\R^3)$, we may identify $\mathcal{L}(\R^{3n}, H)$ with $H^{3n}$ through $L_{i,\alpha} = L e_{i,\alpha}$ where  $e_{i,\alpha}$ denote the canonical basis vectors of $\R^{3n}$.
    Then, with $\langle \cdot , \cdot \rangle$ denoting the pairing between $H^{-s}(\R^3) $ and $H^{s}(\R^3) $, and $(\cdot,\cdot)$ the scalar product in $H$, we use $L_{i,\alpha} = \sum_k ( L_{i,\alpha}, \epsilon_k ) \epsilon_k$, where $(\epsilon_k)_k$ is an orthonormal basis of $H$, in order to get
\begin{align}
   \sum_{i,\alpha} \left \langle \int_0^t L_{i,\alpha}(s) \circ \dd B_{i,\alpha}(s), h \right \rangle 
   &=   \sum_{i,\alpha,k} \left \langle \int_0^t ( L_{i,\alpha}(s), \epsilon_k ) \epsilon_k \circ \dd B_{i,\alpha}(s), h \right \rangle \\
   &=  \sum_{i,\alpha,k}\int_0^t \langle ( L_{i,\alpha}(s), \epsilon_k ) \epsilon_k, h  \rangle  \circ \dd B_{i,\alpha}(s)   \\
     &=  \sum_{i,\alpha}\int_0^t \langle L_{i,\alpha}(s) , h\rangle \circ \dd B_{i,\alpha}(s). 
\end{align}
}
\am{Hence, this yields, with $e_\alpha$ denoting the canonical basis vectors of $\R^3$,
    \begin{align}
    &\langle \phi_n^{-1}u_{n,app},\Delta\varphi \rangle\\
    &= \frac{\sqrt{2 \gamma_{rot}}}{n}\underset{i,\alpha,\beta}{\sum} \int_0^T \int_0^t \langle L_{n,app}(\Xi(s))e_{i,\alpha}, \partial_t \Delta  \varphi(t,\cdot) \rangle \left(\sqrt{\mathcal R_2(\xi_i(s))}\right)_{\alpha,\beta} \circ \dd B_{i,\beta}(s) \dd t \\
&=-\frac{\sqrt{2 \gamma_{rot}}}{n} \underset{i,\alpha,\beta}{\sum} \int_0^T \int_0^t \left([e_{\alpha}]_M+\mS(\xi_i(s))e_{\alpha} \right):\partial_t \nabla \varphi(t,x_i) \left(\sqrt{\mathcal R_2(\xi_i(s))}\right)_{\alpha,\beta} \circ \dd B_{i,\beta}(s) \dd t \\
&= -\frac{\sqrt{2 \gamma_{rot}}}{n} \sum_i\int_0^T \left(\int_0^t   \partial_t D \varphi(t,x_i):\mS(\xi_i(s))\sqrt{\mR_2(\xi_i(s))} \circ \dd B_i(s)  \right.\\
&\qquad \qquad  \qquad  \qquad \qquad   \left.+2 \int_0^t\partial_t\curl  \varphi(t,x_i)\cdot  \sqrt{\mR_2(\xi_i(s))} \circ \dd B_i(s) \right) \amnew{\dd t}\\
&= +\frac{\sqrt{2 \gamma_{rot}}}{n} \sum_i \int_0^T  D \varphi(t,x_i):\mS(\xi_i(t))\sqrt{\mR_2(\xi_i(t))}\circ \dd B_i(t) \\
& +\frac{\sqrt{2 \gamma_{rot}}}{n} \sum_i\int_0^T 2\curl \varphi(t,x_i)\cdot  \sqrt{\mR_2(\xi_i(t))} \circ \dd B_i(t) .
\end{align}
where we used \eqref{eq:distributional.v}, \eqref{eq:skew_curl} and an integration by parts formula for Stratonovitch integrals which follows from the chain rule, Remark \ref{rem:Strato} \ref{it:stratonovich.chain.rule}.
}


Thus,
\begin{align} \label{eq:I_1.I_2}
\Phi_\varphi(\phi_n^{-1}u_{n_k},S_{n_k})=\mathcal{I}_1+\mathcal{I}_2,
\end{align}
where 

\[\mathcal{I}_1=\langle \phi_n^{-1}\Ll(u_{n_k}-u_{n_k,app}\Rr),\Delta \varphi\rangle,\]
and
\begin{align*}
\mathcal{I}_2=:  \frac 1 n \sum_i \mathcal J_i&= \frac 1 n \sum_i\left( \frac{\sqrt{2 \gamma_{rot}}}{n} \int_0^T 2\curl \varphi(t,x_i)\cdot  \sqrt{\mR_2(\xi_i(t))} \circ \dd B_i(t)   \right.\\
&+ \frac{\sqrt{2 \gamma_{rot}}}{n} \int_0^T   D \varphi(t,x_i):\mS(\xi_i(t))\sqrt{\mR_2(\xi_i(t))} \circ \dd B_i(t) \\
& \left. - \gamma_E\int_0^T(3\xi_i \otimes \xi_i-\textrm{Id}):\nabla \varphi(t,x_i) \dd t \right) .
\end{align*}
By  Lemma \ref{lem:bound_un.T_D}, we have 
\begin{align} \label{est:I_1}
        \lim_{n \to \infty} \E\Ll[|\mathcal{I}_1|\Rr] = 0.
\end{align}
To estimate $\mathcal{I}_2$, we use Lemma \ref{lem:Expectations},
which implies $\E[\mathcal J_i] = 0$.
Furthermore, since the particle orientations $\xi_i$ are independent, also the terms $\mathcal J_i$ are independent. \rh{Moreover, they have a bounded variance by the It\^o Stratonovitch conversion formula, Remark \ref{rem:Strato} \ref{it:Conversion.composition} and It\^o isometry,  (see \cite[Theorem 2.3.4]{Kuo}).}  Therefore.
\begin{align}\label{est:I_2}
    \E[\mathcal I_2^2] = \frac 1 {n^2} \sum_i \E[\mathcal J_i^2] \leq \frac{C}{n} \|\nabla \varphi \|^2_{L^\infty}. 
\end{align}
Inserting \eqref{est:I_1}--\eqref{est:I_2} in \eqref{eq:I_1.I_2} yields
\begin{align}\label{est:Phi_phi}
\lim_{n \to \infty} \E\Ll[ |\Phi_\varphi(\phi_n^{-1}u_{n_k},S_{n_k})|\Rr] = 0.
\end{align}

Regarding the second term of \eqref{eq:Q_step0}, Lemma \ref{lem:identity_S_n} implies
\[\Psi_\psi(S_n) - \langle f_0 - S_n(0),\psi(0)\rangle  =  -\frac{1}{n}\sum_{i=1}^n\int_0^T \nabla \psi(t,x_i,\xi_i(t))\sigma_D(\xi_i(t)) \dd B_i(t) =: \frac 1 n \sum_i \tilde{\mathcal J_i}.\]
\rhnew{Using assumption \eqref{ass:initial.convergence}, $\E[|\langle f_0 - S_n(0),\psi(0)\rangle|] \to 0$}. 

Combining this with  independence of $\tilde{\mathcal J_i}$ \rh{and bounded variance} as above, we conclude
\begin{align} \label{est:Psi_psi}
    \lim_{n \to \infty} \E\Ll[\Ll|\Psi_\psi(S_n)\Rr|\Rr] = 0.
\end{align}
Combination of \eqref{est:Phi_phi} and \eqref{est:Psi_psi} yields \eqref{eq:Q_step0}
and the proof is completed.
\end{proof}

To be in position to prove Theorem \ref{th:main.T_D}, 
we state the following uniqueness results which are proved in Appendix \ref{appendixA}.

\begin{thm} \label{thm:uniqueness.instationary}
   Let $f_0 \in  L^2(\S^2\times \R^3)$ and $f \in C([0,T];\mathcal P_1(\R^3 \times \S^2))$ satisfying $\Psi_\psi(f) = 0$ for all
    $\psi  \in C_c^\infty([0,T] \times \R^3 \times \S^2)$.
     \am{Then, $f$ is the unique weak solution to \eqref{eq:Fokker-Planck.instationary} \rhnew{in $C([0,T];L^2(\R^3 \times \S^2))$} such that for almost all $x\in \R^3$, $f(\cdot,\cdot,x)\in L^2(0,T; H^1(\S^2)) $ and $f'(\cdot,\cdot,x)\in L^2(0,T; H^{-1}(\S^2)) $}.
\end{thm}

\begin{thm} \label{thm:uniqueness.Stokes}
\am{Let $f \in L^2((0,T) \times \R^3 \times \S^2)$.} 
Assume $u \in H^{-s}((0,T),H_{\mathfrak s}^{-s}(\R^3))$ satisfies $\Phi_\varphi(u,f) = 0$ for all  $\varphi\in C^\infty_c((0,T),\R^3)$ with $\dv \varphi = 0$.  
Then $u$ is the unique weak solution in
$L^2(0,T;\dot H_{\mathfrak s}^1(\R^3))$ to \eqref{eq:viscoelastic}.
\end{thm}

\begin{proof}[Proof of Theorem \ref{th:main.T_D}]
By Lemmas \ref{lem:tightness_Sn} and \ref{lem:tightness_un}, the law  $\{Q^n\}$ of $(S_n,u_n)$ is tight, and thus we can extract a convergent subsequence $Q^{n_k}$ which converges weakly to some probability measure $Q$ on $H^{-s}((0,T); H^{-s}_{\mathfrak s}( \R^3)\times C([0,T];\mathcal{P}_1(\R^3\times \S^2))$.
We will argue that $Q= \delta_{(u,f)}$, where $(u,f)$ is the unique solution to \eqref{eq:Fokker-Planck.instationary}--\eqref{eq:viscoelastic}. Then, by standard arguments, uniqueness of $(u,f)$ implies weak convergence of the whole sequence $Q^n$, and, since $Q$ is deterministic, we deduce convergence in probability of $(u_n,S_n)$ to $(u,f)$.

It thus remains to show that $\delta_{(u,f)}$ is the only accumulation point of $Q^n$.
Let $Q^{n_k}$ be a converging subsequence and $Q$ its limit.
First we observe that for all  $\varphi\in C^\infty_c((0,T) \times \R^3)$ with $\dv \varphi = 0$  and 
all $\psi\in C^\infty_c([0,T]\times \R^3\times\S^2)$ the functionals $\Psi_\psi$ and $\Phi_\varphi$ are continuous with respect to the topology of $ H^{-s}(0,T; H^{-s}_{\mathfrak s}(\R^3))\times C([0,T],\mathcal{P}_1(\R^3\times \S^2)) $. Therefore, by Portmanteau theorem and Proposition  \ref{prop:Q_sol_delta}, for all $\delta > 0$
\[Q\Ll((v,g)\,:\, |\Phi_\varphi(v,g)|+|\Psi_\psi(g)|>\delta\Rr)\leq\liminf_{k\to\infty} Q^{n_k}\Ll((v,g)\,:\, |\Phi_\varphi(v,g)|+|\Psi_\psi(g)|>\delta\Rr)=0,\]
which implies
\[Q\Ll((v,g)\,\,:\,\, |\Phi_\varphi(v,g)|+|\Psi_\psi(g)|=0\Rr)=1.\]
Since the functionals $\Psi_\psi$ and $\Phi_\varphi$ are also continuous with respect to $\psi$ and $\varphi$, a density argument yields
\begin{align*}
&Q\Ll((v,g)\,:\, |\Phi_\varphi(v,g)|+|\Psi_\psi(g)|=0,\, \forall \,\varphi\in C^\infty_c((0,T),\R^3) \text{ with } \dv \varphi = 0, \Rr. \\
& \qquad \qquad \Ll. \psi\in C^\infty_c([0,T]\times \R^3\times\S^2)\Rr)=1.
\end{align*}

Thus, Theorem \ref{thm:uniqueness.instationary} implies that the law of $S^n$ concentrates indeed on $f$.
Finally, by Theorem \ref{thm:uniqueness.Stokes} we conclude  $Q=\delta_{(u,f)}$.
\end{proof}

 \section{Passage to the limit for very small Deborah numbers}\label{section6}
 
 
 In this section, we prove Theorem \ref{th:main.T_u}, i.e. the passage to the quasi-stationary system \eqref{eq:Fokker-Planck.stationary}, \eqref{eq:viscoelastic} starting from  the microscopic dynamics \eqref{eq:Stokes.micro.T_u}--\eqref{eq:Particles.T_u}.
 Large parts of the proof are analogous to the the proof of Theorem \ref{th:main.T_D} given in the previous section.
 
 The main difference concerns the tightness of the law of the empirical measure $S_n$ (defined in \eqref{eq:empirical_measure}). Indeed, we cannot expect tightness in $C([0,T];\mathcal{P}_1(\R^3 \times \S^2))$ since the solution to the limit problem \eqref{eq:Fokker-Planck.stationary} is discontinous at time $0$. This discontinuity arises from the fast diffusion that induces a boundary layer at the initial time. This fast diffusion also makes tightness in $C([\delta,T];\mathcal{P}_1(\R^3 \times \S^2))$ difficult to obtain.
 We therefore work in a weaker space instead, namely
 $H^{0_-}((0,T);H^{-3/2_-}(\R^3 \times \S^2))$.
\begin{lem}\label{lem:tightness_Sn.T_u}
    Let $T> 0$ and let $S_n$ be the empirical measure defined in \eqref{eq:empirical_measure}. 
Then, the family of laws $\{Q^{S_n}\}$ of the empirical measures $S_n$ is tight in the space  $H^{0_-}((0,T);H^{-3/2_-}(\R^3 \times \S^2))$.
\end{lem}
\begin{proof}
 We note that since 
 $H^{3/2_+}(\R^3 \times \S^2)$ is compactly embedded into $C^{0}(\R^3 \times \S^2)$, we have
 that $S_n$ is uniformly bounded 
 in $L^\infty((0,T);H^{-3/2_-}(\R^3 \times \S^2))$
 and thus compactly embedded in $H^{0_-}((0,T);H^{-3/2_-}(\R^3 \times \S^2))$.
\end{proof}

Since the fast diffusion is balanced by the different scaling of $u_n$, tightness of the law of $u_n$ is completely analogous as in the previous section.
 \begin{lem}\label{lem:bound_un.T_u}
 Let
 \begin{align}\label{def:U_n,app_T_u}
	u_{n,app} &:= U_{n,app}', \\
	   U_{n,app}(t) &:= \frac {\sqrt{\am{2 \gamma_{rot}\phi_n}}} n \int_0^t L_{n,app}(\Xi(s)) \sqrt{\mathfrak R_2(\Xi(s))} \circ \dd (B_1(s), \dots, B_n(s)).
\end{align}
Then, $u_{n,app} \in L^2(\Omega;H^{-s}(0,T;H^{-s}_{\mathfrak s}(\R^3)))$ is well-defined and there exists $N_0 \in \N$ such that for all $s>\frac12$ and all  $n \geq N_0$
\begin{align}\label{eq:lem_bound_un.T_u_1}
   \E\Ll[\norm{u_n}^2_{H^{-s}((0,T); H^{-s}_w(\R^3))}\Rr]&\leq C, \\
\label{eq:lem_bound_un.T_u_2}
    \lim_{n\to\infty}\E\Ll[\norm{u_n-u_{n,app}}^2_{H^{-s}((0,T); H^{-s}_w(\R^3))}\Rr]&=0.
\end{align}

In particular, the family of laws $\{Q^{u_n}\}$ of $u^n$, defined through \eqref{def:U_n}--\eqref{eq:u_n.U_n}, is tight in the space $H^{-s}((0,T),H^{-s}_{\mathfrak s}(\R^3))$.
\end{lem}

\medskip

In order to pass to the limit, we consider again functionals $\bar \Psi_\psi(g), \Phi_\varphi(v,g), \Theta_\theta(g)$ for $(v,g) \in H^{-s}((0,T); H^{-s}_{\mathfrak s}(\R^3))\times  H^{0_-}((0,T);H^{-3/2_-})$  and $\varphi\in C^\infty_c((0,T)\times \R^3)$ with $\dv \varphi = 0$, $\psi\in C^\infty_c((0,T)\times \R^3\times\S^2)$
and $\theta \in C_c^\infty((0,T) \times \R^3)$.
Here,  $\Phi_\varphi$ is still given by \eqref{eq:formal.approximation_weak},  $\bar \Psi_\psi$ and $\Theta_\theta$ are defined as 
\begin{align}\label{eq:empiricalmeasure_weak_stationnary}
\bar \Psi_\psi(g)&:=\langle g,\Delta_\xi \psi \rhnew{+}\nabla_\xi\psi \cdot P_{\xi^\perp} h\rangle, \\
\label{eq:def_Theta}
     \Theta_\theta(g) &:= \langle f_0 -  g, \theta \otimes 1 \rangle.
 \end{align}
Correspondingly to Lemma \ref{lem:identity_S_n}, we have the following.
\begin{lem} \label{lem:identity_S_n.T_u}
For every test function $\psi\in C^\infty_c((0,T)\times\R^3\times \S^2)$, the empirical measures $S_n$ satisfies the following identity,
\begin{align*}
    \bar \Psi_\psi(S_n) &= \rhnew{-} \phi_n\int_0^T\langle S_n(t),\partial_t \psi(t,x,\xi)\rangle \rhnew{-} \int_0^T \frac{1}{n}\sum_{i=1}^n \nabla \psi(t,x,\xi_i(t))\sigma_D(\xi_i(t)) dB_i(t),
\end{align*}
with $\sigma_D$ as in in \eqref{eq:sigma}.
\end{lem}

 The following proposition is the analogous version of Proposition \ref{prop:Q_sol_delta}.
 
\begin{prop}\label{prop:Q_sol_delta.T_u}
For each $\delta>0$, for each $\varphi\in C^\infty_c((0,T),\R^3)$ with $\dv \varphi = 0$, $\psi\in C^\infty_c((0,T)\times \R^3\times\S^2)$ and $\theta \in C_c^\infty((0,T) \times \R^3)$ 
\begin{multline*}
\lim_{n\to\infty} Q^{n}\Big((v,g)\in H^{-s}((0,T); H^{-s}_{\mathfrak s}(\R^3))\times  H^{0_-}((0,T);H^{-3/2_-}(\R^3 \times \S^2)) :\\
|\Phi_\varphi(v,g)|+|\bar \Psi_\psi(g)| + |\Theta_\theta(g)|>\delta\Big)=0.
\end{multline*} 
\end{prop} 
 \begin{proof}
The proof is analogous to the one of Proposition \ref{prop:Q_sol_delta}.\am{ Note in particular that Lemma \ref{lem:Expectations} has to be adapted accordingly since $\Xi$ satisfies \eqref{eq:Particles.T_u}  instead of \eqref{eq:Particles.T_D} which yields the appearance of a factor $\sqrt{\phi_n}^{-1}$ in \eqref{eq:stress.average} that will be compensated by the factor $\sqrt{\phi_n}$ appearing in the definition of $U_{n,app}$, see \eqref{def:U_n,app_T_u}.}
To deal with the additional functional $\Theta_\theta$, we notice that
\begin{align}
    \Theta_\theta(S_n) = \int_0^T \langle \frac 1 n \sum_i \delta_{x_i, \xi_{i,0}} - f_0, \theta(t,\cdot) \otimes 1 \rangle \dd t,
\end{align}
and therefore this functional can be handled due to assumption \eqref{ass:initial.convergence}.
 \end{proof}

\begin{proof}[Proof of Theorem \ref{th:main.T_u}]
The proof is completely analogous to the proof of Theorem \ref{th:main.T_D}.  The only difference is that regarding the uniqueness, we rely on Theorem \ref{thm:uniqueness.stationary} below instead of Theorem \ref{thm:uniqueness.instationary}.
\end{proof}

\begin{thm} \label{thm:uniqueness.stationary}
    Let $f_0\in L^2(\S^2\times \R^3)$ and $g \in H^{0_-}((0,T);H^{-3/2_-})$ satisfy $\bar \Psi_\psi(g) = 0$, $\Theta_\theta(g)= 0$ for all
    $\psi  \in C_c^\infty((0,T) \times \R^3 \times \S^2)$ and
     $\theta \in C_c^\infty((0,T) \times \R^3)$.
     \am{Then, $g$ is the unique weak solution to \eqref{eq:Fokker-Planck.stationary} in $L^2((0,T)\times \S^2\times \R^3)$ such that for almost all $x\in \R^3$, $f(\cdot,\cdot,x)\in C^\infty((0,T)\times \S^2)$.}
\end{thm}
 



\appendix

\section{Uniqueness results for the Stokes and Fokker-Planck equations in negative Sobolev spaces}\label{appendixA}

In this appendix we show the uniqueness results for solutions of the Stokes and the (in-)stationary Fokker-Planck equations as stated in Theorems \ref{thm:uniqueness.Stokes}, \ref{thm:uniqueness.instationary} and \ref{thm:uniqueness.stationary}.

\subsection{Proof of Theorem \ref{thm:uniqueness.Stokes}}

\begin{proof}[Proof of Theorem \ref{thm:uniqueness.Stokes}]
Well-posedness of  \eqref{eq:viscoelastic} in \am{$L^2(0,T;\dot H^1_{\mathfrak{s}}(\R^3))$} is classical.
By linearity, it therefore remains to show that there is at most one function $u \in H^{-s}((0,T),H^{-s}(\R^3))$ that satisfies $\Phi_\varphi(u,0) = 0$ for all  $\varphi\in C^\infty_c((0,T),\R^3)$ with $\dv \varphi = 0$.

Let $g \in C_c^{\infty}((0,T) \times \R^3)$ and define the pair
$v_g \in C_c^\infty((0,T); \dot H_{\mathfrak s}^{2 +s}(\R^3) \cap \dot H_{\mathfrak s}^2(\R^3))$, $p_g \in C_c^\infty((0,T); \dot H^{1 +s}(\R^3) \cap \dot H^1(\R^3))$ to be the solution to the Stokes equations
 \[-\Delta v_g+\nabla p_g=g,\quad \dv v_g=0.\]
 Then, by density of divergence free function of $C^\infty_{c}(\R^3)$ in $\dot H_{\mathfrak s}^{2 +s}(\R^3) \cap \dot H_{\mathfrak s}^2(\R^3)$, we have
 \begin{align}
    0 = \Phi_{v_g}(u,0) = \langle u, g \rangle,
 \end{align}
which finishes the proof.
\end{proof}

\subsection{Proof of Theorem \ref{thm:uniqueness.instationary}}


\begin{proof}[Proof of Theorem \ref{thm:uniqueness.instationary}]
    Let $\varphi \in C_c^\infty([0,T] \times \R^3 \times \S^2)$ and let $\psi  \in C_c^\infty([0,T] \times \R^3 \times \S^2)$ be the unique classical solution to the backwards parabolic equation
    \begin{equation}\left\{
    \begin{array}{rcl}
        - \partial_t \psi -  \Delta_\xi \psi - P_{\xi^\perp} h \cdot \nabla \psi  &=& \varphi , \\
        \psi(T,\cdot) &=& 0.
        \end{array} \right.
    \end{equation}
    By standard regularity theory for parabolic equations, $\psi$ is well-defined.
    
    Thus, for $f$ as in the statement
    \begin{align}
        0 = \Psi_\psi(f) = \langle f_0, \psi(0) \rangle - \int_0^T \langle f, \varphi \rangle \dd t.
    \end{align}
    Since $f_0$ is given, this identity characterizes $f$ uniquely, and therefore $f$ must coincide with the unique classical solution to \eqref{eq:Fokker-Planck.instationary}.
\end{proof}

\subsection{Proof of Theorem \ref{thm:uniqueness.stationary}}



The proof of Theorem \ref{thm:uniqueness.stationary} relies on the following theorem regarding the elliptic operator
\begin{align}
    L v = -\dv (\nabla v - \bar h v),  
\end{align}
where \am{$\bar h \in C^\infty(\S^2;\R^3)$} with
\begin{align} \label{eq:h.in.TS^2}
    \xi \cdot \am{\bar h = 0 }\quad \text{for all }\xi\in \S^2.
\end{align}
The condition \eqref{eq:h.in.TS^2} ensures that $h(t,\cdot,x)$ takes values in the tangent space $T \S^2$.

We denote the formal adjoint by 
\begin{align}
    L^\ast v := - \Delta v -\bar h \cdot \nabla v.
\end{align}

\rh{In what follows we will often deal with functions $h$ depending on a parameter $z \in \R^m$ and denote by $L_{z}$ and $L_{z}^\star$ the corresponding operators as above.}
\begin{thm}\label{th:regularity.f.bar}
\begin{enumerate}[label=(\roman*)., ref=(\roman*)]
    \item \label{it:unique.stationary} 
       Let $\bar h \in C^\infty(\S^2;\R^3)$ satisfy \eqref{eq:h.in.TS^2}. Then $\dim \ker L = 1$ and all elements in $\ker L $ have a sign. In particular, there exists a unique classical solution $\bar f$ to
       \begin{align} \label{eq:stationary.normalized}
           - L \bar f = 0, && \int_{\S^2} \bar f = 1.
       \end{align}
       
       Furthermore, this solution $\bar f$ depends smoothly on $h$. More precisely, let for a smooth family $(z,\xi) \in \R^m \times \S^2 \mapsto \bar h_z(\xi)$, $z \in \R^m$ be smooth such that $h_z$ satisfies \eqref{eq:h.in.TS^2} for all $z \in \R^m$, then the family of solutions $z \mapsto \bar f_z$ is smooth.
       
   \item   \label{it.smooth.dependence.dual}    Let $(z,\xi) \in \R^m \times \S^2 \mapsto h_z(\xi)$ be as above. Then, for all $K \Subset \R^m$, $k \in \N$ there exists $C> 0$ such that for all $z \in K$ and all $\psi \in (\ker L_z)^\perp$, the unique solution $v \in H^1(\S^2)$ to $L_z^\ast v = \psi$ with $\int v = 0$ satisfies 
    \begin{align} \label{est:uniform.spectral.gap}
        \|v\|_{H^{k+2}} \leq C \|\psi\|_{H^k}.
    \end{align}
    
    In particular, if $\varphi \in C_c^\infty(\R^m \times \S^2)$, and $v(z,\cdot)$ is the unique solution with $\int v(z,\cdot) = 0$ to $L^\ast_z v(z,\cdot) = \varphi(z,\cdot)$, then $v \in C_c^\infty(\R^m \times \S^2)$.
    \end{enumerate}
\end{thm}

\begin{proof}
We rely on classical results concerning elliptic PDEs on compact manifolds, for which we refer to \cite{Taylor}.

Existence and uniqueness for the problem \eqref{eq:stationary.normalized} is classical, see e.g. \cite{zeeman1988stability}.
We give a short proof for completeness.
The operator $L$ is a compact perturbation of the selfadjoint operator $-\Delta$. Thus, the index of $L$ is $0$, i.e.
\begin{align}
    \dim \ker L = \dim \ker L^\ast.
\end{align}

Since $L^\ast$ contains no zero-order terms, the maximum principle implies $\ker L^\ast = \{ const\}$.
In particular $\dim \ker L = \dim \ker L^\ast = 1$.
It is easy to see that $g \in \ker L$ implies that also the positive and negative parts $g_+$, $g_-$ lie in $\ker L$.
Thus, $g = g_+$ or $g= g_-$, otherwise $g_+$ and $g_-$ would be linearly independent, contradicting $\dim \ker L = 1$.

\medskip 

    The proof of the assertion that $\bar f_z$  depends smoothly on $z$
    is similar to the proof of \ref{it.smooth.dependence.dual}
    which we prove first.

  Let $v,\psi$ be as in the statement.  From standard elliptic regularity theory it follows that
    \begin{align} \label{est:elliptic.regularity}
        \|v\|_{H^{k+2}} \leq C( \|v\|_{L^2} +  \|\psi\|_{H^k}).
    \end{align}
    for a constant $C$ independent of $z$ in $K$.
    For \eqref{est:uniform.spectral.gap}, it remains to show that 
      \begin{align*}
        \|v\|_{L^2} \leq C  \|\psi\|_{L^2}.
    \end{align*}
    Assume for the sake of contradiction that there exists no such constant. Then, there exists a sequence $z_n \in K$, $\psi_n \in (\ker L_{z_n})^\perp$, $v_n \in H^1(\S^2)$ such that $L^\ast v_n = \psi_n$, $\int v_n = 0$, $\|v_n\|_{L^2} = 1$ and $\|\psi_n\|_{L^2} \leq 1/n$.
    
    
    By \eqref{est:elliptic.regularity}, $v_n$ is bounded in $H^2(\S^2)$, thus by taking a suitable subsequence $v_n \to v_\ast$ in $H^1(\S^2)$ and $z_n \to z_\ast$ for some $v_\ast \in H^2(\S^2)$, $z_\ast \in K$.
    Note that by the smoothness assumption on $h$, we have $h_{z_n} \cdot \nabla v_n  \to \nabla  h_{z_\ast} \cdot v_\ast$ in $L^2$.
  Thus, $L^\ast_{z_\ast}  v_\ast = 0$. Since $\int v_\ast = 0$ and $\|v_\ast\|_{L^2} = 1$ this contradicts $\ker L_{z_\ast}^\ast = \{const\}$. This establishes \eqref{est:uniform.spectral.gap}.
    
    \medskip
    
    For the second part of \ref{it.smooth.dependence.dual}, let $\varphi \in C_c^\infty(\R^m \times \S^2)$, and $v(z,\cdot)$ be the unique solution with $\int v(z,\cdot) = 0$ to $L^\ast_z v(z,\cdot) = \varphi(z,\cdot)$. Note that $v(z,\cdot) = 0$ if $\varphi(z,\cdot) = 0.$ Moreover, for $z_1,z_2 \in \R^m$, we have
    \begin{align*}
        L^\ast_{z_1} (v(z_1,\cdot) - v(z_2,\cdot)) =   ( h_{z_2} -  h_{z_1}) \cdot \nabla v(z_2,\cdot) + \varphi(z_1,\cdot) - \varphi(z_2,\cdot).
    \end{align*}
    Thus, by \eqref{est:uniform.spectral.gap}, we have 
    \begin{align*}
        \|v(z_1,\cdot) - v(z_2,\cdot)\|_{H^{k+2}(\S^2)} &\leq C  \| h_{z_1} - h_{z_2}\|_{W^{k,\infty}(\S^2)} \|\nabla v(z_2,\cdot)\|_{H^k(\S^2)}  \\
        &+  \|\varphi(z_1,\cdot)- \varphi (z_2,\cdot)\|_{H^k(\S^2)} \\
        & \leq C |z_1 - z_2|
    \end{align*}
    Thus, $v \in W^{1,\infty}(\R^m;H^k(\S^2))$. Differentiating the equation with respect to $z_j$
   leads to
    \begin{align}
        L^\ast \tilde v = \tilde \varphi + \tilde h \cdot \nabla v,  
    \end{align}
    where $\tilde v = \partial_{z_k} v, \tilde h = \partial_{z_k} h$ and $\tilde \varphi = \partial_{z_k} \varphi$.
    Note that $\int \tilde v(z,\cdot) = 0$. Thus, repeating the above argument, we find that $\tilde v \in W^{1,\infty}(\R^m;H^k(\S^2))$ and by induction $v \in C_c^\infty(\R^m \times \S^2)$.

    \medskip
    
    To see that $\bar f_z$ as in the second part of assertion \ref{it:unique.stationary} depends smoothly on $z$,
    we proceed analogously:
    Similarly as before, we observe that for any $\psi \in (\ker L^\ast_z)^\perp$ the unique solution to $L_z v = \psi$ with $\int v = 0$ satisfies 
    \begin{align} \label{est:regularity.L}
        \|v\|_{H^{k+2}} \leq C \|\psi\|_{H^k},
    \end{align}
    where $C$ depends only on $k$, uniformly on compact sets $K \subset \R^m$. Note that the only difference of this estimate to \eqref{est:uniform.spectral.gap} is that $\int v = 0$ is not equivalent to $v \in \ker L_z^\perp$.
    However, inspection of the proof above reveals that we never used orthogonality but only  that 
    \begin{align} \label{eq:ker.sign}
        \{v \in \ker L_z : \int v= 0 \} = \{0\}.
    \end{align}
    This still holds true, since $\dim \ker L_z =1$  and  $\bar f_z \in \ker L_z$ with $\int \bar f_z =1$. Thus, \eqref{est:regularity.L} holds.
    
    Next, we observe that $\bar f_{z}$ satisfies 
    \begin{align} \label{est:bar.f}
        \|\bar f_z\|_{H^k} \leq C
    \end{align}
    uniformly on compact sets $z \in K$.
    Indeed, by elliptic regularity corresponding to \eqref{est:elliptic.regularity}, it suffices to treat the case $k =0$. Observe that
    \begin{align}
        \bar f_z = \frac{\tilde f_z}{\int \tilde f_z \dd \xi}
    \end{align}
    for the normalized unique non-negative eigenvector $\tilde f_z \in \ker L_z$ with $\|\tilde f_z\|_2 = 1$.
    Thus, it suffices to show that
    \begin{align}
        \int \tilde f_z \dd \xi \geq c,
    \end{align}
       uniformly on compact sets $z \in K$. Again, we argue by contradiction. Indeed, if such an estimate was not true, we would find a sequence $z_n \to z_\ast$ such that $\tilde f_{z_n} \to f_\ast$ in $H^1(\S^2)$ with $\int f_\ast = 0$. But then, $f_\ast \in \ker L_{z_\ast}$ with $\|f\|_2 =1$ which contradicts \eqref{eq:ker.sign}. This proves \eqref{est:bar.f}.
       
    
    Finally, for $z_1,z_2 \in \R^m$ we have
    \begin{align}
        L_{z_1}( \bar f_{z_1} - \bar f_{z_2}) = \dv((h_{z_1} - h_{z_2})  \bar f_{z_2}). 
    \end{align}
    Using \eqref{est:bar.f}, the right-hand side is bounded by $C |z_1 - z_2|$ in $H^k$ and thus \eqref{est:regularity.L} implies Lipschitz estimates in $z$ as before and we conclude again by differentiation and iteration.
\end{proof}

\begin{proof}[Proof of Theorem \ref{thm:uniqueness.stationary}]
    We denote by $L_{t,x}$ the operator corresponding to $\bar h_{t,x}(\xi) = P_{\xi^\perp} h(t,x,\xi)$, and by $\bar f_{t,x}$ the corresponding unique solution to \eqref{eq:stationary.normalized}.

    We claim that $\bar \Psi_\psi(g) = 0$
    for all $\psi  \in C_c^\infty((0,T) \times \R^3 \times \S^2)$ implies 
    \begin{align} \label{eq:g=bar.f.mu}
        g = \frac{1}{\|\bar f_{\boldsymbol \cdot}\|_{L^2}^2}  \bar f_{\boldsymbol \cdot} \mu
    \end{align}
    for some $\mu \in  (C_c^\infty((0,T) \times \R^3))^\ast$ in the sense that
    \begin{align}
        \langle g, \psi \rangle = \langle \mu, \frac{1}{\|\bar f_{\boldsymbol \cdot}\|_{L^2}^2}  \int_{\S^2}   \bar f_{\boldsymbol \cdot}(\xi) \psi(\cdot,\xi) \dd \xi \rangle 
    \end{align}
    Note that this is well-defined since $(t,x,\xi) \mapsto \bar f_{t,x}(\xi)$ is smooth by Theorem \ref{th:regularity.f.bar} \ref{it:unique.stationary}.
    
    From this identity, $\int_{\S^2} \bar f_{t,x}(\xi) \dd \xi = 1$ and $\Theta_\theta(g)= 0$ for all
     $\theta \in C_c^\infty((0,T) \times \R^3)$, we immediately deduce 
     \begin{align}
         \mu = \|\bar f_{t,x}\|_{L^2}^2 \int_{\S^2} f_0(x,\xi) \dd \xi,
     \end{align}
     which yields uniqueness \am{and smoothness in $(t,\xi)\in(0,T)\times \S^2$ of the solution for almost all $x\in \R^3$.}
     
     It remains to prove \eqref{eq:g=bar.f.mu}.
     Let $\psi \in  C_c^\infty((0,T) \times \R^3 \times \S^2)$ and let 
     \begin{align}
         \varphi(t,x,\xi) = \psi(t,x,\xi) - \frac{1}{\|\bar f_{t,x}\|_{L^2}^2} \int \bar f_{t,x}(\zeta) \psi(t,x,\zeta) \dd \zeta \bar f_{t,x}(\xi).
     \end{align}
     Then, $\varphi(t,x,\cdot) \in (\ker L_{t,x})^\perp$. Hence, by Theorem \ref{th:regularity.f.bar} \ref{it.smooth.dependence.dual}, there exists $\tilde \varphi(t,x,\xi) \in  C_c^\infty((0,T) \times \R^3 \times \S^2)$ 
     such that $L_{t,x}^\ast \tilde \varphi(t,x,\xi) = \varphi(t,x,\xi)$ and $\int_{\S^2} \varphi(t,x,\xi) \dd \xi = 0$.
     Thus, $\langle g, \varphi \rangle = \Psi_{\tilde \varphi}(g) = 0$ and hence
     \begin{align}
         \langle g, \psi \rangle = \langle g,\frac{1}{\|\bar f_{\boldsymbol \cdot}\|_{L^2}^2}  \int \bar f_{\cdot}(\zeta) \psi(\cdot,\zeta) \dd \zeta \bar f_{\boldsymbol \cdot} \rangle
     \end{align}
     which yields the claim with
     \begin{align}
         \langle \mu, v \rangle = \langle g, \frac{1}{\|\bar f_{\boldsymbol \cdot}\|_{L^2}^2}  \bar f_{\boldsymbol \cdot} v \rangle.
     \end{align}
     This concludes the proof of the Theorem.
\end{proof}

\section{Proof of Lemma \ref{lem:Expectations}} \label{sec:Expectations}

\begin{proof}
{We drop the index $i$ in the proof. Let $b\in {C}^\infty_c((0,t), \R^3)$.
Recalling \eqref{eq:strato.xi_i.T_D}, we appeal to the Itô Stratonovitch conversion formula \eqref{eq:Conversion.composition} (after extending the occurring functions of $\xi$ to the whole space $\R^3$) to deduce with $\mathcal R = \sqrt{\mR_2}$
    \begin{align*}
       \E \left( \int_0^t b(s)\cdot  \sqrt{\mathcal{R}_2(\xi(s))} \circ d B(s) \right)&=  \E\left(\int_0^t  \mathcal{R}(\xi(s))b(s)  \circ d B(s) \right)\\
        & = \E\left(\frac{1}{2} \int_0^t tr\left(\nabla_\xi \mathcal{R}(\xi(s)) b(s) \sigma_D(\xi(s)) \right) ds\right)
    \end{align*}
     where we used that the expectation of Îto integrals vanish.
     From \eqref{eq:R_2}, we have
     \begin{align} \label{eq:nabla.R}
         \nabla_\xi \sqrt{\mathcal{R}_2(\xi) }b=(\sqrt{\gamma_{rot,\parallel}} - \sqrt{\gamma_{rot}}) [(\xi\cdot b) \Id +\xi \otimes b]
    \end{align}
    and we recall from \eqref{eq:sigma} that $\sigma_D = \sqrt 2 [\xi]_M$ (with the notation \eqref{eq:[T]_M}).
    In particular, since $\sigma_D$ is skew-symmetric,
    \begin{align}
        tr\left(\nabla_\xi \mathcal{R}(\xi) b \sigma_D(\xi) \right) &= (\sqrt{\gamma_{rot,\parallel}} - \sqrt{\gamma_{rot}}) tr\left((\xi \otimes b) \sigma_D(\xi) \right) \\
        &= \sqrt 2 (\sqrt{\gamma_{rot,\parallel}} - \sqrt{\gamma_{rot}}) ((\xi \times a_\alpha) \cdot b) (\xi \cdot a_\alpha)
    \end{align}
    for any  orthonormal basis $(a_\alpha)$ of $\R^3$. Since $\xi\times a_\alpha (\xi\cdot a_\alpha)=\xi \times \xi=0$, we conclude \eqref{eq:torque.average}.
    \medskip
    
   Let $A\in {C}^\infty_c((0,t), \Sym_0(3))$. 
   From \eqref{eq:R_2} and \eqref{def:mS_i} we deduce that \amnew{for any $T\in \R^3$
   \begin{align}
     A:  \left(\mS(\xi) \sqrt{\mR_2(\xi)}T\right) &=  \gamma_E \gamma_{rot}^{-\frac 1 2} A: \left(T\times\xi \otimes \xi \right)= \gamma_E \gamma_{rot}^{-\frac 1 2} \left([\xi]_M A \xi \right)\cdot T , \\
     \nabla_\xi \left([\xi]_M A \xi \right)&=[\xi]_M A-[A\xi]_M
   \end{align}
   where we used, for the first line, that for any $u,v,w\in\R^3$, 
   \begin{equation}\label{eq:cross_product}
   [u]_Mv=u\times v, \quad u\cdot (v\times w)=v \cdot(w\times u)
   \end{equation}
   }
   Thus, analogously as above
    \begin{align*}
       & \E \left(\int_0^t A(s) : \mathcal{S}(\xi(s)) \sqrt{\mathcal{R}_2(\xi(s))} \circ d B(s)\right) \\&= \frac{\gamma_E}{\sqrt{2\gamma_{rot}}} \int_0^t tr (([\xi(s)]_M A(s)-[A(s)\xi(s)]_M  ) [\xi(s)]_M ) ds.
    \end{align*}
 Using \eqref{eq:cross_product} and $[v]_M^\top=-[v]_M$ for any $v\in\R^3$, 
 we get for any orthonormal basis $(a_\alpha)$ of $\R^3$
 \begin{align*}
      tr  (([\xi]_M A-[A\xi]_M  ) [\xi]_M )&= 
     \left([\xi]_M a_\alpha\right) \cdot \left([A\xi]_M a_\alpha\right)-\left([\xi]_M a_\alpha\right)\cdot \left([\xi]_M a_\alpha\right)\\
     &=(\xi \times a_\alpha)\cdot [(A\xi) \times a_\alpha]-[A (\xi\times a_\alpha)]\cdot (\xi\times a_\alpha)\\
     &=A\xi \cdot [a_\alpha \times (\xi\times a_\alpha)]- A: (\xi\times a_\alpha)\otimes (\xi\times a_\alpha)\\
     &=A : \xi \otimes [a_\alpha \times (\xi\times a_\alpha)]- A: (\xi\times a_\alpha)\otimes (\xi\times a_\alpha)\\
     &=A : (2\xi\otimes \xi)+ A : (\xi \otimes  \xi - \Id)
 \end{align*}
 For the last line, we choose $a_1=\xi$ to get $a_\alpha \times (\xi\times a_\alpha)=\xi$ for $\alpha=2,3$ and 
     \begin{align}
         \sum_\alpha (a_\alpha \times \xi) \otimes (\xi \times a_\alpha)  &= \xi \otimes \xi -\Id .
    \end{align}
 which yields the desired result.
}
\end{proof}

\section{Stochastic calculus}\label{appendix_stochastics}
In this appendix we gather some theory on stochastic integrals that we rely on. In Subsection \ref{sec:Reminder.finite.dimension}, we collect standard  definitions and results on stochastic calculus in Euclidean space. We mostly follow \cite{Kuo} and keep the setting as simple as possible. For further results, one could for instance refer to \cite{karatzas2012brownian}.

At the end of Subsection \ref{sec:Reminder.finite.dimension}, we also precise the notion of solutions to the SDEs for the particle orientations $\xi_i$
in \eqref{eq:Particles.T_D} and \eqref{eq:Particles.T_u} which are SDEs on the manifold $\S^2$.
Instead of dealing with SDEs on manifolds in
an abstract way (see for instance \cite{Hsu}),
we follow a more pedestrian approach.
Namely, we just consider corresponding SDEs in the ambient space $\R^3$ and show a posteriori that the solutions stay on $\S^2$.

\medskip

Finally, Subsection \ref{sec:infinite} is devoted to stochastic calculus in infinite dimensions.
Stochastic calculus in infinite dimensional Hilbert spaces is quite well developed, see for instance the monograph \cite{da2014stochastic} where stochastic integrals of the form $\int f(s) \dd B(s)$
are considered for Hilbert-valued Brownian motions $B(s) \in E$ and  functions $f$
such that $f(s)\colon E \to H$ is a Hilbert-Schmidt operator to another Hilbert space $H$.
Corresponding stochastic integrals in a Stratonvich sense are used for instance in
\cite{MaurelliModinSchmeding19}.

For our purpose, it is sufficient to consider the case where $E = \R^d$ and  $H$ is an infinite dimensional separable Hilbert space. This case is considerably easier than the general
framework. Therefore, we choose to present here the necessary theory in a self-contained way based only on the finite dimensional stochastic calculus in the preceding section.

\subsection{Reminder on stochastic calculus in finite dimensional space} \label{sec:Reminder.finite.dimension}
In what follows we consider a probability space $(\Omega,\mathcal{A}, P)$ and recall that a continuous-time stochastic process is a map $t \mapsto X(t)$, $t\in \mathbb{R}^+$ such that $X(t)$ is a random variable. 
We recall first the definition of a Brownian motion
\begin{defi}[Brownian Motion]
 A continuous-time stochastic process $B : \R^+ \times \Omega \to \mathbb{R}$, $t\geq 0$ is called a Brownian motion starting in $a\in \mathbb{R}$ if and only if
\begin{enumerate}[label=(\roman*).,ref=(\roman*)]
\item  $B(0,\omega) = a$ for each $\omega \in \Omega$.
\item For any partition $0 \leq t_0 <t_1 <\cdots<t_n$, the increments $B(t_{i+1})-B(t_i)$ are independent random variables  with distribution $B(t_{i+1})-B(t_i) \sim N(0,t_{i+1}-t_i)$
\item $P$-almost every sample path $t \mapsto B(t,\omega)$ is continuous.
\end{enumerate}
An $\mathbb{R}^d$-valued stochastic process $B(t,\omega)=(B_1(t,\omega),\cdots,B_d(t,\omega))$ is called a multi-dimensional Brownian motion if and only if the components are independent one-dimensional Brownian motions.

\end{defi}
In order to recall the definition of  stochastic integrals, we introduce first filtrations and the notion of adapted processes.
\begin{defi}[Filtrations and adapted processes]
Let $(\Omega,\mathcal{A}, P)$ be a probability space.
\begin{itemize}
\item A continuous-time filtration on $(\Omega, \mathcal{A})$ is a family $ (\mathcal{F}_t)_{t \geq 0}$, indexed by time, of increasing $\sigma$-algebras $\mathcal{F}_t \subseteq \mathcal{A}$, that is, satisfying $\mathcal{F}_t \subset \mathcal{F}_s$ for $t \leq s$.

\item A real-valued stochastic process $M$ on $(\Omega, \mathcal{A}, P)$ is adapted with respect to $(\mathcal{F})_{t \geq 0}$ if $M(t)$ is $\mathcal{F}_t$ measurable for all $t \geq 0$.

\item The natural (canonical) filtration $(\mathcal{F})_{s \geq 0}$ generated by a stochastic process $M$, denoted by $\sigma(M)$, is the filtration formed, for each $t$, by the smallest $\sigma$-algebra for which the maps $\omega \mapsto M(s,\omega)$, are measurable functions for all $ 0 \leq s \leq t$ .
\end{itemize}
\end{defi}
\begin{defi}[Martingale]
A stochastic process $M$ adapted to a filtration $(\mathcal{F}_t)_{t \geq 0}$ is a martingale if the following conditions hold:
\begin{itemize}
\item $E[M(t)] < + \infty$ for all $t \geq 0$.
\item $E[M(t) \mathcal{F}_s] = M(s)$ for all $ 0 \leq s \leq t $
\end{itemize}
\end{defi}

\subsubsection*{Itô integral with respect to Brownian motion and Itô formula}
We consider below a Brownian Motion and a filtration $(\mathcal{F}_t)_{ t \geq 0}$ such that
\begin{itemize}
\item For all $t \geq 0$, $B(t)$ is $\mathcal{F}_t$ measurable.
\item For all $0 \leq s \leq t$, the random variable $B(t)-B(s)$ is independent of the $\sigma$-algebra $\mathcal{F}_s$.
\end{itemize}

We recall the construction of the Itô integral in the same way as  the Riemann integral, proceeding first for step stochastic processes given by
$$
f(\omega)=\underset{i=0}{\overset{n}{\sum}}A_i(\omega) 1_{]s_i,s_{i+1}]},
$$
for a given decomposition $s_0=0<s_1<\cdots<s_n=T $ of $[0,T]$ such that $ A_i$ is $F_{s_i}$-measurable and $E[A_i^2]<+\infty$. We set 
$$I(f)=\underset{i=0}{\overset{n}{\sum}}A_i (B({s_{i+1}}) -B({s_{i}})).
$$
This defines a linear mapping $I$ and moreover we have $E I(f)=0$ and
\begin{equation}\label{isometry_ito}
E( I(f)^2)=\int_0^T E(|f(t)|^2) dt.
\end{equation}
for each step stochastic process, see \cite[Lemma  4.3.2]{Kuo}. Next, the idea is to extend this definition for each stochastic process in the space $f\in L^2_{\text{ad}}([0,T]\times \Omega)$ defined as the space of stochastic processes $f$ satisfying
\begin{itemize}
\item $ \int_0^T E |f(s)|^2 ds < +\infty $ 
\item $f$ adapted to the filtration $(\mathcal{F}_s)_{0 \leq s \leq T}$.
\end{itemize}
Indeed, each stochastic process $f\in L^2_{\text{ad}}([0,T]\times \Omega)$ can be approximated by a sequence of step stochastic processes $(f_k)_k$ in the following sense
$$
\underset{k \to \infty}{\lim} \int_0^t E | f(s)-f_k(s)|^2 dt = 0,
$$
see \cite[Lemma  4.3.3]{Kuo}. According to this, one can define the Ito integral as the limit 
\begin{equation}\label{def_ito}
I(f)= \underset{k \to \infty}{\lim} I(f_k)
\end{equation}
where the limit does not depend on the choice of the sequence $(f_k)_k$. In particular, thanks to \eqref{isometry_ito}, $ I \colon L^2_{\text{ad}}([0,T]\times \Omega) \to L^2(\Omega)$ is an isometry. We gather below additionnal properties, see \cite[Theorem 4.3.5, Theorem 4.6.1 and Theorem 4.6.2]{Kuo}.
\begin{thm}[Itô Isometry] \label{th:Ito}
Let $f \in L^2_{\text{ad}}([0,T]\times \Omega)$. The stochastic process corresponding to the Itô integral $I(f)$ defined in \eqref{def_ito} and denoted by
$$
X(t)= \int_0^t f(s) dB(s), \: 0 \leq t \leq T,
$$
is a continuous martingale with respect to the filtration $ (\mathcal{F}_t)_{0 \leq t \leq T}$ such that $E[X(t)]=0$ and
$$
E \left| \int_0^t f(s) dB(s) \right|^2=\int_0^t E|f(s)|^2 ds.
$$
\end{thm}

\begin{rem} \label{rem:multiIto}
    The definition of Itô integrals and the Itô isometry extends straightforwardly to the multidimensional case: For $d$ independent Brownian motions $B_i$, $1 \leq i \leq d$ and $ \sigma\in {L}^2_{\text{ad}}([0,t] \times \Omega 
    )$ a $d\times d$ valued function,
    \begin{align} \label{eq:mulit.Ito.isometry}
       \| \int_0^t \sigma(s) dB(s)\|^2_{L^2(\Omega)} = \int_0^t \|\sigma(s)\|_{HS}^2,
    \end{align}
    where $\|A\|_{HS} = \sqrt{\Tr(A^T A)}$ denotes the Hilbert-Schmidt norm of a matrix $A \in \R^{d \times d}$.
\end{rem}

Following \cite[Subsection 7.5]{Kuo} we recall the multidimensional Itô formula.
\begin{defi}[Ito process]
Let $ \sigma\in {L}^2_{\text{ad}}([0,T] \times \Omega)$ a $d\times d$ valued function, $b\in{L}^2_{\text{ad}}([0,T] \times \Omega)$ a $\mathbb{R}^d$ valued function.
Then, the stochastic process
\begin{equation}\label{def_ito_process}
X(t)=X_0 + \int_0^t \sigma(s) dB(s) + \int_0^t b(s) ds 
\end{equation}
is called Itô process with drift $b$ and diffusion $\sigma$.
\end{defi}

\begin{lem}[Itô formula] \label{lem:Ito.formula}

Let $\phi \in {C}(\mathbb{R}\times\mathbb{R}^d ; \mathbb{R}^d)$ with continuous partial derivatives $\partial_t\phi $,$\nabla_x \phi $, $\nabla_x^2\phi$. 
Then $\phi(t,X(t))$ is also an Itô process satisfying
\begin{multline}\label{eq:Ito_formula}
\phi(t,X(t))=\phi(0,X(0))+ \int_0^t \nabla \phi(s,X(s))  \sigma(s) dB(s)\\
+ \int_0^t \left[\partial_t \phi(s,X(s)) + b(s) \cdot \nabla_x\phi(s,X(s))  + \frac{1}{2} \sigma(s)\sigma(s)^\top :\nabla_x^2 \phi(s,X(s)) \right ]ds.
\end{multline}
\end{lem}

\subsubsection*{Stratonovich Integral}

We now turn to the definition of the Stratonovich integral.
We first recall the definition of quadratic variance and covariance.
\begin{defi}[Quadratic variation and covariation]
The quadratic variation of a real-valued stochastic process $M$ is the process written as $[M]$ and defined as the following limit in probability, if it exists
$$
[M](t) = \underset{\|P\|\to 0}{\lim} \underset{i}{\overset{n}{\sum}} (M(t_i)-M(t_{i-1}))^2,
$$
where $P=\{t_0, t_1,\cdots, t_n \}$ is any partition of $[0,t]$ and $\|P\|=\underset{1 \leq i \leq n}{\max}|t_i-t_{i-1}|$. Moreover, for two real-valued stochastic processes $X, Y$, the covariation, denoted by $[X,Y]$, is defined as
$$
[X,Y](t) = \underset{\|P\| \to 0 }{\lim} \underset{i}{\overset{n}{\sum}} (X({t_i})-X({t_{i-1}}))(Y({t_i})-Y({t_{i-1}})).
$$
\end{defi}
\begin{rem}\label{rem_quadratic_variation}
\begin{enumerate}[label=(\roman*).,ref=(\roman*)]
\item  The quadratic variation exists for all continuous finite variation processes, and is zero.
\item The quadratic variation exists for all right continuous, square integrable martingale with left-hand limits, see \cite[Formula (6.4.4)]{Kuo}
\item If $B$ and $W$ are two independent Brownian motions, $f\in L^2_{\text{ad}}([0,t]\times \Omega)$ and $X(t)= \int f(s) dB(s)$ then $[X,W](t)=0$ and $[X,B](t)=\int f(s) ds$.
\end{enumerate}

\end{rem}
\begin{defi}[Stratonovich integral] \label{def:Strato}
Let $X \in L^2_{\text{ad}}([0,T]\times \Omega)$ such that $[X]$ exists. 
The Stratonovich integral is defined as
\begin{equation}\label{eq:Ito_vs_Stratonovich}
\int_0^t X(s) \circ dB(s) = \int_0^t X(s) dB(s) +\frac{1}{2} [X,B](t).
\end{equation}
\end{defi}
\begin{rem}\label{rem:Strato}

\begin{enumerate}[label=(\roman*).,ref=(\roman*)]
\item \label{it:Strato.Ito.process} In particular, if $X$ is  given as an Ito process of the form \eqref{def_ito_process}, independence of $B_i$, $B_j$ for $i \neq j$  and Remark \ref{rem_quadratic_variation} iii) yield
\begin{align*}
\int_0^t X(s) \circ dB(s) &= \int_0^t X(s) dB(s) +\frac{1}{2} \underset{i}{\sum} \left[\int_0^t \sigma_{ij}(s) dB_j,B_i \right ](t)\\
&= \int_0^t X(s) dB(s) +\frac{1}{2} \int_0^t \text{tr} \sigma(s) ds.
\end{align*}
\item  \label{it:Conversion.composition} One can extend the above formula to all $g \in {C}^1([0,T] \times \mathbb{R}^d; \mathbb{R}^d)$ in the following way
\begin{align} \label{eq:Conversion.composition}
\int_0^t g(s,X(s)) \circ dB(s) &= \int_0^t g(s,X(s)) dB(s)  +\frac{1}{2} \int_0^t \text{tr}(\nabla g  \sigma ) ds.
\end{align}
where we used Itô formula \eqref{eq:Ito_formula} in the case $g\in {C}^2$ and a density argument in the case ${g\in {C}^1}$.

\item  \label{it:stratonovich.chain.rule} The chain rule holds for Stratonovitch integrals. 
More precisely, if \rhnew{$\sigma \in L^2_{\text{ad}}([0,T]\times \Omega;\R^{d\times d})$} and
$Y(t) := \int_0^t \sigma(s) \circ \dd B(s) $ and
$\phi : [0,T] \times \R^d \to \R$ is a smooth  function, then
\begin{align}
    \phi(T,Y(T)) = \phi(0,Y(0)) + \int_0^T \partial_t \phi(t,Y(t)) \dd t + \int_0^T \nabla \phi(t,Y(t)) \rhnew{\sigma(t)} \circ \dd B(t).
\end{align}
\end{enumerate}
\end{rem}

\subsubsection*{Stochastic differential equations}\label{section_SDE_existence}
Following \cite[Subsection 10.3]{Kuo}, we summarize existence and uniqueness result for SDEs of the form 
$$
dX(s)=\sigma(s,X(s)) dB(s)+b(s,X(s))ds, \qquad 0\leq s \leq T, \: X(0)= X^0, 
$$
which must be interpreted as the stochastic integral equation 
\begin{equation}\label{eq:SDE}
X(t)=X^0 +\int_0^t \sigma(s,X(s)) dB(s)+\int_0^t b(s,X(s))ds.
\end{equation}
\begin{defi} \label{def:SDE}
We say that a jointly measurable stochastic process $X$ on $[0, T]$, is a solution of the SDE \eqref{eq:SDE} if 
\begin{enumerate}[label=(\roman*).,ref=(\roman*)]
\item The stochastic process $\sigma(s,X(s)) $ lies in $L^2_{\text{ad}}([0,T]\times \Omega)$ and $\displaystyle \int_0^t \sigma(s,X(s)) dB(s)$ is an Itô integral for each $t\in [0,T]$.
\label{it:sigma}
\item Almost all sample paths of the stochastic process $ b(s,X(s))$ belong to $L^1(0,T)$.
\item For each $t\in [0,T]$, equation \eqref{eq:SDE} holds true almost surely. \label{it:SDE}
\end{enumerate}
\end{defi}
A global existence and uniqueness result can be obtained, in a similar way as for ODEs, in the case where both $\sigma(t,\cdot)$ and $b(t,\cdot)$ satisfy a global Lipschitz property on $\mathbb{R}^d$.

\begin{thm}\label{thm_SDE_existence_uniqueness}
Let $(t,x) \mapsto (\sigma_{ij})_{1\leq i,j\leq d}$, $(t,x) \mapsto (b_i)_{1\leq i \leq d}$ measurable on $[0,T] \times \mathbb{R}^d$ satisfying a global Lipschitz property on $\mathbb{R}^d$ uniformly on $[0,T]$. Assume that $X^0$ is an $\mathcal{F}_0$ measurable random variable with $\mathbb{E}[|X^0|^2] <\infty$. Then the SDE \eqref{eq:SDE} has a unique continuous solution $X\in L^2_{\text{ad}}([0,T] \times \Omega)$.
\end{thm}

We end this subsection by making the following observations regarding SDEs on a smooth compact submanifold $M \subset \R^d$. In this case, it is convenient to consider Stratonovitch SDEs which read in integral form
\begin{equation}\label{eq:SDE.M}
X(t)=X^0 +\int_0^t \sigma(s,X(s)) \circ dB(s)+\int_0^t b(s,X(s))dt,\: t \in[0,T].
\end{equation}
for $(\sigma,b) \colon [0,T] \times M \times \to (\R^{d \times d},\R^d)$

We first remark how such a Stratonovitch SDE may be transformed to an Ito SDE:
\begin{rem} \label{rem:SDE.Strato.Ito}
    If $\sigma \in C^1([0,T] \times M)$ and $X$ solves the Ito SDE 
    \begin{align}
X(t)=X^0 +\int_0^t \sigma(s,X(s)) dB(s)+\int_0^t b(s,X(s)) + \frac 1 2 \sigma(s,X(s)) : \nabla \sigma(s,X(s)) ds,
    \end{align}
    where $(\sigma : \nabla \sigma)_i = \sum_{k} \sigma_{k} \cdot \nabla \sigma_{ki}$,
    then, $X$ satisfies \eqref{eq:SDE.M} almost surely.
    
    Indeed, this follows immediately from \eqref{eq:Conversion.composition} with $g = \sigma$.
\end{rem}

\begin{thm} \label{thm:SDE.M}
    Assume that $X^0 \colon \Omega \to M$ is $\mathcal F_0$ measurable, $\sigma \in C^2([0,T] \times M;\R^{d \times d})$ and $b \in C^1([0,T] \times M ;\R^d)$ satisfy $b(s,x) \in T_x M$ and $\sigma(s,x) v \in T_x M$ for all  $(s,x) \in [0,T] \times M$
    and for all $v \in \R^d$. Then, there exists a unique solution $X$ to \eqref{eq:SDE.M}.
\end{thm}
\begin{proof} For simplicity we restrict to the case $M = \S^{d-1}$. 

    Since $M \subset \R^d$ is compact and smooth, there exist extensions $\tilde \sigma \in C_c^2([0,t] \times \R^d)$
    and $\tilde b \in C_c^1([0,T] \times \R^d)$ such that $(\tilde \sigma, \tilde b)|_M = (\sigma,b)$. 
    By Theorem \ref{thm_SDE_existence_uniqueness} (and the remark above), there exists a unique solution to the SDE
    \eqref{eq:SDE.M} with $M$ replaced by $\R^d$ and $\sigma$ replaced by  $\tilde \sigma$.

    By the Stratonovitch chain rule, Remark \ref{rem:Strato} 
    \ref{it:stratonovich.chain.rule}, and the assumptions on $\sigma$ and $b$  we have
    \begin{align}
        \dd |X|^2 = 0,
    \end{align}
    which ensures that $X$ takes values in $M=\S^d$. Thus, $X$ is a solution to \eqref{eq:SDE.M}.
    
    Conversely, every solution to \eqref{eq:SDE.M} also satisfies the Ito SDE with $\tilde \sigma$ on $\R^d$ for which we already know uniqueness.
\end{proof}

\subsection{Extension to infinite dimensions}
\label{sec:infinite}

Let $H$ an infinite dimensional separable Hilbert space, and let $(e_k)_k$ be an orthonormal basis of $H$. Let
$(\Omega, \mathcal{A},P)$ a probability space and $(\mathcal{F}_t)_t$ a filtration as in the previous subsection.
Let us denote by $L^2_{\text{ad}}([0,T]\times \Omega;\mathcal{L}(\mathbb{R}^d,H))$ the space of operator valued adapted processes $X:[0,T]\times \Omega \to \mathcal{L}(\mathbb{R}^d,H)$ such that
$$
\|X\|_{L^2_{\text{ad}}([0,T]\times \Omega;\mathcal{L}(\mathbb{R}^d,H))}^2:=\mathbb{E} \int_0^T \| X(t)\|_{HS}^2 dt < \infty
$$
where 
\begin{align}
    \| X\|^2_{HS}=\text{tr} (XX^\star)=\underset{k}{\sum}\| X^\star e_k\|_{\R^d}^2 = \sum_l \| X a_l\|_{H}^2 
\end{align}

for the standard basis $(a_l)_l$ of $\R^d$.
Note that we can identify $\L(\R^d,H) = H^d$ through $A b = \sum_i A_i b_i$ for $A \in \L(\R^d,H)$, $b \in \R^d$, where $A_i = A a_i$. Then $A^\star e_k = \sum_i (A_i,e_k)a_i$.

We extend the definition of the Itô integral in the following
way
\begin{defi}
	Let $B$ be a $d$-dimensional Brownian motion and $X \in L^2_{\text{ad}}([0,T]\times \Omega;\mathcal{L}(\mathbb{R}^d,H))$. Then, for $0 \leq t \leq T$
	\begin{align} \label{def:Ito.infinite}
		\int_0^t X(s) \dd B(s) := \sum_{k \in \N} \sum_{i=1}^d \left(\int_0^t (X_i(s),e_k) \dd B_i(s) \right)e_k.
	\end{align}
\end{defi}
By density, we have that this is well defined and the Îto isometry carries over.
\begin{prop}  \label{prop:Ito.isometry.infinite}
	Let $B$ and $X$ be as above. Then $\int_0^t X(s) \dd B(s) \in L^2(\Omega;H)$
	and for all $0 \leq t \leq T$
	\begin{align} \label{eq:Ito.isometry.infinite}
		\E\left\|\int_0^t X(s) \dd B(s) \right\|^2_H = \int_0^t \E\|X(s)\|_{HS}^2 \dd s.
	\end{align}
\end{prop}
\begin{proof}
	We observe that 
	\begin{align}
	\int_0^t X(s) \dd B(s) &= \sum_k \lambda_k e_k, \\
	\lambda_k &= \sum_i\int_0^t (X_i(s),e_k) \dd B_i(s).
	\end{align}
	Then, by {\cite[Theorem 2.3.4]{Kuo}}
	\begin{align}
		\E \lambda_k^2  = \sum_{i=1}^d \int_0^t\E|(X_i(s),e_k)|^2 \dd s =  \int_0^t\E|(X(s),e_k)|^2 \dd s,
	\end{align}
and thus
\begin{align}
	\E \sum_k \lambda_k^2 = \int_0^t \E\|X(s)\|_{HS}^2 \dd s.
\end{align}
Therefore, by Parseval's identity and dominated convergence theorem, the integral \eqref{def:Ito.infinite} is well-defined and \eqref{eq:Ito.isometry.infinite} holds.
\end{proof}

Similarly, we can define the Stratonovitch integral as follows.
\begin{defi} \label{defi:Strato.infinite}
	Let $B$ be a $d$-dimensional Brownian motion and $X \in L^2_{\text{ad}}([0,t]\times \Omega;\mathcal{L}(\mathbb{R}^d,H))$. Moreover, assume that
	\begin{align}
		\sum_{k \in \N}  \sum_{i=1}^d \E \left[[(X_i,e_k),B_i]^2(T)\right] < \infty.
	\end{align}
	
	Then, we define the Stratonovitch integral for $ 0\leq t \leq T$
	\begin{align} \label{def:Strato.infinite}
		\int_0^t X(s) \circ \dd B(s) := \sum_{k \in \N} \sum_{i=1}^d \left(\int_0^t (X_i(s),e_k) \dd B_i(s)  + [(X_i,e_k),B_i](t) \right)e_k.
\end{align}
\end{defi}
\begin{rem}
    Clearly Proposition \ref{prop:Ito.isometry.infinite}
    implies that under the assumptions of this definition, the Stratonovitch integral is well-defined and satisfies
	\begin{align}
		\E\left\|\int_0^t X(s) \circ \dd B(s) \right\|^2_H \leq 2 \int_0^t \E\|X(s)\|_{HS}^2 \dd s + \frac 1 2  \sum_{k \in \N}   \E \left[ \left( \sum_{i=1}^d [(X_i,e_k),B_i](t)\right)^2 \right] .
	\end{align}
\end{rem}

We are especially concerned with Stratonovitch integrals when the integrand is of the form $X(s) = g(Y(s))$, where $Y$ is an Itô process.

\begin{prop} \label{pro:Strato.infinite}
Let $0 < t \leq T$, let $B$ be a $d$-dimensional Brownian motion, and let $g \in C^1(\R^d ; \L(\R^d,H))$ with $\|g\|_{C^1(\R^d ;\L(\R^d,H))} < \infty$, $\sigma  \in L^2_{\text{ad}}([0,T]\times \Omega;\R^{d \times d})$ and  $b  \in L^2_{\text{ad}}([0,T]\times \Omega;\R^{d}))$.
Let $Y \colon \Omega \times [0,T] \to \R^d$ be the Itô process
\begin{align} \label{eq:Y.Ito}
	Y(t) = Y(0) + \int_0^t \sigma(s) \dd B(s) + \int_0^t b(s) \dd s.
\end{align}
Then $g \circ Y$ satisfies the assumptions of Definition \ref{defi:Strato.infinite} and
 	\begin{align} \label{eq:Strato.composition.infinite}
		Z(t):= \int_0^t g(Y(s)) \circ \dd B(s) = \int_0^t g(Y(s)) dB(s)+ \sum_i \sum_j \frac{1}{2}\int_0^t  \partial_i g_j(Y(s))\sigma_{ij} \dd s.
\end{align}

Moreover, $Z \in H^{\frac 1 2_-}((0,T) \times H)$ and
\begin{align} \label{est:Z.H^s}
   \E\left[ \|Z\|_{H^{\frac 1 2_-}((0,T) \times H)}^2  \right]^{\frac 1 2 } \leq C_T \left(\| g \|_{L^\infty(\R^d; H^d)} + \|  \sigma\|_{L^\infty(\Omega \times [0,T];\R^{d \times d}_{HS})} \|\nabla  g\|_{L^\infty( \R^d;H^{d \times d}_{HS})} \right),
\end{align}
for a constant $C_T$ that depends only on $T$.
\end{prop}

\begin{proof}
	Clearly $g \circ Y \in L^2_{\text{ad}}([0,T]\times \Omega;\mathcal{L}(\mathbb{R}^d,H))$.
	Moreover, for all $1 \leq i \leq d$ and all $k \in \N$ Definition \ref{def:Strato} and
	\eqref{eq:Conversion.composition} imply
	\begin{align}
		\sum_{i=1}^d	[(g_i \circ Y,e_k),B_i](t) = \sum_{i=1}^d \sum_{j=1}^d \int_0^t  (\partial_j  g_i(Y(s)),e_k) \sigma_{ij} \dd s.
	\end{align}
	In particular,
	\begin{align} \label{est:Z.H^s.1}
		\sum_{k \in \N}  \E \left[ \sum_i [(g \circ Y,e_k),B_i]^2(t) \right] 
		&\leq t^2 \left\| \sum_{i,j = 1}^d  \sigma_{ij} \partial_i  g_j \circ Y  \right\|^2_{L^\infty(\Omega \times (0,t);H)} \\
		&\leq t^2 \|\nabla g\|^2_{L^\infty((0,t);H^{d\times d})} \|\sigma\|^2_{L^\infty(\Omega \times (0,t))}.
	\end{align}
	Thus, $Z$ is well-defined and satisfies \eqref{eq:Strato.composition.infinite}.
	
	For $s \in (0,1)$, recall
	\begin{align}
	    \|f\|_{H^s(0,T);H)} = \|f\|_{L^2(0,T);H)} + \left(\int_0^T \int_0^T \frac{\|f(t_1) - f(t_2)\|^2_H}{|t_1-t_2|^{1+2s}} \right)^{\frac 1 2}.
	\end{align}
	Thus, \eqref{est:Z.H^s} follows from \eqref{eq:Ito.isometry.infinite} and \eqref{est:Z.H^s.1} (adapted to $(t_1,t_2)$).
\end{proof}

\begin{rem} \label{rem:manifold}
This proposition can be directly extended to the case, when $M \subset \R^d$ is a smooth submanifold, $g \in C^1(M; \L(\R^d,H))$ with $\|g\|_{C^1(M;\L(\R^d,H))} < \infty$ and $Y \colon \Omega \times [0,t] \to M$ is the Îto process as above with the constraint $\sigma (s) v \in T_{Y(s)}M$ for all $v \in \R^d$ (which is in fact a necessary condition for $Y$ to stay in $M$).

Indeed, we can find $\tilde g \in C^1(\R^d; \L(\R^d,H))$ such that $g = \tilde g$ on $M$ and apply the above proposition. Due to the condition $\sigma (s) v \in T_{Y(s)}M$, we have $\sum_i \sum_j \partial_i \tilde g_j \sigma_{i,j} = \sum_j D_{\sigma_j} g_j $, where $D_{\sigma_j}$ is the derivative in the direction $\sigma_j = (\sigma_{ij})_i \in T_{Y}M$. Therefore \eqref{eq:Strato.composition.infinite} becomes
 	\begin{align} \label{eq:Strato.composition.infinite.manifold} 
		Z(t) := \int_0^t g(Y(s)) \circ \dd B(s) = \int_0^t g(Y(s)) dB(s)+  \sum_j \frac{1}{2}\int_0^t   D_{\sigma_j} g_j(Y(s)) \dd s,
\end{align}
and \eqref{est:Z.H^s} still holds true.
\end{rem}

\section{Some embeddings in weighted Sobolev spaces}\label{sec:weighted sobolev space}
In this appendix we show some embeddings  for the  weighted fractional Sobolev space introduced in Subsection \ref{subsec:notatation}. 

\begin{lem} \label{lem:embedding.Lebesque.weighted}
Let $K$ be an open bounded set in $\mathbb{R}^3$. Consider the non negative function $w(x)=(1+|x|)^a$ with $a\geq 0$. Then, for $p\geq\frac{6}{3+2s}$, \rh{we have the continuous embedding}
\[L^{p,2}_{w,K} \hookrightarrow H^{-s}_w(\R^3).\]
\end{lem}
\begin{proof}
   By recalling the definition of the $H^s(K)$ norm and $H^s_w(K)$ norm and noting that $w$ is bounded below on $K$, we have
    \begin{equation}
    \norm{f}_{H^s(K)}\leq \norm{f}_{H^s_{1/w}(\R^3)}.
    \end{equation}
By Sobolev embedding we have that $H^s(K)\hookrightarrow L^{p'}(K)$ since $p' \leq \frac{6}{3-2s}$. Hence, 
\begin{equation}\label{eq:embedding_1}
H^s_{1/w}(\R^3)\hookrightarrow H^s(K)\hookrightarrow L^{p'}(K)
\end{equation}
Moreover, since $1/w\leq 1$, it is straightforward to prove that
\begin{equation}\label{eq:embedding_2}
    H^s_{1/w}(\R^3)\hookrightarrow L^2_{1/w}(\R^3\setminus K).
\end{equation}

Now we observe that we can identify ${L^{p,2}_{w,K}} = L^2_{w}(\R^3\setminus K)\times L^{p}(K)$ and thus the dual $({L^{p,2}_{w,K}})^\ast = L^2_{w}(\R^3\setminus K)\times L^{p}(K)$.
Combining \eqref{eq:embedding_1} and \eqref{eq:embedding_2} we get that
\[ H^s_{1/w}(\R^3)\hookrightarrow L^2_{1/w}(\R\setminus K)\times L^{p'}(K).\]
 Since $V \hookrightarrow W$ implies $W^\ast \hookrightarrow V^\ast$ and $\left(L^2_{1/w}(\R^3\setminus K)\times L^{p'}(K)\right)^*=L^2_{w}(\R^3\setminus K)\times L^{p}(K)$
we conclude the desired embedding.
\end{proof}
\begin{lem} \label{lem:compact.weighted}
For \rh{$z>s>\frac{1}{2}$},  the embedding $H^{-s}_w(\R^3)\hookrightarrow H^{-z}(\R^3) $ is compact.
\end{lem}
\begin{proof}
\rh{By Schauder's Theorem, it suffices to show compactness of the embedding $H^z(\R^3)\hookrightarrow H^{s}_{1/w}(\R^3)$ }

To this end, we will prove that the unit ball of radius  in $H^{z}(\R^3)$ is precompact in $H^{s}_{1/w}(\R^3)$.  
Fixed $L>0$, we introduce $B(0,L)$, a ball of radius $L$ in $\mathbb{R}^3$
and a cut-off function $\eta:\R^3\to\R$ such that $\norm{\eta}_{C^1(\R^3)}\leq C$ and
\[\eta(x)=\begin{cases}
1\,\,\textrm{if}\quad x\in B(0,L),\\
0\,\,\textrm{if}\quad x\in B(0,L+2)^c.
\end{cases}\]
\rh{Let $\phi \in H^z(\R^3)$.} Then, 
\begin{align}\label{eq:norm_etaphi}
\norm{\phi \eta}_{H^z(\R^3)}&\leq C \norm{\phi}_{H^z(\mathbb{R}^3)}\norm{\eta}_{C^1(\mathbb{R}^3)}, \\
\label{eq:norm_etaphi2}
    \norm{\phi(1-\eta)}_{H^s_{\frac{1}{w}}(\R^3)}&\leq \frac{C}{(1+L)^{a/2}} \norm{\phi}_{H^s(\R^3)}\norm{(1-\eta)}_{C^1(\R^3)}
\end{align}
Indeed,
\begin{align*}
    [\phi \eta]_{H^z(\R^3)}^2 &= \int_{\R^3}\int_{\R^3} \frac{|(\phi\eta)(x)-(\phi\eta)(y)|^2}{|x-y|^{3+2z}}dxdy \\
    &\leq 2 \norm{\eta}^2_{C(\mathbb{R}^3)} [\phi]^2_{H^z(\R^3)}+2\int_{\R^3}\int_{\R^3}\frac{|\phi(y)|^2|\eta(x)-\eta(y)|^2}{|x-y|^{3+2z}}dxdy,
\end{align*}
and the second term is further estimated as
\begin{align*}
    \int_{\R^3}\int_{\R^3}\frac{|\phi(y)|^2|\eta(x)-\eta(y)|^2}{|x-y|^{3+2z}}dxdy 
&\leq\norm{\eta}^2_{C^1(\R)^3}\int_{\R^3}\int_{|x-y|\leq 1}\frac{|\phi(y)|^2}{|x-y|^{1+2z}}dxdy\\
&+\norm{\eta}^2_{C^1(\R^3)} \int_{\R^3}\int_{|x-y|> 1}\frac{|\phi(y)|^2}{|x-y|^{3+2z}}dxdy\\
&\leq C \norm{\eta}^2_{C^1(\R)^3}\norm{\phi}^2_{L^2(\R^3)}.
\end{align*}
Moreover,
\begin{align*}
    \|\phi (1-\eta)\|_{H^s_{1/w}(\R^3)}^2 &=\int_{B(0,L)^c} \frac{|\phi(x)(1-\eta(x))|^2}{w(x)} \\
    &+\int_{B(0,L)^c}\int_{B(0,L)^c}\frac{|(\phi(1-\eta))(x)-(\phi(1-\eta))(y)|^2}{|x-y|^{3+2s} w(x) w(y)}dxdy\\
    &\leq \int_{\R^3} |\phi(x)(1-\eta(x))|^2\frac{1}{(1+L)^a}dx\\
    &+\frac{1}{(1+L)^{2a}}\int_{\R^3}\int_{\R^3}\frac{|(\phi(1-\eta))(x)-(\phi(1-\eta))(y)|^2}{|x-y|^{3+2s}}dxdy\\
    &\leq \norm{\phi(1-\eta)}^2_{H^s(\R^3)}\frac{1}{(1+L)^{a}}.
\end{align*}
\rh{
Let $\eps >0$. Since $H^{z}_0(B(0,{L+2}))$ is compactly embedded into $H^{s}_0(B(0,{L+2}))$, there exists $N \in \N$ and functions $u_i \in H^{s}_0(B(0,{L+2}))$, $1 \leq i \leq N$, such that
\begin{align}
    B_{H^{z}_0(B(0,{L+2}))}(0,C) \subset \bigcup_{i=1}^N B_{H^{s}_0(B(0,{L+2}))}(\eps,u_i),
\end{align}
where $C$ is the constant from \eqref{eq:norm_etaphi}.
Thus, (extending the functions $u_i$ by $0$ to functions in $H^{z}(\R^3)$ that vanish in $\R^3 \setminus B(0,L+2)$), for each $v \in H^{z}(\R^3)$ with $\|v\|_{H^{z}(\R^3)} \leq 1$, there exists $1 \leq i \leq N$ such that
\begin{align}
    \|v - u_i\|_{H^{s}_{1/w}(\R^3)} \leq \|\eta v- u_i\|_{H^{s}(\R^3)} + \|(1-\eta) v\|_{H^{s}_{1/w}(\R^3)} \leq \eps + \frac{C}{(1+L)^{a/2}}.
\end{align}

Choosing $L$ sufficiently large finishes the proof.
}
\end{proof}

\section*{Acknowledgements}
\rhnew{
We would like to thank Francesco de Vecchi for insights about Brownian motions on spheres.
We are grateful to Immanuel Zachhuber for pointing out references on stochastic integrals in infinite dimensional spaces.
Moreover, we thank Henrik Matthiesen for helpful discussions on some aspects related to global analysis.
}

R.H. is supported  by the German National Academy of Science Leopoldina, grant LPDS 2020-10.
Moreover, R.H. acknowledges support by the Agence Nationale de la Recherche,  Project BORDS, grant ANR-16-CE40-0027-01
and by the Deutsche Forschungsgemeinschaft (DFG, German Research Foundation) 
through the collaborative research center ``The Mathematics of Emerging Effects'' (CRC 1060, Projekt-ID 211504053) 
and the Hausdorff Center for Mathematics (GZ 2047/1, Projekt-ID 390685813). 
\amnew{M.L. is supported by by the Italian Ministry of Education, University and
Research (MIUR), in the framework of PRIN projects 2017FKHBA8 001. 
A.M. is supported by the SingFlows project, grant ANR-18-CE40-0027 of the French National Research Agency (ANR).}

\begin{refcontext}[sorting=nyt]
\printbibliography
\end{refcontext}
\end{document}